\documentclass[11pt,envcountsect,envcountsame]{llncs}
\usepackage{makeidx}
\usepackage{amsfonts}
\usepackage{amsmath}
\usepackage{amssymb}
\usepackage{amscd}
\usepackage[dvips]{graphicx}
\usepackage{algorithmic}
\usepackage{algorithm}
\usepackage{comment}
\usepackage{color}
\usepackage{array,amssymb,amsmath,amsbsy}
\usepackage{enumerate}
\usepackage[bottom]{footmisc}
\usepackage[all,cmtip]{xy}
\usepackage{fullpage}
\usepackage{compress}
\usepackage{simplemargins}
\setallmargins{1in}

\usepackage{hyperref}

\makeatletter
\renewcommand*\l@author[2]{}
\renewcommand*\l@title[2]{}
\makeatletter
\newcommand{\nocontentsline}[3]{}
\newcommand{\tocless}[2]{\bgroup\let\addcontentsline=\nocontentsline#1{#2}\egroup}

\def\Aut{{\rm Aut}}
\def\ct{{\rm CT}}
\def\Fix{{\rm Fix}}
\def\Ker{{\rm Ker}}
\def\O{\calO}
\def\id{{\rm id}}
\def\bo{\partial} 
\def\int{\mathring} 
\def\hookto{\hookrightarrow}

\def\rcover{\textsc{RegularCover}}
\def\cover{\textsc{$H$-Cover}}
\def\hom{\textsc{$H$-Hom}}

\def\ivmatch{{\sffamily IV}-\textsc{Matching}}
\def\ivsubgraph{{\sffamily IV}-subgraph}
\def\sat{\textsc{$3$-Sat}}

\newenvironment{packed_enum}{
	\begin{enumerate}
		\setlength{\itemsep}{1pt}
	    \setlength{\parskip}{0pt}
		\setlength{\parsep}{0pt}
}{\end{enumerate}}

\newenvironment{packed_itemize}{
	\begin{itemize}
		\setlength{\itemsep}{1pt}
	    \setlength{\parskip}{0pt}
		\setlength{\parsep}{0pt}
}{\end{itemize}}

\newenvironment{packed_head_enum}[1]{
	\begin{enumerate}[#1]
		\setlength{\itemsep}{1pt}
	    \setlength{\parskip}{0pt}
		\setlength{\parsep}{0pt}
}{\end{enumerate}}

\newcommand{\heading}[1]{\medskip\par\noindent{\bf #1}}

\def\computationproblem#1#2#3{
	\begin{center}
	\begin{tabular}{rp{10cm}}
	{\bf Problem:\enspace}&#1\\
	{\bf Input:\enspace}&#2\\
	{\bf Output:\enspace}&#3\\
	\end{tabular}
	\end{center}
}

\def\gS{\mathbb{S}} \def\gC{\mathbb{C}} \def\gA{\mathbb{A}} \def\gD{\mathbb{D}}

\def\calA{{\cal A}}  \def\calC{{\cal C}} 
   \def\calH{{\cal H}}
   \def\calL{{\cal L}}
\def\calM{{\cal M}}  \def\calO{{\cal O}} \def\calP{{\cal P}}
 \def\calR{{\cal R}}

   \def\frakL{{\mathfrak L}}

\def\cNP{\hbox{\rm \sffamily NP}}
\def\cGI{\hbox{\rm \sffamily GI}}

\def\cFPT{\hbox{\rm \sffamily FPT}}

\iftrue

\newcounter{lth}
\title{Algorithmic Aspects of Regular Graph Covers\\with Applications to Planar Graphs\thanks{
The conference version of this paper appeared in ICALP 2014.
The first three authors are supported by the ESF Eurogiga project GraDR as GA\v{C}R GIG/11/E023,
the fourth author by the ESF Eurogiga project GReGAS as the APVV project ESF-EC-0009-10 and by VEGA 1/0621/11.
The first author is also supported by the project Kontakt LH12095 and the second author by GA\v{C}R 14-14179S.
The second and the third authors are also supported by Charles University as GAUK 196213.}}
\author{Ji\v{r}\'i Fiala\inst{1}
		\and Pavel Klav\'ik\inst{2}
		\and Jan Kratochv\'il\inst{1}
		\and Roman Nedela\inst{3}
		}

\institute{Department of Applied Mathematics,\\
Faculty of Mathematics and Physics, Charles University,\\
Malostransk{\'e} n{\'a}m{\v e}st{\'\i} 25, 118 00 Prague, Czech Republic.\\
E-mails: \texttt{\{fiala,honza\}@kam.mff.cuni.cz}.\\[0.5em]\and
Computer Science Institute,\\
Faculty of Mathematics and Physics, Charles University,\\
Malostransk{\'e} n{\'a}m{\v e}st{\'\i} 25, 118 00 Prague, Czech Republic.\\
E-mail: \texttt{klavik@iuuk.mff.cuni.cz}.\\[0.5em]\and
Institute of Mathematics and Computer Science SAS\\
and Matej Bel University,\\
\v{D}umbierska 1, 974 11 Bansk\'a Bystrica, Slovak republic.\\
Email: \texttt{nedela@savbb.sk}.}

\begin{document}
\maketitle

\begin{abstract}
A graph $G$ \emph{covers} a graph $H$ if there exists a locally bijective homomorphism from $G$ to
$H$. We deal with \emph{regular covers} in which this locally bijective homomorphism is prescribed
by an action of a subgroup of $\Aut(G)$.  Regular covers have many applications in constructions and
studies of big objects all over mathematics and computer science.

\hskip 0.5cm We study \emph{computational aspects} of regular covers that have not been addressed
before. The decision problem $\rcover$ asks for two given graphs $G$ and $H$ whether $G$ regularly
covers $H$.  When $|H|=1$, this problem becomes Cayley graph recognition for which the complexity is
still unresolved. Another special case arises for $|G| = |H|$ when it becomes the graph isomorphism
problem. Therefore, we restrict ourselves to graph classes with polynomially solvable graph
isomorphism.  

\hskip 0.5cm Inspired by Negami, we apply the structural results used by Babai in the 1970's to study
automorphism groups of graphs. Our main result is the following \cFPT\ meta-algorithm: Let $\calC$ be a class
of graphs such that the structure of automorphism groups of 3-connected graphs in $\calC$ is simple.
Then we can solve $\rcover$ for $\calC$-inputs $G$ in time $\O^*(2^{e(H)/2})$ where $e(H)$ denotes
the number of the edges of $H$. As one example of $\calC$, this meta-algorithm applies to planar
graphs. In comparison, testing general graph covers is known to be \cNP-complete for planar inputs
$G$ even for small fixed graphs $H$ such as $K_4$ or $K_5$. Most of our results also apply to
general graphs, in particular the complete structural understanding of regular covers for 2-cuts.
\end{abstract}

\newpage
\tableofcontents
\newpage

\section{Introduction}

The notion of \emph{covering} originates in topology as a notion of local similarity of two
topological surfaces. For instance, consider the unit circle and the real line. Globally, these two
surfaces are not the same, they have different properties, different fundamental groups, etc.
But when we restrict ourselves to a small part of the circle, it looks the same as a small part
of the real line; more precisely the two surfaces are locally homeomorphic, and thus they share the
local properties. The notion of covering formalizes this property of two surfaces being
\emph{locally the same}.

\begin{figure}[b!]
\centering
\includegraphics{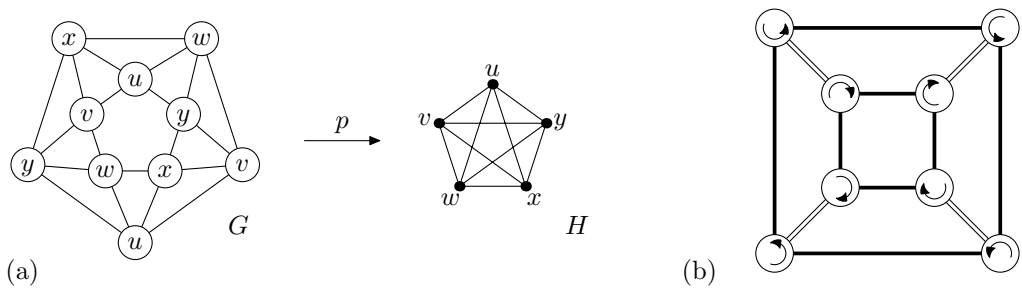}
\caption{(a) A covering projection $p$ from a graph $G$ to a graph $H$. (b) The Cayley
graph of the dihedral group $\gD_4$ generated by the $90^\circ$ rotations (in black) and the
reflection around the $x$-axis (in white).}
\label{fig:big_picture}
\end{figure}

More precisely, suppose that we have two  topological spaces: a big one $G$ and a small one $H$. We
say that $G$ \emph{covers} $H$ if there exists a mapping called a \emph{covering projection} $p : G
\to H$ which locally preserves the structure of $G$.  For instance, the mapping $p(x) = (\cos x,\sin
x)$ from the real line to the unit circle is a covering projection. The existence of a covering
projection ensures that $G$ looks locally the same as $H$; see Figure~\ref{fig:big_picture}a.

In this paper, we study coverings of graphs in a more restricting version called \emph{regular
covering}, for which the covering projection is described by an action of a group; see
Section~\ref{sec:preliminaries} for the formal definition. If $G$ regularly covers $H$, then we say
that $H$ is a \emph{quotient} of $G$.

\subsection{Applications of Graph Coverings}

Suppose that $G$ covers $H$ and we have some information about one of the objects. How much
knowledge does translate to the other object? It turns out that quite a lot, and this makes
covering a powerful technique with many diverse applications. The big advantage of regular coverings
is that they can be efficiently described and many properties easily translate between the objects.
We sketch some applications now.

\heading{Powerful Constructions.} The reverse of covering called \emph{lifting} can be applied to
small objects in order to construct large objects of desired properties. For instance, the
well-known Cayley graphs are large objects which can be described easily by a few elements of a
group. Let $G$ be a Cayley graph generated by elements $g_1,\dots,g_e$ of a group $\Gamma$. The
vertices of $G$ correspond to the elements of $\Gamma$ and the edges are described by actions of
$g_1,\dots,g_e$ on $\Gamma$ by left multiplication; each $g_i$ defines a permutation on $\Gamma$ and
we put edges along the cycles of this permutation. See Figure~\ref{fig:big_picture}b for an example.
Cayley graphs were originally invented to study the structure of groups~\cite{cayley}.

\begin{figure}[t!]
\centering
\includegraphics{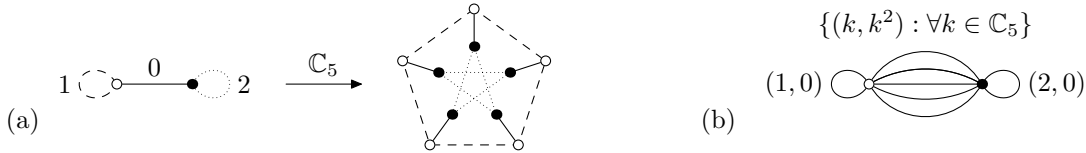}
\caption{(a) A construction of the Petersen graph by lifting with the group $\gC_5$. (b) By lifting
the described graph with the group $\gC_5^2$, we get the Hoffman-Singleton graph. The five parallel
edges are labeled $(0,0)$, $(1,1)$, $(2,4)$, $(3,4)$ and $(4,1)$.}
\label{fig:lifting_graphs}
\end{figure}

In the language of coverings, every Cayley graph $G$ with an involution-free generating set can be
described as a lift of a one vertex graph $H$ with $e$ loops labeled $g_1,\dots,g_e$. Regular covers
can be viewed as a generalization of Cayley graphs where the small graph $H$ can contain more then
one vertex. For example, the famous Petersen graph can be constructed as a lift of a two vertex
graph $H$, see Figure~\ref{fig:lifting_graphs}a. These two vertices are necessary as it is known
that Petersen graph is not a Cayley graph. Figure~\ref{fig:lifting_graphs}b shows a simple
construction~\cite{hoffman_construction} of the Hoffman-Singleton graph~\cite{hoffman_singleton}
which is a 7-regular graph with 50 vertices. Notice that from this construction it is immediately
apparent that the Hoffman-Singleton graph contains many induced copies of the Petersen graph. (In
fact, it contains 525 copies of it.)

The Petersen  and the Hoffman-Singleton graphs are extremal graphs for the
\emph{degree-diameter} problem: Given integers $d$ and $k$, find a maximal graph $G$ with diameter
$d$ and degree $k$. In general, the size of $G$ is not known. Many currently best constructions are
obtained using the covering techniques~\cite{miller2005moore}. 

Further applications employ the fact that nowhere-zero flows, vertex and edge colorings, eigenvalues
and other graph invariants lift along a covering projection. In history, two main applications are
the solution of the Heawood map coloring problem~\cite{ringel_youngs,gross_tucker} and construction
of arbitrary large highly symmetrical graphs~\cite{biggs}.

\heading{Models of Local Computation.} These and similar constructions have many practical
applications in designing highly efficient computer networks~\cite{network0,network1,network2,%
network3,network4,network5,network6,network7} since these networks can be efficiently
described/constructed and have many strong properties. In particular, networks based on covers of
simple graphs allow fast parallelization of computation as described e.g.
in~\cite{bodlaender,Ang1,Ang2}.

\heading{Simplifying Objects.} Regular covering can be also applied in the opposite way, to project
big objects onto smaller ones while preserving some properties. One way is to represent a class of
objects satisfying some properties as quotients of the universal object of this property.  For
instance, this was used in the study of arc-transitive cubic graphs~\cite{goldschmidt}, and the key
point is that universal objects are much easier to work with. This idea is commonly used in fields
such as Riemann surfaces~\cite{farkas_kra} and theoretical physics~\cite{katanaev}.

\subsection{Complexity Aspects}

In the constructions we described, the covers are regular and satisfy additional algebraic
properties.  The reason is that regular covers are easier to describe.  In this paper, we initiate
the study of the computational complexity of regular covering.

\computationproblem
{\rcover}
{Connected graphs $G$ and $H$.}
{Does $G$ regularly cover $H$?}

\heading{Relations to Covers.} This problem is closely related to the complexity of general covering
which was widely studied before. We try to understand how much the additional algebraic structure
changes the computational complexity. Study of the complexity of general covers was pioneered by
Bodlaender~\cite{bodlaender} in the context of networks of processors in parallel computing, and for
fixed target graph was first asked by Abello et al.~\cite{AFS}. The problem \cover\ asks for an
input graph $G$ whether it covers a fixed graph $H$. The general complexity is still unresolved but
papers~\cite{kratochvil97,fiala00} show that it is \cNP-complete for every $r$-regular graph $H$
where $r \ge 3$. (For a survey of the complexity results, see~\cite{FK}.)

The complexity results concerning graph covers are mostly \cNP-complete. In our impression, the
additional algebraic structure of regular graph covers makes the problem easier, as shown by the
following two contrasting results. The problem \cover\ remains \cNP-complete for several small
fixed graphs $H$ (such as $K_4$, $K_5$) even for planar inputs $G$~\cite{planar_covers}. On the
other hand, our main result is that the problem \rcover\ is fixed-parameter tractable in the
number of edges of $H$ for every planar graph $G$ and for every $H$.

Two additional problems of finding lifts and quotients, closely related to \rcover, are considered
in Section~\ref{sec:complexity_properties}.

\heading{Relations to Cayley Graphs and Graph Isomorphism.}
The notion of regular covers builds a bridge between two seemingly different problems. If the graph
$H$ consists of a single vertex, it corresponds to recognizing Cayley graphs which is still open; a
polynomial-time algorithm is known only for circulant graphs~\cite{circulant_recognition}. When both
graphs $G$ and $H$ have the same size, we get graph isomorphism testing. Our results are far from
this but we believe that better understanding of regular covering can also shed new light on these
famous problems.

Theoretical motivation for studying graph isomorphism is very similar to \rcover. For practical
instances, one can solve the isomorphism problem very efficiently using various heuristics. But
polynomial-time algorithm working for all graph is not known and it is very desirable to understand
the complexity of graph isomorphism.  It is known that testing graph isomorphism is equivalent to
testing isomorphism of general mathematical structures~\cite{hedrlin}.  The notion of isomorphism is
widely used in mathematics when one wants to show that two seemingly different mathematical
structures are the same. One proceeds by guessing a mapping and proving that this mapping is an
isomorphism. The natural complexity question is whether there is a better way in which one
algorithmically derives an isomorphism. Similarly, regular covering is a well-known mathematical
notion which is algorithmically interesting and not understood.

Further, a regular covering is described by a semiregular subgroup of the automorphism group
$\Aut(G)$. Therefore it seems to be closely related to computation of $\Aut(G)$ since one should
have a good understanding of this group first to solve the regular covering problem. The problem of
computing automorphism groups is known to be closely related to graph isomorphism.

\heading{Homomorphisms and CSP.}
Since regular covering is a restricted locally bijective homomorphism, we give an overview of
complexity results concerning homomorphisms and general coverings.  Hell and
Ne\v{s}et\v{r}il~\cite{hell_nesetril} studied the problem \hom\ which asks whether there exists a
homomorphism between an input graph $G$ and a fixed graph $H$. Their celebrated dichotomy result
for simple graphs states that the problem \hom\ is polynomially solvable if $H$ is bipartite and it is
\cNP-complete otherwise. Homomorphisms can be described in the language of constraint satisfaction (CSP), 
and the famous dichotomy conjecture~\cite{csp_dichotomy} claims that every CSP is either
polynomially solvable, or \cNP-complete.

\subsection{Our Results}

Let $\calC$ be a class of connected multigraphs. By $\calC / \Gamma$ we denote the class of all
regular quotients of graphs of $\calC$ (note that $\calC \subseteq \calC / \Gamma$). For instance for
$\calC$ equal to planar graphs, the class $\calC / \Gamma$ is -- according to the Negami's
Theorem~\cite{negami} -- the class of projective planar graphs. We consider four properties of
$\calC$, for formal definitions see Section~\ref{sec:preliminaries}:
\begin{packed_head_enum}{(P1)}
\item[(P0)] The class $\calC$ is  closed under taking subgraphs and under replacing connected
components attached to 2-cuts by edges. 
\item[(P1)] The graph isomorphism problem is solvable in polynomial time for $\calC$ and $\calC /
\Gamma$.
\item[(P2)] For a 3-connected graph $G \in \calC$, the group $\Aut(G)$ and all its semiregular
subgroups $\Gamma$ can be computed in polynomial time. Here by semiregularity, we mean that the
action of $\Gamma$ has no non-trivial stabilizers of the vertices.
\item[(P3)] Let $G$ and $H$ be 3-connected graphs of $\calC / \Gamma$, possibly with colored and
directed edges. Let the vertices of $G$ be further colored by $c(u), u \in V(G)$, and let $H$ be
equipped with a list $\frakL(u)$ of possible colors for each vertex $u \in V(H)$ (the coloring is
not necessarily proper).  We can test in polynomial time whether there exists a color-compatible
isomorphism $\pi : G \to H$, i.e. an isomorphism such that the colors and orientations of edges are
preserved and for every $u \in V(G)$, we have $c(u) \in \frakL(\pi(u))$. (The existence of such an
isomorphism is denoted by $G\hookto H$.)
\end{packed_head_enum}
As we prove in Section~\ref{sec:planar_graphs}, these four properties are tailored for the class of
planar graphs. (But the proof of the property (P3) is non-trivial, based on the result
of~\cite{pp_iso}.) Negami's Theorem~\cite{negami} dealing with regular covers of planar graphs is
one of the oldest results in topological graph theory; therefore we decided to start the study of
computational complexity of \rcover\ for planar graphs. Our algorithm, however, applies to a wider
class of graphs.

We use the complexity notation $f = \O^*(g)$ which omits polynomial factors. Our main result is
the following \cFPT\ meta-algorithm:

\begin{theorem} \label{thm:metaalgorithm}
Let $\calC$ be a class of graphs satisfying (P0) to (P3). Then there is an \cFPT\ algorithm for
\rcover\ for $\calC$-inputs $G$ in the parameter $e(H)$, running in time $\O^*(2^{e(H)/2})$ where
$e(H)$ is the number of edges of $H$.
\end{theorem}

It is important that most of our results apply to general graphs. We wanted to generalize the result
of Babai~\cite{babai1975automorphism} which states that it is sufficient to solve graph isomorphism
for 3-connected graphs. Our main goal was to understand how regular covering behaves with respect to
vertex 1-cuts and 2-cuts. Concerning 1-cuts, regular covering behaves non-trivially only on the
central block of $G$, so they are easy to deal with. But we show that regular covering can behave
highly complex on 2-cuts. From structural point of view, we give a complete description of this
behaviour. Algorithmically, we solve computation only partially and we need several other
assumptions to get an efficient algorithm.

Planar graphs are very important and also well studied in connection to coverings. Negami's
Theorem~\cite{negami} dealing with regular covers of planar graphs is one of the oldest results in
topological graph theory; therefore we decided to start the study of computational complexity of
\rcover\ for planar graphs.  In particular, our theory applies to planar graphs since they satisfy
(P0) to (P3).

\begin{corollary} \label{cor:planar_rcover}
For a planar graph $G$, \rcover\ can be solved in time $\O^*(2^{e(H)/2})$.
\end{corollary}

\heading{Our Approach.}
We quickly sketch our approach. The meta-algorithm proceeds by a series of \emph{reductions}
replacing parts of the graphs by edges. These reductions are inspired by the approach of
Negami~\cite{negami} and turn out to follow the same lines as the reductions introduced by Babai for
studying automorphism groups of planar graphs~\cite{babai1975automorphism,babai1996automorphism}.
Since the key properties of the automorphism groups are preserved by the reductions, computing
automorphism groups can be reduced to computing them for 3-connected graphs~\cite{babai1975automorphism}.
In~\cite{hopcroft_tarjan_dividing,hopcroft_tarjan_planar_iso}, this is used to compute automorphism
groups of planar graphs since the automorphism groups of 3-connected planar graphs are the
automorphism groups of tilings of the sphere, and are well-understood.

The \rcover\ problem is more complicated, and we use the following novel approach.  When the
reductions reach a 3-connected graph, the natural next step is to compute all its quotients; there
are polynomially many of them according to (P2). What remains is the most difficult part: To test
for each quotient whether it corresponds to $H$ after unrolling the reductions. This process is
called \emph{expanding} and the issue here is that there may be exponentially many different ways to
expand the graph, so we have to test in a clever way whether it is possible to reach $H$. Our
algorithm consists of several subroutines, most of which we indeed can perform in polynomial time.
Only one subroutine (finding a certain ``generalized matching'') we have not been able to
solve in polynomial time.

This slow subroutine can be avoided in some cases:

\begin{corollary} \label{cor:simple_cases}
If the $\calC$-graph $G$ is 3-connected or if $k=|G|/|H|$ is odd, then the meta-algorithm of
Theorem~\ref{thm:metaalgorithm} can be modified to run in polynomial time.
\end{corollary}

\heading{Structure.} This paper is organized as follows. In Section~\ref{sec:preliminaries}, we
introduce the formal notation used in this paper. In Section~\ref{sec:atoms}, we introduce
\emph{atoms} which are the key objects of this paper. In Section~\ref{sec:reduction_and_expansion},
we describe structural properties of reductions via atoms, and expansions of constructed quotient
graphs. In Section~\ref{sec:algorithm}, we use these structural properties to create the
meta-algorithm of Theorem~\ref{thm:metaalgorithm}. Finally, in Section~\ref{sec:planar_graphs} we
deal with specific properties of planar graphs and show that the class of planar graphs satisfies
(P0) to (P3). In Conclusions, we describe open problems and possible extensions of our results. See
Section~\ref{sec:structure} for more detailed overview of the main steps.

\section{Definitions and Preliminaries} \label{sec:preliminaries}

A multigraph $G$ is a pair $(V(G),E(G))$ where $V(G)$ is a set of vertices and $E(G)$ is a multiset
of edges. We denote $|V(G)|$ by $v(G)$ and $|E(G)|$ by $e(G)$.  The graph can possibly contain
parallel edges and loops, and each loop at $u$ is incident twice with the vertex $u$. (So it
contributes by two to the degree of $u$.) Each edge $e = uv$ gives rise to two half-edges, one
attached to $u$ and the other to $v$. We denote by $H(G)$ the collection of all half-edges. We
denote $|H(G)|$ by $h(G)$ and clearly $h(G) = 2e(G)$.  As quotients, we sometime obtain graphs
containing half-edges with free ends (missing the opposite half-edges).

We consider graphs with colored edges and also with three different edge types (directed
edges, undirected edges and a special type called halvable edges).  It might seem strange to
consider such general objects. But when we apply reductions, we replace parts of the graph by edges
and the colors encode isomorphism classes of replaced parts. This allows the algorithm to work with
smaller reduced graphs and deduce some structure of the original large graph. So even if the input
graphs $G$ and $H$ are simple, the more complicated multigraphs are naturally constructed.

\subsection{Automorphisms and Groups}

\heading{Automorphisms.} We state the definitions in a very general setting of multigraphs and
half-edges. An  \emph{automorphism} $\pi$ is fully described by a permutation $\pi_h : H(G) \to
H(G)$ preserving edges and incidences between half-edges and vertices. Thus, $\pi_h$ induces two
permutations $\pi_v : V(G) \to V(G)$ and $\pi_e : E(G) \to E(G)$ connected together by the very
natural property $\pi_e(uv) = \pi_v(u) \pi_v(v)$ for every $uv \in E(G)$.  In most of situations, we
omit subscripts and simply use $\pi(u)$ or $\pi(uv)$. In addition, when we work with colored graphs,
we require that an automorphism preserves the colors.

\heading{Groups.} We assume that the reader is familiar with basic properties of groups;
otherwise see~\cite{rotman}. We denote groups by Greek letters as for instance $\Gamma$.  We use the following
notation for standard groups:
\begin{packed_itemize}
\item $\gS_n$ for the symmetric group of all $n$-element permutations,
\item $\gC_n$ for the cyclic group of integers modulo $n$,
\item $\gD_n$ for the dihedral group of the symmetries of a regular $n$-gon, and
\item $\gA_n$ for the group of all even $n$-element permutations.
\end{packed_itemize}

\heading{Automorphism Groups.} Groups are quite often studied in the context of group actions, since
their origin is in studying symmetries of mathematical objects. A group $\Gamma$ \emph{acts} on a
set $S$ in the following way. Each $g \in \Gamma$ permutes the elements of $S$, and the \emph{action}
is described by a mapping $\cdot : \Gamma \times S \to S$ usually satisfying further properties that 
arise naturally from the structure of $S$.

In the language of graphs, an example of such an action is the group of all automorphisms of $G$,
denoted by $\Aut(G)$.  Each element $\pi \in \Aut(G)$ acts on $G$, permutes its vertices, edges and
half-edges while it preserves edges and incidences between the half-edges and the vertices.

The \emph{orbit} $[v]$ of a vertex $v \in V(G)$ is the set of all vertices $\{\pi(v) \mid \pi \in
\Gamma\}$, and the orbit $[e]$ of an edge $e \in E(G)$ is defined similarly as $\{\pi(e) \mid \pi
\in \Gamma\}$. The \emph{stabilizer} of $x$ is the subgroup of all automorphisms which fix $x$.
An action is called \emph{semiregular} if it has no non-trivial (i.e., non-identity) stabilizers of
both vertices and half-edges. Note that the stabilizer of an edge in a semiregular action may be
non-trivial, since it may contain an involution transposing the two half-edges.  We say that a group
is \emph{semiregular} if the associated action is semiregular.

\subsection{Coverings}

\begin{figure}[b!]
\centering
\includegraphics{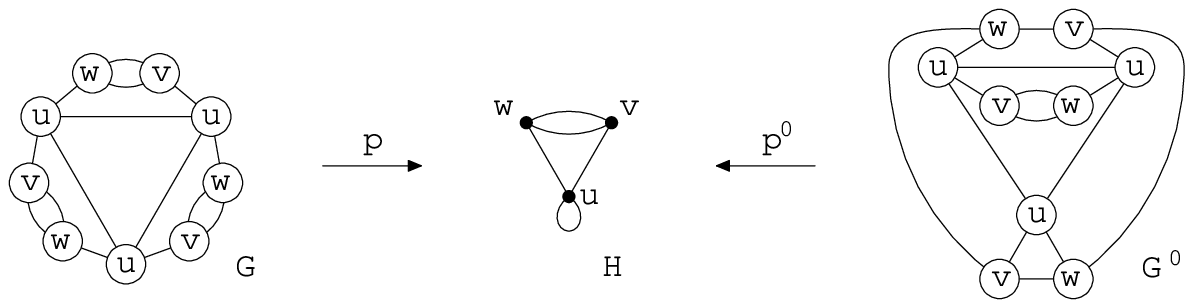}
\caption{Two covers of $H$. The projections $p_v$ and $p'_v$ are written inside of the vertices, and
the projections $p_e$ and $p'_e$ are omitted. Notice that each loop is realized by having two
neighbors labeled the same, and parallel edges are realized by having multiple neighbors labeled the
same. Also covering projections preserve degrees.}
\label{fig:cover_examples}
\end{figure}

A graph $G$ \emph{covers} a graph $H$ (or $G$ is a \emph{cover} of $H$) if there exists a locally
bijective homomorphism $p$ called a \emph{covering projection}.  A homomorphism $p$ from $G$ to $H$
is given by a mapping $p_h : H(G) \to H(H)$ preserving edges and incidences between half-edges and
vertices. It induces two mappings $p_v : V(G) \to V(H)$ and $p_e : E(G) \to E(H)$ such that $p_e(uv)
= p_v(u)p_v(v)$ for every $uv \in E(G)$.  The property to be local bijective states that for every
vertex $u \in V(G)$ the mapping $p_h$ restricted to the half-edges incident with $u$ is a bijection.
Figure~\ref{fig:cover_examples} contains two examples of graph covers. Again, we mostly omit
subscripts and just write $p(u)$ or $p(e)$.  A \emph{fiber} of a vertex $v \in V(H)$ is the set
$p^{-1}(v)$, i.e., the set of all vertices $V(G)$ that are mapped to $v$, and similarly for fibers
over edges and half-edges.

\begin{figure}[t!]
\centering
\includegraphics{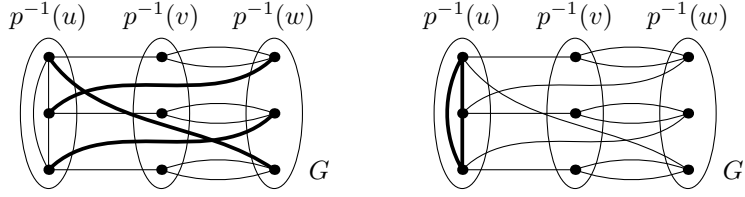}
\caption{The graph $G$ from Figure~\ref{fig:cover_examples} depicted by fibers of $p$. On the left,
we show that $p^{-1}(uw)$ gives a matching between the fibers of $u$ and $w$. On the right, the loop
around $u$ gives a cycle in $p^{-1}(u)$.}
\label{fig:fibers_and_lifts}
\end{figure}

\heading{The Unique Walk Lifting Property.} Let $uv \in E(H)$ be an edge which is not a loop. Then
the set $p^{-1}(uv)$ corresponds to a perfect matching between the fibers $p^{-1}(u)$ and
$p^{-1}(v)$. And if $uu \in E(H)$ is a loop, then the set $p^{-1}(uu)$ is a union of disjoint cycles
which cover exactly $p^{-1}(u)$. Figure~\ref{fig:fibers_and_lifts} shows examples. In general for a
subgraph $H'$ of $H$, the correspondence between $H'$ and its preimages $p^{-1}(H')$ in $G$ is
called \emph{lifting}.

Let $W$ be a walk $u_0e_1u_1e_2\dots e_nu_n$ in $H$.  Then $p^{-1}(W)$ consists of $k$ copies of
$W$. Suppose that we fix as the start one vertex in the fiber $p^{-1}(u_0)$. Due to the local
bijectiveness of $p$, there is exactly one edge incident with it which $p$ maps to $e_1$, so
$p^{-1}(u_1)$ is now uniquely determined, and so on for $p^{-1}(u_2)$ and the other vertices of $W$.
When we proceed in this way, the rest of the walk is determined.  This important property of every
covering is called the \emph{unique walk lifting property}.

We adopt the standard assumption that both $G$ and $H$ are connected.  Then as a simple corollary we
get that all fibers of $p$ are of the same size. To see that, observe that a path in $H$ is lifted
to disjoint paths in $G$. For $u,v \in V(H)$, consider a path $P$ between them.  Then the paths in
$p^{-1}(P)$ define a bijection between the fibers $p^{-1}(u)$ and $p^{-1}(v)$. In other words, $|G|
= k|H|$ for some $k \in \mathbb N$ which is the size of each fiber, and we say that $G$ is a
\emph{$k$-fold cover} of $H$.

\heading{Covering Transformation Groups.} Every covering projection $p$ defines a special subgroup of
$\Aut(G)$ called \emph{covering transformation group} $\ct(p)$. It consists of all automorphisms
$\pi$ which preserve the fibers of $p$, i.e., for every $u \in V(G)$, the vertices $u$ and $\pi(u)$
belong to one fiber. Consider the graphs from Figure~\ref{fig:cover_examples}. For the graph $G$, we
have $\Aut(G) = \gD_3$ and $\ct(p) = \gC_3$. And $G'$ has $\Aut(G') = \gC_2$ but $\ct(p')$ is
trivial; and we note that it is often the case that a covering projection $p$ has only
one fiber-preserving automorphism, the trivial one.

Now suppose that $\pi \in \ct(p)$.  Observe that a single choice of the image $\pi(u)$ of one vertex
$u \in V(G)$ fully determines $\pi$. This follows from the unique walk lifting property.  Let $v \in
V(G)$ and consider some path $P_{u,v}$ connecting $u$ and $v$ in $G$. This path corresponds to a path
$P = p(P_{u,v})$ in $H$. Now we lift $P$ and according to the unique walk lifting property, there
exists a unique path $P_{\pi(u),x}$ which starts in $\pi(u)$. But since $\pi$ is an automorphism and
it maps $P_{u,v}$ to $P_{\pi(u),x}$, then $x$ has to be equal $\pi(v)$. In other words, we just proved
that $\ct(p)$ is semiregular.

\heading{Regular Coverings.} We want to consider coverings which are highly symmetrical. For
examples from Figure~\ref{fig:cover_examples}, the covering $p$ is more symmetric than $p'$. The size
of $\ct(p)$ is a good measure of symmetricity of the covering $p$. Since $\ct(p)$ is semiregular,
it easily follows that $|\ct(p)| \le k$ for any $k$-fold covering $p$. A covering $p$ is
\emph{regular} if $|\ct(p)| = k$. In Figure~\ref{fig:cover_examples}, the covering $p$ is regular
since $|\ct(p)| = 3$, and the covering $p'$ is not regular since $|\ct(p')| = 1$.

We use the following definition of regular covering. Let $\Gamma$ be any semiregular subgroup
of $\Aut(G)$. It defines a graph $G / \Gamma$ called a \emph{regular quotient} (or simply
\emph{quotient}) of $G$ as follows: The vertices of $G / \Gamma$ are the orbits of the action
$\Gamma$ on $V(G)$, the half-edges of $G / \Gamma$ are the orbits of the action $\Gamma$ on $H(G)$.
A vertex-orbit $[v]$ is incident with a half-edge-orbit $[h]$ if and only if the vertices of $[v]$
are incident with the half-edges of $[h]$. (Because the action of $\Gamma$ is semiregular, each
vertex of $[v]$ is incident with exactly one half-edge of $[h]$, so this is well defined.) We
naturally construct $p : G \to G / \Gamma$ by mapping the vertices to its vertex-orbits and
half-edges to its half-edge-orbits. Concerning an edge $e \in E(G)$, it is mapped to an edge of $G /
\Gamma$ if the two half-edges belong to different half-edge-orbits of $\Gamma$. If they belong to
the same half-edge-orbits, it corresponds to a standalone half-edge of $G / \Gamma$.

\begin{figure}[b!]
\centering
\includegraphics{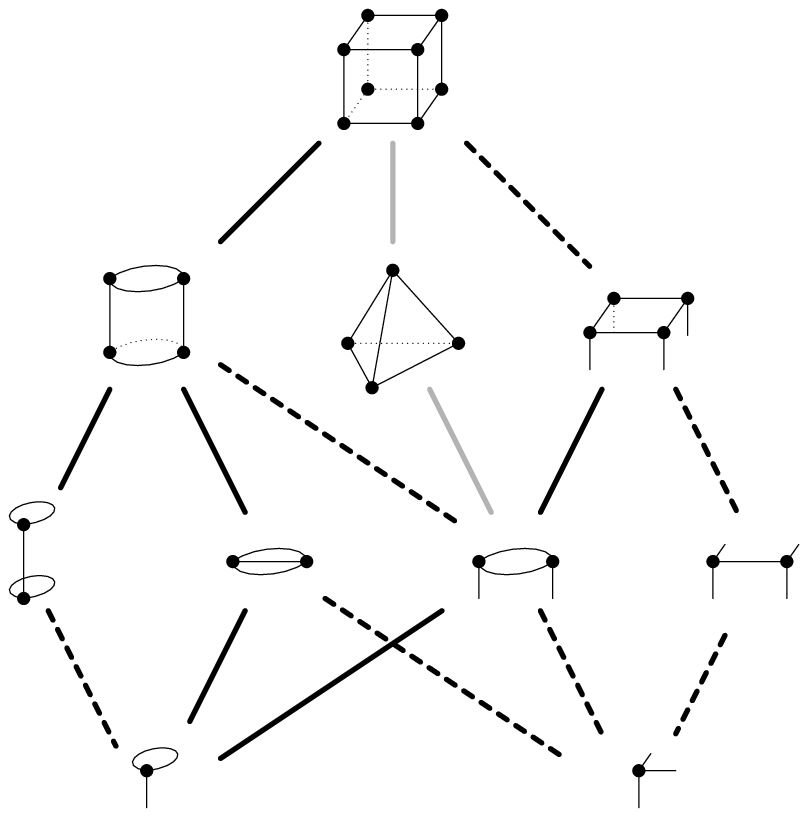}
\caption{The Hasse diagram of all quotients of the cube graph depicted in a geometric way. When
semiregular actions fix edges, the quotients contain half-edges. The quotients connected by bold
edges are obtained by 180 degree rotations. The quotients connected by dashed edges are obtained by
reflections. The tetrahedron is obtained by the antipodal symmetry of the cube, and its quotient is
obtained by a 180 degree rotation with the axis going through the centers of two non-incident edges
of the tetrahedron.}
\label{fig:quotients_of_cube}
\end{figure}

Since $\Gamma$ acts semiregularly on $G$, one can prove that $p$ is a $|\Gamma|$-fold regular
covering with $\ct(p) = \Gamma$. For the graphs $G$ and $H$ of Figure~\ref{fig:cover_examples}, we
get $H \cong G / \Gamma$ for $\Gamma \cong \gC_3$ which ``rotates the cycle by three vertices''.
As a further example, Figure~\ref{fig:quotients_of_cube} geometrically depicts all quotients of the
cube graph.

\subsection{Fundamental Complexity Properties of Coverings} \label{sec:complexity_properties}

We establish fundamental complexity properties of regular covering and also general covering. Our
goal is to highlight similarities with the graph isomorphism problem. Also we discuss other variants
of the \rcover\ problem.

\heading{Belonging to NP.} The general \cover\ problem is clearly in \cNP\ since one can just test in
polynomial time whether a given mapping is a locally bijective homomorphism. Not so obviously, the
same holds for the \rcover\ problem. 

\begin{lemma} \label{lem:rcover_in_np}
The problem \rcover\ is in \cNP.
\end{lemma}

\begin{proof}
One just needs to use a suitable definition of regular covering. As stated above, $G$ regularly
covers $H$, if and only if there exists a semiregular subgroup $\Gamma$ of $\Aut(G)$ such that $G /
\Gamma \cong H$. As the certificate, we just give $k$ permutations, one for each element of
$\Gamma$, and the isomorphism between $G / \Gamma$ and $H$. We can easily check whether these $k$
permutations define a group $\Gamma$, and whether $\Gamma$ acts semiregularly on $G$.  Further, the
given isomorphism allows to check whether the constructed $G / \Gamma$ is isomorphic to $H$.
Clearly, this certificate is polynomially large and can be verified in polynomial time.\qed
\end{proof}

One can prove even a stronger result:

\begin{lemma} \label{lem:testing_regularity}
For a mapping $p : G \to H$, we can test whether it is a regular covering in polynomial time.
\end{lemma}

\begin{proof}
Testing whether $p$ is a covering can clearly be done in polynomial time. It remains to test
regularity. Choose an arbitrary spanning tree $T$ of $H$. Since $p$ is a covering, then $p^{-1}(T)$
is a disjoint union of $k$ isomorphic copies $T_1,\dots,T_k$ of $T$. We number the vertices of the
fibers according to the spanning trees, i.e., $p^{-1}(v) = \{v_1,\dots,v_k\}$ such that $v_i \in
T_i$. This induces a numbering of the half-edges of each fiber over a half-edge of $H(H)$, following
the incidences between half-edges and vertices.  For every half-edge $h \notin H(T)$, we define a
permutation $\sigma_h$ of $\{1,\dots,k\}$ taking $i$ to $j$ if there is a half-edge $h$ in
$p^{-1}(h)$ such that the edge $u_iv_j$ corresponds to $h$.

It remains to test whether the size of the group $\Theta$ generated by all $\sigma_h$, where $h \notin H(T)$,
is of size exactly $k$. From the theory of permutation groups, since $G$ is assumed to be connected, it follows
that the action of $\Theta$ is transitive. Therefore its size is at least $k$, and the action is
regular if and only if it is exactly $k$.\qed
\end{proof}

The constructed permutations $\sigma_h$ associated with $p$ are known in the
literature~\cite{gross_tucker} as \emph{permutation voltage assigments} associated with $p$.

\heading{Other Variants.} In the \rcover\ problem, the input gives two graphs $G$ and $H$ and we ask
for an existence of a regular covering from $G$ to $H$. There are two other reasonable variants of
this problem we discuss now. The input can specify only one of the two graphs and ask for existence
of the other graph of, say, a given size.

First, suppose that only $H$ is given and we ask whether a $k$-fold cover $G$ of $H$ exists.  This
is called \emph{lifting} and the answer is always positive. The theory of covering describes a
technique called voltage assignment which can be applied to generate all $k$-folds $G$.  We do not
deal with lifting in this paper, but there are nevertheless many interesting computational questions
with important applications. For instance, one can try to generate efficiently all lifts up to
isomorphism; this is not trivial since different voltage assignments might lead to isomorphic
graphs. Also, one may ask for existence of a lift with some additional properties.

The other variant gives only $G$ and asks for existence of a quotient $H$ which is regularly covered
by $G$ and $|G| = k|H|$.  This problem is \cNP-complete even for fixed $k=2$, proved in a different
language by Lubiw~\cite{lubiw}. Lubiw shows that testing existence of a fixed-point free involutory
automorphism is \cNP-complete which is equivalent to existence of a half-quotient $H$. We sketch
hers reduction from \sat. Each variable is represented by an even cycle attached to the rest of the
graph.  Each cycle has two possible regular quotients, either a cycle of half length (obtained by
the $180^\circ$ rotation), or a path of half length with attached half-edges (obtained by a
reflection through opposite edges). Each of these quotients represents one truth assignment of the
corresponding variable. To distinguish variables, distinct gadgets are attached to the cycles. These
variable gadgets are attached to clause gadgets. Naturally, one can construct a quotient of the
clause gadget if and only if at least one literal is satisfied.

One should ask whether this reduction also implies \cNP-completeness for the \rcover\ problem. Since
the input gives also a graph $H$, one can decode the assignment of the variables from it, and thus
this reduction does not work. We conjecture that no similar reduction with a fixed $k$ can be
constructed for \rcover\ since we believe that for a fixed $k$ the problem of counting the number
of regular coverings between $G$ and $H$ can be solved using polynomially-many instances of \rcover.
In complexity theory, it is believed that the counting version of no \cNP-complete problem satisfies
this. Similar evidence was used by Mathon~\cite{mathon_isocount} to show that graph isomorphism is
unlikely \cNP-complete, and as a work in progress we believe that a similar argument can be applied
to \rcover.


The results of this paper also show that the reduction of Lubiw cannot be modified for planar inputs
$G$. Our algorithmic and structural insights allow an efficient enumeration of all quotients $H$
of a given planar graph $G$. On the other hand, the hardness result of Lubiw states that to solve the
\rcover\ problem in general, one has to work with both graphs $G$ and $H$ from beginning. Our
algorithm starts only with $G$ and tries to match its quotients to $H$ only in the end.
Nevertheless, some modifications in this directions, not necessary for planar graphs, would be
possible.

\subsection{Overview of the Main Steps} \label{sec:structure}

We now give a quick overview of the paper.

The main idea is the following. If the input graph $G$ is 3-connected, using our assumptions the
\rcover\ problem is trivially solvable. Otherwise, we proceed by a series of reductions, replacing
parts of the graph by edges, essentially forgetting details of the graph. We end-up with a primitive
graph which is either 3-connected, or very simple (a cycle or $K_2$). The reductions are done in
such a way that no essential information of semiregular actions is lost.

Inspired by Negami~\cite{negami} and Babai~\cite{babai1996automorphism}, we introduce in
Section~\ref{sec:atoms} the most important definition of an \emph{atom}. Atoms are
inclusion-minimal subgraphs which cannot be further simplified and are essentially
3-connected.  Our strategy for reductions is to detect the atoms and replace them by edges. In
this process, we remove details from the graph but preserve its overall structure. Our definition of
atoms is quite technical, dealing with many details necessary for the next sections.

When the graph $G$ is not 3-connected, we consider its block-tree. The central block plays the
key role in every regular covering projection. The reason is that the covering behaves non-trivially
only on this central block; the remaining blocks are isomorphically preserved in $H$. Therefore the
atoms are defined with respect to the central block.  We distinguish three types of atoms:
\begin{packed_itemize}
\item \emph{Proper atoms} are inclusion-minimal subgraphs separated by a 2-cut inside a
block.
\item \emph{Dipoles} are formed by the sets of all parallel edges joining two vertices.
\item \emph{Block atoms} are blocks which are leaves of the block-tree, or stars of all pendant edges
attached to a vertex. The central block is never a block atom.
\end{packed_itemize}

In Section~\ref{sec:reduction_and_expansion}, we deal with two transformations of graphs called
\emph{reduction} and \emph{expansion}. The graph $G$ is reduced by replacing its atoms by edges; for
proper atoms and dipoles inside the blocks, for block atoms by pendant edges. Further these edges
are colored according to isomorphism classes of the atoms. This way the reduction omits unimportant
details from the graph but its key structure is preserved. We apply a series of reductions $G =
G_0,\dots,G_r$ till we obtain a graph $G_r$ called \emph{primitive} which contains no atoms. We show
that the reduction preserves essentially the automorphism group. More precisely, $\Aut(G_i)$ is a
factor-group of $\Aut(G_{i-1})$; the action of $\Aut(G_i)$ inside the atoms is lost by the
factorization. 

The other transformation called expansion is applied to the quotient graphs.  The goal of
expansion is to revert the reduction, so it replaces colored edges back by atoms. To do this, we
have to understand how regular covering behaves with respect to atoms. Inspired by
Negami~\cite{negami}, we show that each proper atom/dipole has three possible types of quotients
that we call an \emph{edge-quotient}, a \emph{loop-quotient} and a \emph{half-quotient}. The
edge-quotient and the loop-quotient are uniquely determined but an atom may have many non-isomorphic
half-quotients.

When the primitive graph $G_r$ is reached, all semiregular subgroups $\Gamma_r$ of $\Aut(G_r)$ are
computed and for each one a quotient $H_r = G_r / \Gamma_r$ is constructed. Our goal is to
understand all graphs $H_0$ to which $H_r$ can be expanded, as depicted in the following diagram:
$$
\begin{gathered}
\xymatrix{
G_0 \ar[d]_{\Gamma_0} \ar[r]^{\rm red.} & G_1 \ar[d]_{\Gamma_1} \ar[r]^{\rm red.} &
\cdots \ar[r]^{\rm red.}&G_i \ar[d]_{\Gamma_i} \ar[r]^{\rm red.} &
G_{i+1}\ar[d]_{\Gamma_{i+1}} \ar[r]^{\rm red.} & \cdots \ar[r]^{\rm red.} &
G_r \ar[d]_{\Gamma_r}\\%
H_0 									& H_1 \ar[l]_{\rm exp.} &
\cdots \ar[l]_{\rm exp.}&H_i \ar[l]_{\rm exp.} &
H_{i+1}\ar[l]_{\rm exp.}					   & \cdots \ar[l]_{\rm exp.} &
H_r \ar[l]_{\rm exp.}\\%
}
\end{gathered}
$$

The constructed quotients contain colored edges, loops and half-edges corresponding to atoms.  Each
half-edge in $H_r$ is created from a halvable edge if an automorphism of $\Gamma_r$ fixes this
halvable edge and exchanges its endpoints. Roughly speaking it corresponds to cutting the edge in
half. We show that every possible expansion of a quotient $H_{i-1}$ from $H_i$ can be constructed by
replacing the colored edges by the edge-quotients of the atoms, the colored loops by the
loop-quotients and the color half-edges by some choices of half-quotients. This gives the complete
structural description of all graphs which can be reached from $H_r$ by expansion. Half-edges of
$H_{i-1}$ can arise only in expansions of half-edges of $H_i$.

In Section~\ref{sec:algorithm}, we describe the meta-algorithm itself. From algorithmic point of
view, the key difficulty arises from the fact that a graph $H_r$ can have exponentially many
pairwise non-isomorphic expansions $H_0$. Therefore we cannot test all of them and we proceed in the
opposite way. We start with the graph $H$ and try to reach $H_r$ by a series of reductions. But here
the reductions are non-deterministic since a part of the graph $H$ can correspond to many
different subgraphs of $G$. Therefore, we keep lists of possible correspondences when we replace
atoms in $H$.  We proceed with the reductions and compute further lists using dynamic programming. There
is only one slow subroutine which takes time $\O^*(2^{e(H)/2})$ which we describe in detail in
Section~\ref{sec:comb_int_star_atoms}.

In Section~\ref{sec:planar_graphs} we deal with specific properties of planar graphs and show that
the meta-algorithm applies to them. It is a key observation that the \rcover\ problem is trivially
solvable for 3-connected inputs $G$ since the automorphism group $\Aut(G)$ is spherical; it is
either cyclic, dihedral or a subgroup of one of the three special groups.  One can just enumerate
all quotients of $G$ and test graph isomorphism to $H$. The reason is that 3-connected planar graphs
and their quotients behave geometrically. On the other hand, dealing with general planar
graphs is rather non-trivial, since the geometry of the sphere is lost and it requires all the theory
built in this paper. We establish that planar graphs satisfy the conditions (P0) to (P3) where
especially the proof for (P3) is not straightforward. 

\section{Structural Properties of Atoms} \label{sec:atoms}

In this section, we introduce special inclusion-minimal subgraphs of $G$ called atoms. We show their
structural properties such as that they behave nicely with respect to any covering projection.

\subsection{Block-trees and Their Automorphisms} \label{sec:block_trees}

\begin{figure}[b!]
\centering
\includegraphics{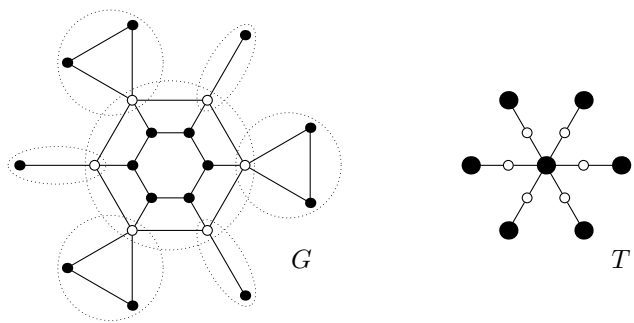}
\caption{On the left, an example graph $G$ with denoted blocks. On the right, the corresponding
block-tree $T$ is depicted. The white vertices correspond to the articulations and the big black
vertices correspond to the blocks.}
\label{fig:example_of_block_tree}
\end{figure}

The \emph{block-tree} $T$ of $G$ is a tree defined as follows. Consider all articulations in $G$
and all maximal $2$-connected subgraphs which we call \emph{blocks} (with bridge-edges also counted
as blocks).  The block-tree $T$ is the incidence graph between the articulations and the blocks. For
an example, see Figure~\ref{fig:example_of_block_tree}.

There is the following well-known connection between $\Aut(G)$ and $\Aut(T)$:

\begin{lemma}
Every automorphism $\pi \in \Aut(G)$ induces an automorphisms $\pi' \in \Aut(T)$.
\end{lemma}

\begin{proof}
First, observe that every automorphism $\pi$ of $G$ maps the articulations to the articulations and
the blocks to the blocks which gives the induced mapping $\pi'$. It remains to show that $\pi'$ is
an automorphism of $T$. Let $a$ be an articulation adjacent to a block $B$ in the tree. Then $a$ is
contained in $B$. Therefore $\pi'(a)$ is contained in $\pi'(B)$ and vice versa, which implies that
$\pi'$ is an automorphisms of the block-tree $T$.\qed
\end{proof}

We note that there is no direct relation between the structure of $\Aut(G)$ and $\Aut(T)$. First,
$\Aut(T)$ may contain some additional automorphisms not induced by anything in $\Aut(G)$. Second,
several distinct automorphisms of $\Aut(G)$ may induce the same automorphism of $\Aut(T)$. For
example in Figure~\ref{fig:example_of_block_tree}, $\Aut(G) \cong \gD_3 \times \gC_2^3$ and $\Aut(T)
\cong \gS_6$.

\heading{The Central Block.} The \emph{center} of a graph is a subset of the vertices which
minimize the maximum distance to all vertices of the graph. For a tree, its center is either the
central vertex or the central pair of vertices of a longest path, depending on the parity of its
length. Every automorphism of a graph preserves its center.

\begin{lemma} \label{lem:central_block}
If $G$ has a non-trivial semiregular automorphism, then $G$ has a central block.
\end{lemma}

\begin{proof}
For the block-tree $T$, all leaves are blocks, so each longest path is of an even length. Therefore
$\Aut(T)$ preserves the central vertex. The central vertex can be either a \emph{central
articulation}, or a \emph{central block}. If the central vertex is an articulation $u$, then every
automorphism fixes $u$ which contradicts the assumptions.\qed
\end{proof}

In the following, we shall assume that $T$ contains a central block.  We orient edges of the
block-tree $T$ towards the central block; so the block-tree becomes rooted. A \emph{subtree} of the
block-tree is defined by any vertex different from the central block acting as \emph{root} and by
all its descendants.

Let $u$ be an articulation contained in the central block. By $R_u$ we denote the subtree attached
to the central block at $u$.

\begin{lemma} \label{lem:semiregular_action_on_blocks}
Let $\Gamma$ be a semiregular subgroup of $\Aut(G)$. If $u$ and $v$ are two articulations of the
central block and of the same orbit of $\Gamma$, then $R_u \cong R_v$. Moreover there is a unique
$\pi \in \Gamma$ which maps $R_u$ to $R_v$.
\end{lemma}

\begin{proof}
Notice that either $R_u = R_v$, or $R_u \cap R_v = \emptyset$. Since $u$ and $v$ are in the same
orbit of $\Gamma$, there exists $\pi \in \Gamma$ such that $\pi(u) = v$. Consequently $\pi(R_u)
= R_v$.  Suppose that there exist $\pi,\sigma \in \Gamma$ such that $\pi(R_u) = \sigma(R_u) = R_v$.
Then $\pi\cdot\sigma^{-1}$ is an automorphism of $\Gamma$ fixing $u$. Since $\Gamma$ is semiregular,
$\pi = \sigma$.\qed
\end{proof}

\subsection{Definition and Basic Properties of Atoms}

Let $u$ and $v$ be two distinct vertices of degree at least three joined by at least two parallel
edges. Then the subgraph induced by $u$ and $v$ is called a \emph{dipole}.  Let $B$ be one block of
$G$, so $B$ is a 2-connected graph. Two vertices $u$ and $v$ form a \emph{2-cut} $U = \{u,v\}$ if $B
\setminus U$ is disconnected. We say that a 2-cut $U$ is \emph{non-trivial} if $\deg(u) \ge 3$ and
$\deg(v) \ge 3$.

\begin{lemma} \label{lem:minimal_cut_property}
Let $U$ be a $2$-cut and let $C$ be a component of $B \setminus U$. Then $U$ is uniquely determined
by $C$.
\end{lemma}

\begin{proof}
If $C$ is a component of $B \setminus U$, then $U$ has to be the set of all neighbors of $C$ in $B$.
Otherwise $B$ would not be 2-connected, or $C$ would not be a component of $B \setminus U$.\qed
\end{proof}

\heading{The Definition.} We first define a set $\calP$ of subgraphs of $G$ which we call
\emph{parts}:
\begin{packed_itemize}
\item A \emph{block part} is a subgraph non-isomorphic to $K_2$ induced by the blocks of a subtree
of the block-tree. 
\item A \emph{proper part} is a subgraph $S$ of $G$ defined by a non-trivial 2-cut $U$ of a block
$B$ not containing the central block. The subgraph $S$ consists of a connected component $C$ of $G
\setminus U$ together with $u$ and $v$ and all edges between $\{u,v\}$ and $C$.
\item A \emph{dipole part} is any dipole.
\end{packed_itemize}

\begin{figure}[t!]
\centering
\includegraphics{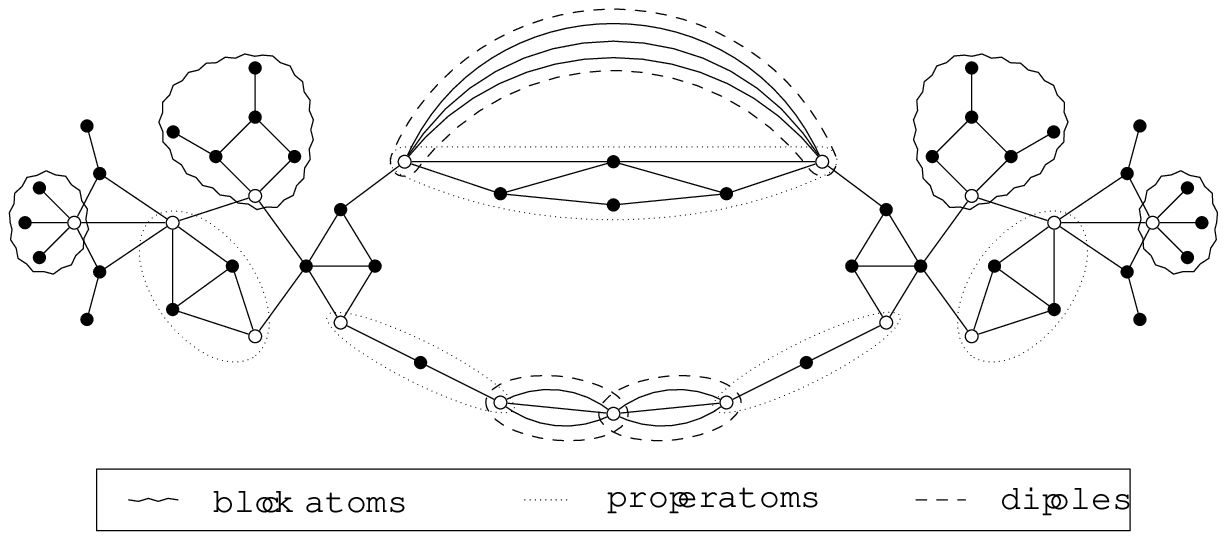}
\caption{An example of a graph with denoted atoms. The white vertices belong to the boundery of some
atom, possibly several of them.}
\label{fig:atoms_examples}
\end{figure}

\noindent The inclusion-minimal elements of $\calP$ are called \emph{atoms}. We distinguish \emph{block
atoms}, \emph{proper atoms} and \emph{dipoles} according to the type of the defining part.  Block
atoms are either pendant stars, or pendant blocks possibly with single pendant edges attached to it.
Also each proper atom or dipole is a subgraph of a block. For an example, see
Figure~\ref{fig:atoms_examples}. 

We use topological notation to denote the \emph{boundary} $\bo A$  and the \emph{interior} $\int A$
of an atom $A$.  If $A$ is a dipole, we set $\bo A = V(A)$. If $A$ is a proper or block atom, we put
$\bo A$ equal the set of vertices of $A$ which are incident with an edge not contained in $A$. For
the interior, we use the standard topological definition $\int A = A \setminus \bo A$ where we only
remove the vertices $\bo A$, the edges adjacent to $\bo A$ are kept.

Note that $|\bo A| = 1$ for any block atom $A$, and $|\bo A| = 2$ for a proper atom or dipole $A$.
The interior of a dipole is a set of free edges.  We note that dipoles are automatically atoms and
they are exactly the atoms with no vertices in their interiors. Observe for a proper atom $A$ that
the vertices of $\bo A$ are exactly the vertices $\{u,v\}$ of the non-trivial 2-cut used in the
definition of proper parts. Also the vertices of $\bo A$ of a proper atom are never adjacent.
Further, no block or proper atom contains parallel edges; otherwise a dipole would be its subgraph.

\heading{Properties of Atoms.}
Our goal is to replace atoms by edges, and so it is important to know that the atoms cannot overlap
too much. The reader can see in Figure~\ref{fig:atoms_examples} that the atoms only share their
boundaries. This is true in general, and we are going to prove it in two steps now. 

\begin{lemma} \label{lem:nonintersecting_interiors}
The interiors of atoms are pairwise disjoint.
\end{lemma}

\begin{proof}
For contradiction, let $A$ and $A'$ be two distinct atoms with non-empty intersections of $\int A$
and $\int A'$. First suppose that $A$ is a block item. Then $A$ corresponds to a subtree of the
block-tree which is attached by an articulation $u$ to the rest of the graph. If $A'$ is a block
atom then it corresponds to some subtree, and we can derive that $A \subseteq A'$ or $A' \subseteq
A$. And if $A'$ is proper atom or dipole, then it is a subgraph of a block, and thus subgraph of
$A$. In both cases, we get contradiction with minimality. Similarly, if one atom is a dipole, we can
easily argue contradiction with minimality.

The last case is that both $A$ and $A'$ are proper atoms. Since the interiors are connected and the
boundaries are defined as neighbors of the interiors, it follows that both $W' = A \cap \bo A'$
and $W = A' \cap \bo A$ are nonempty. We have two cases according to the sizes of these
intersections depicted in Figure~\ref{fig:intersecting_interiors}.

\begin{figure}[t!]
\centering
\includegraphics{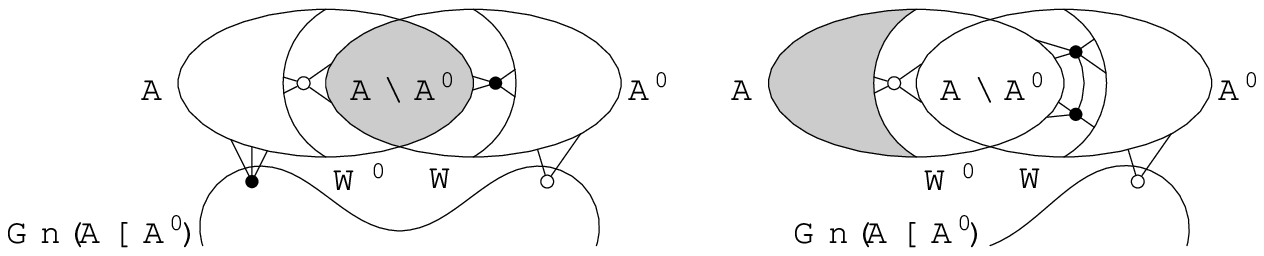}
\caption{We depict the vertices of $\bo A$ in black and the vertices of $\bo A'$ in white.
In both cases, we found a subset of $A$ belonging to $\calP$ (its interior is highlighted in gray).}
\label{fig:intersecting_interiors}
\end{figure}

If $|W| = |W'| = 1$, then $W \cup W'$ is a 2-cut separating $\int A \cap \int A'$ which
contradicts minimality of $A$ and $A'$. And if, without loss of generality, $|W| = 2$, then there is
no edge between $\int A \setminus (\int A' \cup W')$ and the remainder of the graph $G
\setminus (\int A \cup \int A')$. Therefore, $\int A \setminus (\int A' \cup
W')$ is separated by a 2-cut $W'$ which again contradicts minimality of $A$. We note that in both
cases the constructed 2-cut is non-trivial since it is formed by vertices of non-trivial cuts
$\bo A$ and $\bo A'$.\qed
\end{proof}

Next we show a stronger version of the previous lemma which states that two atoms can intersect only
in their boundaries.

\begin{lemma} \label{lem:nonintersecting_atoms}
Let $A$ and $A'$ be two atoms. Then $A \cap A' = \bo A \cap \bo A'$.
\end{lemma}

\begin{proof}
We already know from Lemma~\ref{lem:nonintersecting_interiors} that $\int A \cap \int A'=
\emptyset$. It remains to argue that, say, $\bo A \cap \int A' = \emptyset$. If $A$ is a
block atom, then $\bo A$ is the articulation separating $A$. No atom can contain this
articulation as its interior. Similarly, if $A'$ is a block atom, then $A$ has to be contained
in the interior of $A'$ or vice versa which contradicts minimality.

Let $\bo A = \{u,v\}$ and $\bo A' = \{u',v'\}$.  First we deal with dipoles. The situation
where $A'$ is a dipole is trivial. And if $A$ is a dipole with $u \in \int A'$, then either $v
\in A$ which contradicts minimality of $A$, or $\bo A'$ is not a 2-cut.  It remains to deal
with both $A$ and $A'$ being proper atoms. Recall that in such a case $\bo A$ is defined as
neighbors of $\int A$ in $G$, and that $\bo A'$ are neighbors of $\int A'$ in $G$.

\begin{figure}[b!]
\centering
\includegraphics{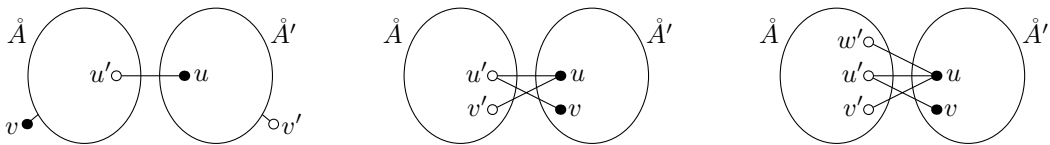}
\caption{An illustration of the main steps of the proof.}
\label{fig:intersecting_atoms}
\end{figure}

The proof is illustrated in Figure~\ref{fig:intersecting_atoms}.  Suppose for contradiction that
$\int A \cap \bo A' \ne \emptyset$ and let $u' \in \int A$. Since $u'$ has at least
one neighbor in $\int A'$, then without loss of generality $u \in \int A'$ and $uu' \in E(G)$.
Since $A$ is a proper atom, the set $\{u',v\}$ is not a 2-cut, so there is another neighbor of $u$ in
$\int A$, which has to be equal $v'$. Symmetrically, $u'$ has another neighbor in $\int
A'$ which is $v$.  So $\bo A \subseteq \int A'$ and $\bo A' \subseteq \int A$.
If $\bo A = \int A'$ and $\bo A' = \int A$, the graph is $K_4$ (since the
minimal degree of cut-vertices is three) which contradicts existence of 2-cuts and atoms. And if for
example $\int A \ne \bo A'$, then $\bo A'$ does not cut a subset of $\int A$, so
there has to a third neighbor $w'$ of $\int A'$, which contradicts that $\bo A'$ cuts
$\int A'$ from the rest of the graph.\qed
\end{proof}

\heading{Connectivity of Atoms.}
We call a graph \emph{essentially 3-connected} if it is a 3-connected graphs with possibly single
pendant edges attached to it. For instance, every block atom is essentially 3-connected. A proper
$A$ might not be essentially 3-connected. Let $\bo A = \{u,v\}$. We define $A^+$ as $A$ with the
additional edge $uv$. Notice that the property (P0) ensures that $A^+$ belongs to $\calC$. It is
easy to see that $A^+$ is essentially 3-connected graph. Additionally, we put $A^+ = A$ for a block
atom or dipole.

\begin{lemma} \label{lem:aut_ess_3-conn}
Let $A$ be an essentially 3-connected graph, and we construct $B$ from $A$ by removing the single
pendant edges of $A$. Then $\Aut(A)$ is a subgroup of $\Aut(B)$. 
\end{lemma}

\begin{proof}
These pendant single edges behave like markers, giving a 2-partition of $V(G)$ which $\Aut(A)$ has
to preserve.\qed
\end{proof}

Further, if $\Aut(B)$ is of polynomial size, we can easily check which permutations preserve this
2-partition, and thus give $\Aut(A)$. Also, similar relation holds for any group $\Gamma$ acting on
$B$ and its subgroup $\Gamma'$ preserving the 2-partition.

It is important that we can code the 2-partition by coloring the vertices of $B$, and work with such
colored 3-connected graph using (P2) and (P3).

\subsection{Symmetry Types of Atoms}

We distinguish three symmetry types of atoms which describe how symmetric each atom is. When such an
atom is reduced, we replace it by an edge carrying the type. Therefore we have to use multigraphs
with three edge types: \emph{halvable edges}, \emph{undirected edges} and \emph{directed edges}. We
consider only the automorphisms which preserve these edge types and indeed the orientation of
directed edges.

Let $A$ be a proper atom or dipole with $\bo A = \{u,v\}$. Then we distinguish the following
three symmetry types, see Figure~\ref{fig:types_of_atoms}:

\begin{figure}[b!]
\centering
\includegraphics{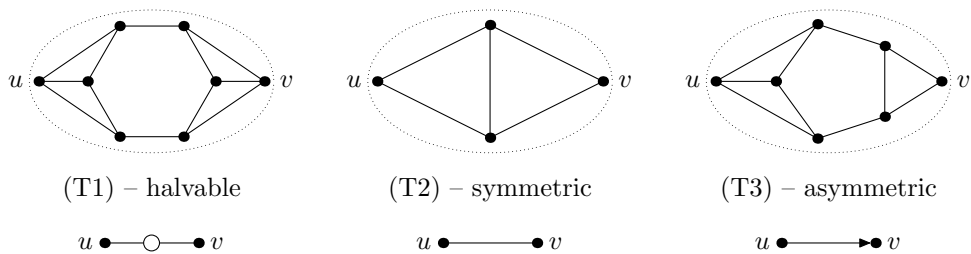}
\caption{The three types of atoms and the corresponding edge types which we use in the reduction.}
\label{fig:types_of_atoms}
\end{figure}

\begin{packed_itemize}
\item \emph{The halvable atom.} There exits an semiregular involutory automorphism $\pi$
which exchanges $u$ and $v$. More precisely, the automorphism $\pi$ fixes no vertices and no edges
with an exception of some halvable edges.
\item \emph{The symmetric atom.} The atom is not halvable, but there exists an automorphism
which exchanges $u$ and $v$.
\item \emph{The asymmetric atom.} The atom which is neither halvable nor symmetric.
\end{packed_itemize}
If $A$ is a block atom, then it is by definition symmetric.

\begin{lemma} \label{lem:dipole_type}
For a given dipole $A$, it is possible to determine its type in polynomial time.
\end{lemma}

\begin{proof}
The type depends only on the quantity of distinguished types of the
parallel edges.  We have directed edges from $u$ to $v$, directed edges from $v$ to $u$, undirected
edges and halvable edges. We call a dipole \emph{balanced} if the number of directed edges in the
both directions is the same. Observe that:
\begin{packed_itemize}
\item The dipole $A$ is halvable if and only if it is balanced and has an even number of undirected
edges.
\item The dipole $A$ is symmetric if and only if it is balanced and has an odd number of undirected
edges.
\item The dipole $A$ is asymmetric if and only if it is unbalanced.
\end{packed_itemize}
This clearly can be tested in polynomial time.\qed
\end{proof}

\begin{lemma} \label{lem:atom_type_algorithm}
For a given proper atom $A$ of $\calC$ satisfying (P2), it is possible to determine its type in polynomial
time.
\end{lemma}

\begin{proof}
Let $\bo A = \{u,v\}$. Recall that $A^+$ is an essentially 3-connected graph. Let $B$ be the
3-connected graph created from $A^+$ by removing pendant edges, where existence of pendant edges is
coded by colors of $V(B)$.  Using (P3), we can check whether there is a color-preserving
automorphism exchanging $u$ to $v$ as follows, see Figure~\ref{fig:testing_symmetry}. We take two
copies of $B$. In one copy, we color $u$ by a special color, and $v$ by another special color. In
the other copy, we swap the colors of $u$ and $v$. Using (P3) on these two copies, we can check
whether there is an automorphism which exchanges $u$ and $v$. If not, then $A$ is asymmetric. If
yes, we check whether $A$ is symmetric or halvable.

\begin{figure}[b!]
\centering
\includegraphics{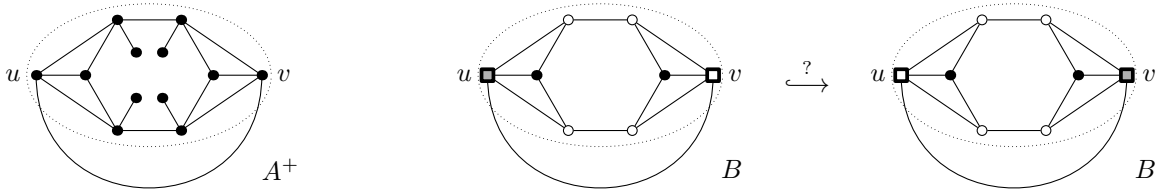}
\caption{For the depicted atom $A$, we test using (P3) whether $B \hookto B$. In this case yes, so
$A$ is either symmetric, or halvable.}
\label{fig:testing_symmetry}
\end{figure}

Using (P2), we generate polynomially many semiregular involutions of order two acting on $B$.  For
each semiregular involution, we check whether it transposes $u$ to $v$, and whether it preserves the
colors of $V(B)$ coding pendant edges. If such a semiregular involution exists, then $A$ is halvable,
otherwise it is just symmetric.\qed
\end{proof}

\subsection{Automorphisms of Atoms} \label{sec:automorphisms_of_atoms}

We start with a simple lemma which states how automorphisms behave with respect to atoms.

\begin{lemma} \label{lem:atom_automorphisms}
Let $A$ be an atom and let $\pi \in \Aut(G)$. Then the following
holds:
\begin{packed_enum}
\item[(a)] The image $\pi(A)$ is an atom isomorphic to $A$. Further $\pi(\bo A) = \bo
\pi(A)$ and $\pi(\int A) = \int \pi(A)$.
\item[(b)] If $\pi(A) \ne A$, then $\pi(\int A) \cap \int A = \emptyset$.
\item[(c)] If $\pi(A) \ne A$, then $\pi(A) \cap A = \bo A \cap \bo \pi(A)$.
\end{packed_enum}
\end{lemma}

\begin{proof}
(a) Every automorphism permutes the set of articulations and non-trivial 2-cuts. So $\pi(\bo A)$
separates $\pi(\int A)$ from the rest of the graph. It follows that $\pi(A)$ is an
atom, since otherwise $A$ would not be an atom. And $\pi$ clearly preserves the boundaries and the
interiors.

For the rest, (b) follows from Lemma~\ref{lem:nonintersecting_interiors} and (c) follows from
Lemma~\ref{lem:nonintersecting_atoms}.\qed
\end{proof}

Therefore, for an automorphism $\pi$ of an atom $A$, we require that $\pi(\bo A) = \bo \pi(A)$. If
a block or proper atom $A \in \calC$ satisfying (P2), then we can compute $\Aut(A)$ according to
Lemma~\ref{lem:aut_ess_3-conn} in polynomial time.

\begin{figure}[t!]
\centering
\includegraphics{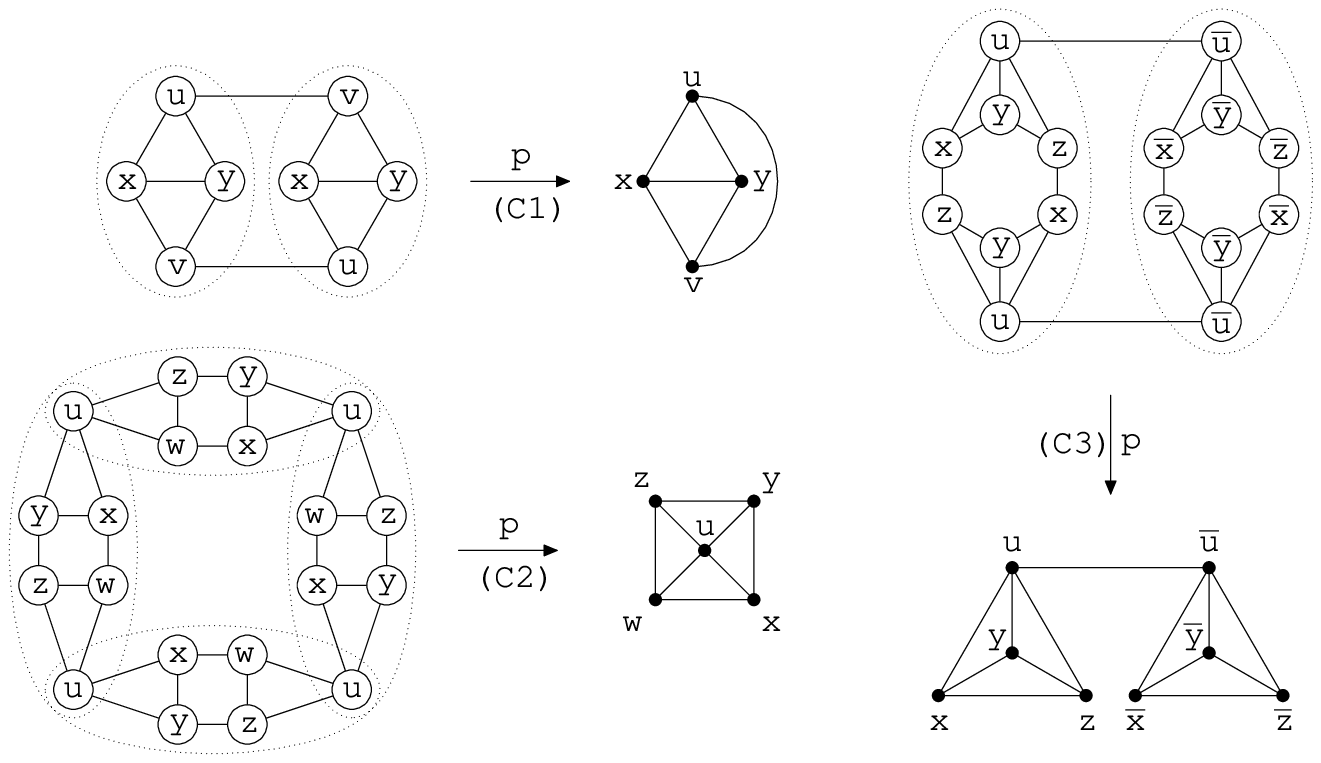}
\caption{The three cases for mapping of atoms (depicted in dots). Notice that for the third
graph, a projection of the type (C1) could also be applied which would give a different quotient.}
\label{fig:atom_projection_cases}
\end{figure}

\heading{Projections of Atoms.}
Let $\Gamma$ be a semiregular subgroup of $\Aut(G)$, which defines a regular covering projection $p
: G \to G / \Gamma$. Negami~\cite[p.~166]{negami} investigated possible projections of proper atoms,
and we investigate this question in more details.  For a proper atom or a dipole $A$ with $\bo
A = \{u,v\}$, we get one of the following three cases illustrated in
Figure~\ref{fig:atom_projection_cases}.

\begin{packed_head_enum}{(C1)}
\item[(C1)] The atom $A$ is preserved in $G / \Gamma$, meaning $p(A) \cong A$. Notice that $p(A)$ is
just a subgraph of $G / \Gamma$. For a proper atom, it can happen that $p(u)p(v)$ is adjacent, even
through $uv \notin E(G)$, as in Figure~\ref{fig:atom_projection_cases}.
\item[(C2)] The interior of the atom $A$ is preserved and the vertices $u$ and $v$ are
identified, i.e., $p(\int A) \cong \int A$ and $p(u) = p(v)$.
\item[(C3)] The covering projection $p$ is a $2k$-fold cover. There exists an involutory permutation
$\pi$ in $\Gamma$ which exchanges $u$ and $v$ and preserves $A$. The projection $p(A)$ is a halved
atom $A$. This can happen only when $A$ is a halvable atom.
\end{packed_head_enum}

\begin{lemma} \label{lem:atom_covering_cases}
For every atom $A$ and every semiregular subgroup $\Gamma$ defining covering projection $p$, one of
the cases (C1), (C2) and (C3) happens. Moreover, for a block atom we have exclusively the case (C1).
\end{lemma}

\begin{proof}
For a block atom $A$, Lemma~\ref{lem:semiregular_action_on_blocks} implies that $p(A) \cong A$, so
the case (C1) happens. It remains to deal with $A$ being a proper atom or a dipole, and let
$\bo A = \{u,v\}$. According to Lemma~\ref{lem:atom_automorphisms}b every automorphism $\pi$
either preserves $\int A$, or $\int A$ and $\pi(\int A)$ are disjoint. If there
exists a non-trivial $\pi \in \Gamma$ which preserves $\int A$, we get (C3); otherwise we get
(C1) or (C2).

Let $\pi$ be a non-trivial automorphism of $\Gamma$ preserving $\int A$. We know $\pi(\bo A) = \bo
A$ and by semiregularity, $\pi$ has to exchange $u$ and $v$.  Then the fiber containing $u$ and $v$
has to be of an even size, with $\pi$ being an involution reflecting $k$ copies of $A$, and
therefore the covering $p$ is a $2k$-fold cover. This proves (C3).

Suppose there is no non-trivial automorphism which preserves $\int A$. The only difference between
(C1) and (C2) is whether $u$ and $v$ are contained in one fiber of $p$, or not. First suppose that
for every non-trivial $\pi \in \Gamma$ we get $A \cap \pi(A) = \emptyset$. Then no fiber contains
more than one vertex of $A$, and we get (C1), i.e, $A \cong p(A)$.  And if there exists $\pi \in
\Gamma$ such that $A \cap \pi(A) \ne \emptyset$. By Lemma~\ref{lem:atom_automorphisms}c, we get
$A \cap \pi(A) = \bo A \cap \bo \pi(A)$, so $u$ and $v$ belong to one fiber of $p$, which gives
(C2).\qed
\end{proof}

\section{Graph Reductions and Quotient Expansions} \label{sec:reduction_and_expansion}

We start with a quick overview. The reduction initiates with a graph $G$ and produces a sequence of
graphs $G = G_0,G_1,\dots,G_r$. To produce $G_i$ from $G_{i-1}$, we find a collection of all atoms
$\calA$ and replace each of them by an edge of the corresponding type. We stop at step $r$ when
$G_r$ contains no further atoms, and we call such a graph \emph{primitive}. We call this sequence of
graphs starting with $G$ and ending with a primitive graph $G_r$ as the \emph{reduction series} of
$G$.

Now suppose that $H_r = G_r / \Gamma_r$ is some quotient of $G_r$. To revert the reductions applied
to obtain $H_r$, we revert the reduction series on $H_r$ and produce an \emph{expansion series}
$H_r,H_{r-1},\dots,H_0$ of $H_r$. We obtain semiregular subgroups $\Gamma_0,\dots,\Gamma_r$ such
that $H_i = G_i / \Gamma_i$. The entire process is depicted in the following diagram:

\begin{equation} \label{eq:red_and_exp_diagram}
\begin{gathered}
\xymatrix{
G_0 \ar[d]_{\Gamma_0} \ar[r]^{\rm red.} & G_1 \ar[d]_{\Gamma_1} \ar[r]^{\rm red.} &
\cdots \ar[r]^{\rm red.}&G_i \ar[d]_{\Gamma_i} \ar[r]^{\rm red.} &
G_{i+1}\ar[d]_{\Gamma_{i+1}} \ar[r]^{\rm red.} & \cdots \ar[r]^{\rm red.} &
G_r \ar[d]_{\Gamma_r}\\%
H_0 									& H_1 \ar[l]_{\rm exp.} &
\cdots \ar[l]_{\rm exp.}&H_i \ar[l]_{\rm exp.} &
H_{i+1}\ar[l]_{\rm exp.}					   & \cdots \ar[l]_{\rm exp.} &
H_r \ar[l]_{\rm exp.}\\%
}
\end{gathered}
\end{equation}

In this section, we describe structural properties of reductions and expansions. We study changes of
automorphism groups done by reductions. Indeed, $\Aut(G_{i+1})$ can differ from $\Aut(G_i)$.
But the reduction is done right and the important information of $\Aut(G_i)$ is preserved in
$\Aut(G_{i+1})$ which is key for expansions. The issue is that expansions are unlike reductions not
uniquely determined. From $H_{i+1}$, we can construct multiple $H_i$. In this section, we
characterize all possible ways how $H_i$ can be constructed from $H_{i+1}$. 

\subsection{Reducing Graphs Using Atoms}

The reduction produces a series of graphs $G = G_0,\dots,G_r$. To construct $G_i$ from $G_{i-1}$, we
find the collection of all atoms $\calA$ and determine their types, using
Lemma~\ref{lem:atom_type_algorithm}. We replace a block atom $A$ by a pendant edge of some color
based at $u$ where $\bo A = \{u\}$. We replace each proper atom or dipole $A$ with $\bo A = \{u,v\}$
by a new edge $uv$ of some color and of one of the three edge types according to the type of $A$.
According to Lemma~\ref{lem:nonintersecting_atoms}, the replaced parts for the atoms of $\calA$ are
pairwise disjoint, so the reduction is well defined. We stop in the step $r$ when $G_r$ contains no
atoms. We show in Lemma~\ref{lem:primitive_graphs} that a primitive graph is either $3$-connected, a
cycle, or $K_2$ possibly with attached single pendant edges.

To be more precise, we consider graphs with colored vertices, colored edges and with three edge
types. We say that two graphs $G$ and $G'$ are isomorphic if there exists an isomorphism which
preserves all colors and edge types, and we denote this by $G \cong G'$. We note that the results
built in Section~\ref{sec:atoms} transfers to colored graphs and colored atoms without any problems.
Two atoms $A$ and $A'$ are isomorphic if there exists an isomorphism which maps $\bo A$ to
$\bo A'$. We obtain isomorphism classes for the set of all atoms $\calA$ such that $A$ and
$A'$ belong to the same class if and only if $A \cong A'$. To each color class, we assign one new
color not yet used in the graph. When we replace the atoms of $\calA$ by edges, we color the edges
according to the colors assigned to the isomorphism classes.

\begin{figure}[t!]
\centering
\includegraphics{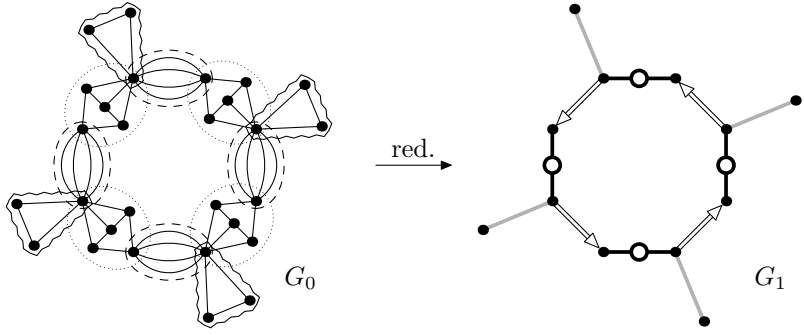}
\caption{On the left, we have a graph $G_0$ with three isomorphism classes of atoms, each having
four atoms.  The dipole atoms are halvable, the block atoms are symmetric and the proper atoms are
asymmetric.  We reduce $G_0$ to $G_1$ which is an eight cycle with pendant leaves, with four
black halvable edges, four gray undirected edges, and four white directed edges. The reduction
series ends with $G_1$ since it is primitive.}
\label{fig:example_of_reduction}
\end{figure}

It remains to say that for an asymmetric atom we choose an arbitrary orientation, but consistently
with $\Aut(G_{i-1})$ for the entire isomorphism class. For an example of the reduction, see
Figure~\ref{fig:example_of_reduction}.

\begin{figure}[t!]
\centering
\includegraphics{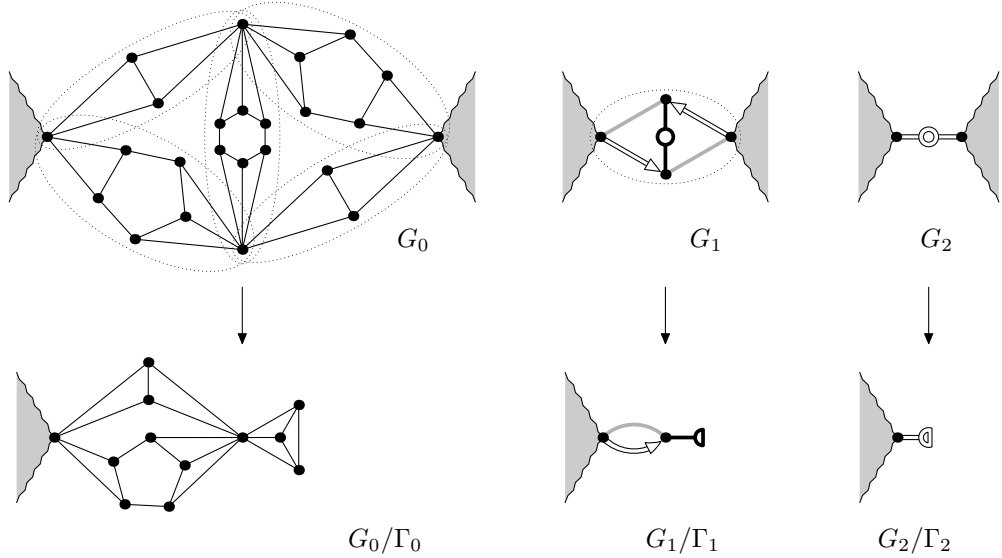}
\caption{We reduce a part of a graph in two steps. In the first step, we replace five atoms by five
edges of different types. As the result we obtain one halvable atom which we further reduce to one
halvable edge. Notice that without considering edge types, the resulting atom would be just
symmetric. In the bottom we show the corresponding quotient graphs when $\Gamma_i$ are generated by
semiregular involutory automorphism $\pi$ from the definition of the halvable atom.}
\label{fig:reducing_atoms}
\end{figure}

The symmetry type of atoms depends on the types of edges the atom contains; see
Figure~\ref{fig:reducing_atoms} for an example. Also, the figure depicts a quotient $G_2 / \Gamma_2$
of $G_2$, and its expansions to $G_1 / \Gamma_1$ and $G_0 / \Gamma_0$. The resulting quotients $G_1
/ \Gamma_1$ and $G_2 / \Gamma_2$ contain half-edges because $\Gamma_1$ and $\Gamma_2$ fixes some
halvable edges but $G_0 / \Gamma_0$ contains no half-edges. This example shows that for reductions
and expansions we need to consider half-edges even when the input $G$ and $H$ are simple graphs.

\heading{Properties of Reduction.}
Consider the groups $\Aut(G_i)$ and $\Aut(G_{i+1})$. There exists a natural homomorphism $\Phi_i:
\Aut(G_i) \to \Aut(G_{i+1})$ which we define as follows. Let $\pi \in \Aut(G_i)$. The graph
$G_{i+1}$ is constructed from $G_i$ by replacing interiors of all atoms by colored edges. For the
common vertices and edges of $G_i$ and $G_{i+1}$, we define the image in $\Phi_i(\pi)$ exactly as in
$\pi$. If $A$ is an atom of $G_i$, then according to Lemma~\ref{lem:atom_automorphisms}a, $\pi(A)$
is an atom isomorphic to $A$. In $G_{i+1}$, we replace the interiors of both $A$ and $\pi(A)$ by the
edges $e_A$ and $e_{\pi(A)}$ of the same type and color. Therefore, we define $\Phi_i(\pi)(e_A) =
e_{\pi(A)}$.

More precisely for purpose of Section~\ref{sec:quotients}, we define $\Phi_i$ on the half edges. Let
$e_A = uv$ and let $h_u$ and $h_v$ be the half-edges composing $e_A$, and similarly let $h_{\pi(u)}$
and $h_{\pi(v)}$ be the half-edges composing $e_{\pi(A)}$. Then we define $\Phi_i(\pi)(h_u) =
h_{\pi(u)}$ and $\Phi_i(\pi)(h_v) = h_{\pi(v)}$. 

\begin{proposition} \label{prop:reduction_homomorphism}
The mapping $\Phi_i$ satisfies the following:
\begin{packed_enum}
\item[(a)] The mapping $\Phi_i$ is a group homomorphism.
\item[(b)] The mapping $\Phi_i$ is surjective.
\item[(c)] Moreover, $\Aut(G_{i+1}) = \Phi_i(\Aut(G_i))$ monomorphically embeds into $\Aut(G_i)$. 
\item[(d)] For a semiregular subgroup $\Gamma$ of $\Aut(G_i)$, the mapping $\Phi_i |_\Gamma$ is an
isomorphism. Moreover, the subgroup $\Phi_i(\Gamma)$ remains semiregular.
\end{packed_enum}
\end{proposition}

\begin{proof}
(a) It is easy to see that each $\Phi_i(\pi) \in \Aut(G_{i+1})$.  The well-known Homomorphism
Theorem states that $\Phi_i$ is a homomorphism if and only if the kernel $\Ker(\Phi_i)$, i.e., the set
of all $\pi$ such that $\Phi_i(\pi) = \id$, is a normal subgroup of $\Aut(G_i)$. It is easy to see
that the kernel $\Ker(\Phi_i)$ has the following structure. If $\pi \in \Ker(\Phi_i)$, it fixes
everything except for the interiors of the atoms. Further, $\pi(\int A) = \int \pi(A)$, so
$\pi$ can non-trivially act only inside the interiors of the atoms.

Let $\sigma \in \Aut(G_i)$ and $\pi \in \Ker(\Phi_i)$. We need to show that $\sigma \pi \sigma^{-1}
\in \Ker(\Phi_i)$. Let $A$ be an atom. Then $\sigma(A) = A'$ is an isomorphic atom. The composition
clearly permutes the interior $\int A$. Moreover, the part of the graph outside of interiors is fixed
by the composition. Hence it belongs to $\Ker(\Phi_i)$ and by the Homomorphism Theorem $\Phi_i$ is a
homomorphism.

(b) Let $\pi' \in \Aut(G_{i+1})$, we want to extend $\pi'$ to $\pi \in \Aut(G_i)$ such that
$\Phi_i(\pi) = \pi'$. We just need to describe this extension on a single edge $e = uv$. If $e$ is
an original edge of $G$, there is nothing to extend. Suppose that $e$ was created in $G_{i+1}$ from
an atom $A$ in $G_i$. Then $e' = \pi'(e)$ is an edge of the same color and the same type as $e$, and
therefore $e'$ is constructed from an isomorphic atom $A'$ of the same symmetry type. The
automorphism $\pi'$ prescribes the action on the boundary $\bo A$. We need to show that it is
possible to define an action on $\int A$ consistently.

If $A$ is a block atom, then the both edges $e$ and $e'$ are pendant, attached by articulations $u$
and $u'$. We just define $\pi$ using an isomorphism from $A$ to $A'$ which takes $u$ to $u'$.
It remains to deal with proper atoms and dipoles.

First suppose that $A$ is an asymmetric atom. Then by definition the orientation of $e$ and $e'$ is
consistent with respect to $\pi'$. Since $\int A \cong \int A'$, we define $\pi$ on
$\int A$ according to one such isomorphism.

Secondly suppose that $A$ is symmetric or halvable. Let $\sigma$ be an isomorphism of $A$
and $A'$. Either $\sigma$ maps $\bo A$ exactly as $\pi'$, and then we can use $\sigma$ for
defining $\pi$. Or we compose $\sigma$ with the automorphism of $A$ exchanging the two vertices of
$\bo A$. (We know that such an automorphism exists since $A$ is not antisymmetric.)

(c) Let $(A_1,A_2,\dots,A_\ell)$ be an orbit of $\pi$ on edges representing atoms
$A_1,\dots,A_\ell$. We construct the extension as above, choosing any isomorphism from $A_i$ to
$A_{i+1}$, where $i = 1,\dots,\ell-1$, and properly the isomorphism $A_\ell$ to $A_1$. Here,
properly means that composition of all these isomorphisms is the identity on $A_1$.  Repeating this
procedure for every orbit of $\pi$, we determine an extension $\tilde\pi$ of $\pi$ defining a
monomorphic embedding $\pi \mapsto \tilde\pi$.

(d) We note that for any subgroup $\Gamma$, the restricted mapping $\Phi_i |_\Gamma$ is a group
homomorphism with $\Ker(\Phi_i |_\Gamma) = \Ker(\Phi_i) \cap \Gamma$. If $\Gamma$ is semiregular,
then $\Ker(\Phi_i) \cap \Gamma$ is trivial. The reason is that $G_i$ contains at least one atom $A$,
and the boundary $\bo A$ is fixed by $\Ker(\Phi_i)$. Hence $\Phi_i |_\Gamma$ is an isomorphism.

For the semiregularity of $\Phi_i(\Gamma)$, let $\pi'$ be an automorphism of $G_{i+1}$. Since
$\Phi_i|_\Gamma$ is an isomorphism, there exists the unique $\pi \in \Gamma$ such that $\Phi_i(\pi) =
\pi'$. If $\pi'$ fixes a vertex $u$, then $\pi$ fixes $u$ as well, so it is the identity, and $\pi'
= \Phi_i(\id) = \id$. And if $\pi'$ only fixes an edge $e = uv$, then $\pi'$ exchanges $u$ and $v$.
Since $\pi$ does not fix $e$, then there is an atom $A$ replaced by $e$ in $G_{i+1}$. Then $\pi|_A$
is an involutary semiregular automorphism exchanging $u$ and $v$, so $A$ is halvable. But then $e$
is a halvable edge, and thus $\pi'$ can fix it.\qed
\end{proof}

The above statement is an example of a phenomenon known in permutation group theory. Interiors
of atoms behave as \emph{blocks of imprimitivity} in the action of $\Aut(G_i)$. It is well-known
that the kernel of the action on the imprimitivity blocks is a normal subgroup of $\Aut(G_i)$.
For the example of Figure~\ref{fig:reducing_atoms}, we get $\Aut(G_1) = \gC_3$ and $\Ker(\Phi_1) =
\gC_2^3 \times \gS_3^3$. As a simple corollary, we get:

\begin{corollary}
We get
$$\Aut(G_r) = \biggl(\cdots\Bigl(\bigl(\Aut(G_0) / \Ker(\Phi_1)\bigr) / \Ker(\Phi_2)\Bigr) \cdots /
		\Ker(\Phi_r)\biggr).$$
\end{corollary}

\begin{proof}
We already proved that $\Aut(G_{i+1}) = \Aut(G_i) / \Ker(\Phi_i)$.\qed
\end{proof}

Actually, one can prove much more, that $\Aut(G_i) = \Aut(G_{i+1}) \ltimes \Ker(\Phi_i)$. First, we
describe the structure of $\Ker(\Phi_i)$.

\begin{lemma} \label{lem:kernel_char}
The group $\Ker(\Phi_i)$ is the direct product $\prod_{A \in \calA} \Fix(A)$ where $\Fix(A)$ is the
point-wise stabilizer of $G_i \setminus \int A$ in $\Aut(G_i)$. 
\end{lemma}

\begin{proof}
According to Lemma~\ref{lem:nonintersecting_interiors}, the interiors of the atoms are pairwise
disjoint, so $\Ker(\Phi_i)$ acts independently on each interior. Thus we get $\Ker(\Phi_i)$ as
the direct product of actions on each interior $\int A$ which is precisely $\Fix(A)$.\qed
\end{proof}

Alternatively, $\Fix(A)$ is isomorphic to the point-wise stabilizer of $\bo A$ in $\Aut(A)$.
Let $A_1,\dots,A_s$ be pairwise non-isomorphic atoms in $G_i$, each appearing with the multiplicity
$m_i$. According to Lemma~\ref{lem:kernel_char}, we get $\Ker(\Phi_i) \cong \Fix(A_1)^{m_1} \times
\cdots \Fix(A_s)^{m_s}$.

\begin{proposition} \label{prop:semidirect_product}
We get 
$$\Aut(G_i) \cong \Aut(G_{i+1}) \ltimes \Ker(\Phi_i).$$
\end{proposition}

\begin{proof}
According to Proposition~\ref{prop:reduction_homomorphism}c, we know that $\Ker(\Phi_i) \lhd
\Aut(G_i)$ has a complement isomorphic to $\Aut(G_{i+1})$. Actually, this already proves that
$\Aut(G_i)$ has the structure of the semidirect product. We give more details into its structure.

Each element of $\Aut(G_i)$ can be written as a pair $(\pi,\sigma)$ where $\pi \in \Aut(G_i)$ and
$\sigma \in \Ker(\Phi_i)$. We first apply $\pi$ and permute $G_i$, mapping interiors of the atoms
as blocks. Then $\sigma$ permutes the interiors of the atoms, preserving the remainder of $G_i$.

It remains to understand how composition of two automorphisms $(\pi,\sigma)$ and
$(\hat\pi,\hat\sigma)$ works. We get this as a composition of four automorphisms $\hat\sigma \circ
\hat\pi \circ \sigma \circ \pi$, which we want to write as a pair $(\tau,\rho)$.
Therefore, we need to swap $\hat\pi$ with $\sigma$. This clearly preserves $\hat\pi$,
since the action $\hat\sigma$ on the interiors does not influence it; so we get $\tau = \hat\pi
\circ \pi$.

But $\sigma$ is changed by this swapping. According to Lemma~\ref{lem:kernel_char}, we get $\sigma =
(\sigma_1,\dots,\sigma_s)$ where each $\sigma_i \in \Fix(A_i)^{m_i}$. Since $\pi$ preserves the
isorphism classes of atoms, it acts on each $\sigma_i$ independently and permutes the isomorphic
copies of $A_i$. Suppose that $A$ and $A'$ are two isomorphic copies of $A_i$ and $\pi(A) =
A'$. Then the action of $\sigma_i$ on the interior of $A$ corresponds after the swapping to the same
action on the interior of $A' = \pi(A)$. This can be described using the semidirect product, since
each $\pi$ defines an automorphism of $\Ker(\Phi_i)$ which permutes the coordinates of each
$\Fix(A_i)^{m_i}$.\qed
\end{proof}

We note that in a similar manner, Babai~\cite{babai1975automorphism,babai1996automorphism}
characterized automorphism groups of planar graphs.

\begin{lemma} \label{lem:preserved_center}
Let $G$ admit a non-trivial semiregular automorphism $\pi$. Then each $G_i$ has a central block
which is obtained from the central block of $G_{i-1}$ by replacing its atoms by colored edges.
\end{lemma}

\begin{proof}
By Proposition~\ref{prop:reduction_homomorphism}d existence of a semiregular automorphism is
preserved during the reduction.  Thus by Lemma~\ref{lem:central_block}, each $G_i$ has a central
block. Since we replace only proper atoms and dipoles in the central block, it remains as a block
after reduction. We argue by induction that it remains central as well.

Let $B$ be the central block of $G_i$ and let $B'$ be this block in $G_{i+1}$. Consider the subtree
$R'_u$ of the block tree $T'$ of $G_{i+1}$ attached to $B$ in $u$ containing the longest path in
$T'$ from $B$.  This subtree corresponds to $R_u$ in $G_i$. (See Section~\ref{sec:block_trees} for
the definition of $R_u$.) Let $\pi$ be a non-trivial semiregular automorphism in $G_i$. Then $\pi(u)
= v$, and by Lemma~\ref{lem:semiregular_action_on_blocks} we have $R_v \cong R_u$. Then $R'_v$
corresponds in $G_{i+1}$ to $R_v$ after reduction and $R'_u \cong R'_v$. Therefore $B'$ is the
central block of $G_{i+1}$.\qed
\end{proof}

\heading{Primitive Graphs.}
Recall that a graph is called primitive if it contains no atoms. If $G$ has a non-trivial
semiregular automorphism, then according to Lemma~\ref{lem:preserved_center} the central block is
preserved in the primitive graph $G_r$. We shall assume in the following that every primitive graph
has a central block.

\begin{lemma} \label{lem:primitive_graphs}
Let $G$ be a graph with a central block. Then the graph $G$ is primitive if and only if it is
isomorphic to a 3-connected graph, to a cycle $C_n$ for $n \ge 2$, or to $K_2$, or can be obtained
from these graphs by taking $U \subseteq V(G)$ such that $|U| \ge 2$ and attaching a single pendant
edge to each vertex of $U$.
\end{lemma}

\begin{proof}
The primitive graphs are depicted in Figure~\ref{fig:primitive_graphs} and clearly such graphs are
primitive. For the other implication, the graph $G$ contains a central block. All blocks attached to
it have to be single pendant edges, otherwise $G$ would contain a block atom. By removal of all
pendant edges, we get the 2-connected graph $B$ consisting of only the central block. We argue that
$B$ is isomorphic to one of the graphs above.

\begin{figure}[t!]
\centering
\includegraphics{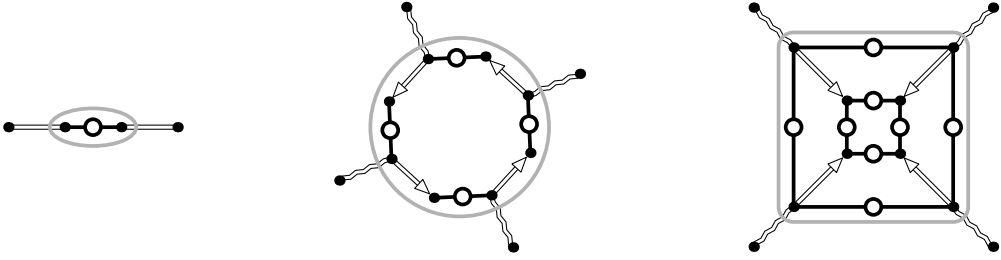}
\caption{A primitive graph with a central block is either $K_2$, $C_n$, or a 3-connected graph, in
all three cases with possible single pendant edges attached to it.}
\label{fig:primitive_graphs}
\end{figure}

Now, let $u$ be a vertex of the minimum degree of $B$. If $\deg(u)=1$, the graph $B$ has to be
$K_2$, otherwise it would not be 2-connected. If $\deg(u)=2$, then either the graph $B$ is a cycle
$C_n$, or $u$ is an inner vertex of a path connecting two vertices $x$ and $y$ of degree at least
three such that all inner vertices are of degree two. But then this path is an atom. And if $\deg(u)
\ge 3$, then every 2-cut is non-trivial, and since $B$ contains no atoms, it has to be
3-connected.\qed
\end{proof}

We note that if existence of a central block is not required, and we define atoms with respect to
the central articulation then in addition the primitive graph can be $K_1$.

\subsection{Quotients and Their Expansion} \label{sec:quotients}

Let $G_0,\dots,G_r$ be the reduction series of $G$ and let $\Gamma_0$ be a semiregular subgroup
of $\Aut(G_0)$. By repeated application of Proposition~\ref{prop:reduction_homomorphism}d, we get
the uniquely determined semiregular subgroups $\Gamma_1,\dots,\Gamma_r$ of $\Aut(G_1),\dots,\Aut(G_r)$,
each isomorphic to $\Gamma_0$. Let $H_i = G_i / \Gamma_i$ be the quotients where we preserve colors of
edges in the quotients, and let $p_i$ be the corresponding covering projection from $G_i$ to $H_i$.
Recall that $H_i$ can contain edges, loops and half-edges; depending on the action of $\Gamma_i$ on
the half-edges corresponding to the edges of $G_i$.

\begin{figure}[t!]
\centering
\includegraphics{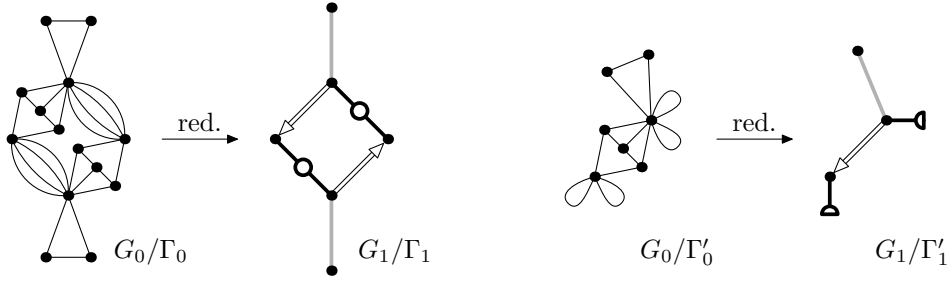}
\caption{Example of two quotients of the graph $G_0$ from Figure~\ref{fig:example_of_reduction} with
the corresponding quotients of the reduced graph $G_1$. Here $\Gamma_1 = \Phi_1(\Gamma_0)$ and
$\Gamma'_1 = \Phi_1(\Gamma'_0)$.}
\label{fig:computed_quotients}
\end{figure}

\heading{Quotients Reductions.}
Consider $H_i = G_i / \Gamma_i$ and $H_{i+1} = G_{i+1} / \Gamma_{i+1}$. We investigate relations
between these quotients. Let $A$ be an atom of $G_i$ represented by a colored edge $e$ in $G_{i+1}$.
According to Lemma~\ref{lem:atom_covering_cases}, we have three possible cases (C1), (C2) and (C3)
for the projection $p_i(A)$. It is easy to see that $\Phi_i$ is defined exactly in the way that
$p_{i+1}(e)$ corresponds to an edge in the case (C1), to a loop in the case (C2) and to a half-edge
in the case (C3).  See Figure~\ref{fig:computed_quotients} for examples. In other words, we get
the following commuting diagram: 
\begin{equation} \label{eq:red_diagram}
\begin{gathered}
\xymatrix{
G_i \ar[d]_{\Gamma_i} \ar[r]^{\rm red.} &G_{i+1}\ar[d]^{\Gamma_{i+1}}\\
H_i \ar[r]^{\rm red.} &H_{i+1}}
\end{gathered}
\end{equation}
So we can construct the graph $H_{i+1}$ from $H_i$ by replacing the projections of atoms in $H_i$ by
the corresponding projections of the edges replacing the atoms. We get the following.

\begin{lemma} \label{lem:group_reduction}
Every semiregular subgroup $\Gamma_i$ of $\Aut(G_i)$ corresponds to a unique semiregular subgroup
$\Gamma_{i+1}$ of $\Aut(G_{i+1})$.\qed
\end{lemma}

\heading{Overview of Quotients Expansions.}
Our goal is to reverse the horizontal edges in Diagram~(\ref{eq:red_diagram}), i.e, to understand:
\begin{equation} \label{eq:exp_diagram}
\begin{gathered}
\xymatrix{
G_i \ar[d]_{\Gamma_i}&G_{i+1} \ar[l]_{\rm exp.} \ar[d]^{\Gamma_{i+1}}\\
H_i&H_{i+1} \ar[l]_{\rm exp.}}
\end{gathered}
\end{equation}
Now we investigate the opposite relations. There are two fundamental questions we address in this
section in full details:
\begin{packed_itemize}
\item \emph{Question 1.} Given a group $\Gamma_{i+1}$, how many different semiregular groups $\Gamma_i$ do
we have such that $\Phi_i(\Gamma_i) = \Gamma_{i+1}$? Notice that all these groups $\Gamma_i$ are
isomorphic to $\Gamma_{i+1}$ as abstract groups, but they correspond to different actions on $G_i$.
\item \emph{Question 2.} Let $\Gamma_i$ and $\Gamma'_i$ be two groups extending $\Gamma_{i+1}$.
Under which conditions are the quotients $H_i = G_i / \Gamma_i$ and $H'_i = G_i / \Gamma'_i$
different graphs?
\end{packed_itemize}

\heading{Extensions of Group Actions.} We first deal with Question 1. Let $\Gamma_i$ and $\Gamma_{i+1}$ be the
semiregular groups such that $\Phi_i(\Gamma_i) = \Gamma_{i+1}$. Then we call $\Gamma_{i+1}$ a
\emph{reduction} of $\Gamma_i$, and $\Gamma_i$ an \emph{extension} of $\Gamma_{i+1}$.

\begin{lemma} \label{lem:group_expansion}
For every semiregular group $\Gamma_{i+1}$, there exists an extension $\Gamma_i$.
\end{lemma}

\begin{proof}
First notice that $\Gamma_{i+1}$ determines the action of $\Gamma_i$ everywhere on $G_i$ except for the
interiors of the atoms, so we just need to define it there. Let $e=uv$ be one edge of $G_{i+1}$ replacing
an atom $A$ in $G_i$. First, we assume that $A$ is not a block atom. Let $|\Gamma_{i+1}| = k$. We
distinguish two cases. Either the orbit $[e]$ contains exactly $k$ edges, or it contains $k \over 2$
edges. See Figure~\ref{fig:group_extension} for an overview.

\begin{figure}[t!]
\centering
\includegraphics{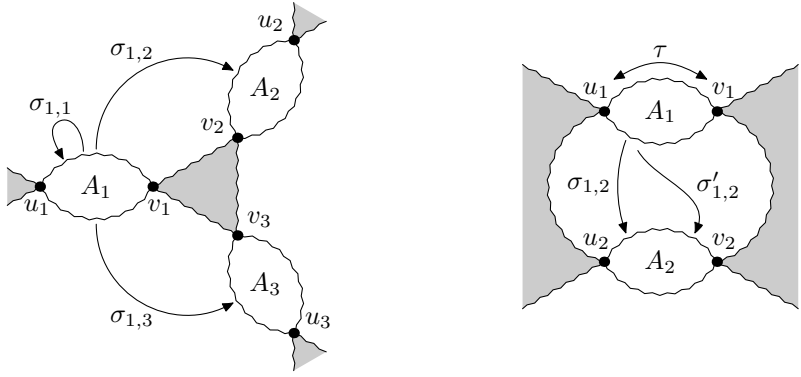}
\caption{Case 1 is depicted on the left for three edges corresponding to isomorphic atoms $A_1$,
$A_2$ and $A_3$. The depicted isomorphism are used to extend $\Gamma_{i+1}$ on interiors of these
atoms. Case 2 is on the right, with additional semiregular involution $\tau$ from the definition of
halvable atoms which transposes $u_1$ and $v_1$.}
\label{fig:group_extension}
\end{figure}

\emph{Case 1: The orbit $[e]$ contains exactly $k$ edges.}  Let $e_1,\dots,e_k$ be the orbit $[e]$
and let $u_i = \pi'(u)$ and $v_i = \pi'(v)$ for the unique $\pi'$ mapping $e$ to $e_i$. (We know
that $\pi'$ is unique because $\Gamma_{i+1}$ is semiregular.) Let $A_1,\dots,A_k$ be the atoms
corresponding to $e_1,\dots,e_k$.  The edges $e_1,\dots,e_k$ have the same color and type, and thus
the atoms $A_i$ are pairwise isomorphic and of the same type.

We define the action of $\Gamma_i$ on the interiors of $A_1,\dots,A_k$ as follows. Let
$\sigma_{1,i}$ denote any isomorphism from $A_1$ to $A_i$ such that $\sigma_{1,i}(u_1) = u_i$ and
$\sigma_{1,i}(v_1) = v_i$, with $\sigma_{1,1}$ being the identity on $A_1$. Such isomorphism exists
trivially for symmetric and halvable atoms, and they also exists for asymmetric atoms since the
action of $\Gamma_{i+1}$ preserves the orientation of $e_1,\dots,e_k$.  Then we define
$\sigma_{i,j} = \sigma_{1,j} \sigma^{-1}_{1,i}$. Let $\pi' \in \Gamma_{i+1}$ and we define the
extension $\pi$ as follows. If $\pi'$ maps $e_i$ to $e_j$, we set $\pi |_{\int A_i} =
\sigma_{i,j}$.

\emph{Case 2: The orbit $[e]$ contains exactly $\ell = {k \over 2}$ edges.} Then we have $k$
half-edges in one orbit, so in $H_i$ we get one half-edge. Let $e_1,\dots,e_\ell$ be the edges of
$[e]$. They have to be halvable, and consequently the corresponding atoms $A_1,\dots,A_\ell$ are
halvable. Let $u_i$ be an arbitrary endpoint of $e_i$ and let $v_i$ be the second endpoint of $e_i$.
Let $\tau$ be any involutory semiregular automorphism of $A_1$ which maps $u_1$ to $v_1$; we know
that such $\tau$ exists since $A_1$ is a halvable atom.

Similarly as above, we set $\sigma_{1,i}$ to be any isomorphism mapping $A_1$ to $A_i$ such that
$\sigma_{1,i}(u_1) = u_i$, and we put $\sigma_{1,1} = \id$. Moreover, we put $\sigma'_{1,1} = \tau$
and $\sigma'_{1,i} = \sigma_{1,i} \tau$. Then we put $\sigma_{i,j} = \sigma_{1,j}
\sigma^{-1}_{1,i}$, and $\sigma'_{i,j} = \sigma'_{1,j} \sigma^{-1}_{1,i}$. Let $\pi' \in
\Gamma_{i+1}$ and $\pi'(e_i) = e_j$. To define the extension $\pi$, we set $\pi|_{\int A_i}$ equal
$\sigma_{i,j}$ if $\pi'(u_i) = u_j$, and $\sigma'_{i,j}$ if $\pi'(u_i) = v_j$. 

We deal with block atoms in a similar manner as in Case 1, except the orbit $[u]$ consists of
articulations, and the orbit $[v]$ consists of leaves.  It is easy to observe that by semiregularity
of $\Gamma_{i+1}$ the constructed group $\Gamma_i$ acts semiregularly on $G_i$, as well.\qed
\end{proof}

\begin{figure}[b!]
\centering
\includegraphics{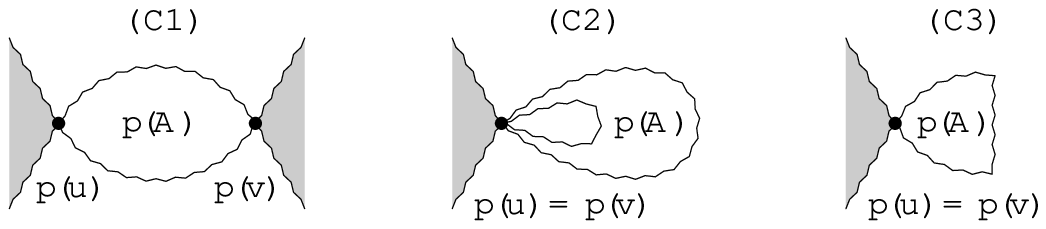}
\caption{How can $p(A)$ look in $G_i / \Gamma_i$, depending on the cases (C1), (C2) and (C3).}
\label{fig:projections_of_atoms}
\end{figure}

\begin{corollary} \label{cor:all_group_expansions}
The construction in the above proof gives all possible extensions of $\Gamma_{i+1}$.
\end{corollary}

\begin{proof}
We get all possible choices for $\Gamma_i$ in Case 1 by different choices of $\sigma_{1,i}$, and in
Case 2 by different choices of $\sigma_{1,i}$ and $\tau$.\qed
\end{proof}

\heading{Quotients of Atoms.} To answer Question 2, we first need to understand possible quotients
of an atom $A$. In Section~\ref{sec:automorphisms_of_atoms}, we stated that that for
each regular covering projection $p : G_i \to H_i$, the projection $p(A)$ satisfies one of the
three cases (C1), (C2) and (C3). Figure~\ref{fig:projections_of_atoms} shows how $p(A)$ can look in
$H_i$ depending on which of the three cases happens. If $A$ is a block atom, it is always projected
as in the case (C1).

So we get three types of quotients $p(A)$ of $A$. For (C1), we call this quotient an
\emph{edge-quotient}, for (C2) a \emph{loop-quotient} and for (C3) a \emph{half-quotient}. The
reason lying behind these names is that $p(A)$ is in $H_{i+1}$ represented by an edge, a loop or a
half-edge respectively. The following lemma allows to say ``the'' edge- and ``the'' loop-quotient of
an atom.

\begin{lemma} \label{lem:unique_quotients}
For every atom $A$, there is the unique edge-quotient and the unique loop-quotient up to isomorphism.
\end{lemma}

\begin{proof}
For the cases (C1) and (C2), we have $\int A \cong \int p(A)$, so the quotients are unique.\qed
\end{proof}

For half-quotients uniqueness does not hold. First, an atom $A$ has to be halvable to admit a
half-quotient. Then each half-quotient is determined by an involutory automorphism, and we denote
$\tau$ its restriction to $A$; recall (C3). There is a one-to-many relation between
non-isomorphic half-quotients and automorphisms $\tau$, i.e., several different automorphisms $\tau$
may give the same half-quotient.

For a proper atom, we can bound the number of non-isomorphic half-quotients by the number of
different semiregular involutions of 3-connected graphs.

\begin{lemma} \label{lem:proper_atom_quotients}
Let $A$ be a proper atom of the class $\calC$ satisfying (P2). Then there are polynomially many
non-isomorphic half-quotients of $A$.
\end{lemma}

\begin{proof}
The graph $A^+$ is essentially 3-connected graph and belongs to $\calC$. According to (C2), the
number of different semiregular subgroups of order two is polynomial in the size of $A^+$. Each
half-quotient is defined by one of these semiregular involutions which fix the edge $uv$ and
transpose $u$ and $v$.\qed
\end{proof}

For dipoles, we get the following result valid for general graphs:

\begin{lemma} \label{lem:dipole_quotients}
Let $A$ be a dipole. Then the number of pairwise non-isomorphic half-quotients is bounded by
$2^{\lfloor {e(A) \over 2}\rfloor}$ and this bound is achieved. 
\end{lemma}

\begin{proof}
Figure~\ref{fig:exponentially_many_quotients} shows a construction of dipoles achieving the upper
bound. It remains to argue correctness of the upper bound.

\begin{figure}[b!]
\centering
\includegraphics{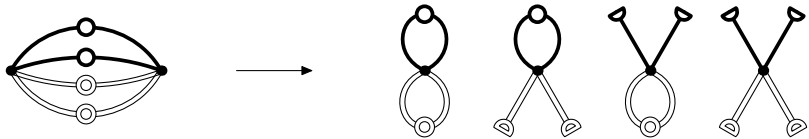}
\caption{An example of a dipole with a pair of white halvable edges and a pair of black halvable
edges (corresponding to two isomorphism classes of halvable atoms). There exist four pairwise
non-isomorphic half-quotiens. This example can easily be generalized to exponentially many pairwise
non-isomorphic quotients by introducing more pairs of halvable edges of additional colors.}
\label{fig:exponentially_many_quotients}
\end{figure}

First, we derive the structure of all involutory semiregular automorphisms $\tau$ acting on $\int
A$. We have no freedom concerning the non-halvable edges of $A$: The undirected edges of each color
class has to paired by $\tau$ together. Further, each directed edge has to be paired with a directed
edges of the opposite direction and the same color. It remains to describe possible action of $\tau$
on the remaining at most $e(A)$ halvable edges of $A$. These edges belong to $c$ color classes
having $m_1,\dots,m_c$ edges.  Each automorphism $\tau$ has to preserve the color classes, so it
acts independently on each class.

We concentrate only for one color class having $m_i$ edges. We bound the number $f(m_i)$ of pairwise
non-isomorphic quotients of this class. Then we get the upper bound
\begin{equation} \label{eq:dipole_bound}
\prod_{1 \le i \le c} f(m_i)
\end{equation}
for the number of non-isomorphic half-quotients of $A$.

An edge $e$ fixed in $\tau$ is mapped into a half-edge of the given color in the half-quotient of
$A$. And if $\tau$ maps $e$ to $e' \ne e$, then we get a loop in the half-quotient. The resulting
half-quotient only depends on the number of fixed edges and fixed two-cycles in the considered color
class.  We can construct at most $f(m_i) = \lfloor{m_i  \over 2}\rfloor + 1$ pairwise non-isomorphic
quotients, since we may have zero to $\lfloor{m_i \over 2}\rfloor$ loops with the complementing
number of half-edges.

The bound~(\ref{eq:dipole_bound}) is maximized when each class contains exactly two edges. (Except
for one class containing three or one edge if $m$ is odd.)\qed
\end{proof}

This bound plays the key role for the complexity of our meta-algorithm of Section~\ref{sec:algorithm}; in one
subroutine, we iterate over all half-quotients of a dipole. Also the structure of all possible
quotients is important.

\heading{Quotient Expansion.}
When we know all quotients of atoms, we can construct from given $H_{i+1}$ all quotients $H_i$ as
follows. We say that two quotients $H_i$ and $H'_i$ extending $H_{i+1}$ are different if there exists no
isomorphism of $H_i$ and $H'_i$ which fixes the vertices and edges common with $H_{i+1}$. (But $H_i$
and $H'_i$ still might be isomorphic.)

\begin{proposition} \label{prop:quotient_expansion}
Every quotient $H_i$ of $G_i$ can be constructed from some quotient $H_{i+1}$ of $G_{i+1}$ by
replacing each edge, loop and half-edge corresponding to an atom of $G_i$ by an edge-, loop-, or
half-quotient respectively. Moreover, for different choices of $H_{i+1}$ and of half-quotients we
get different graphs $H_i$.
\end{proposition}

\begin{proof}
Let $H_{i+1} = G_{i+1} / \Gamma_{i+1}$ and let $H_i$ be constructed in the above way. We first argue
that $H_i$ is a quotient of $G_i$, i.e., it is equal to $G_i / \Gamma_i$ for some $\Gamma_i$ extending
$\Gamma_{i+1}$. To see this, it is enough to construct $\Gamma_i$ in the way described in the proof
of Lemma~\ref{lem:group_expansion}. We choose $\sigma_{1,i}$ arbitrarily, and the involutory
permutations $\tau$ are prescribed by chosen half-quotients replacing half-edges. It is easy to see
that the resulting graph is the constructed $H_i$. We note that only the choices of $\tau$ matter,
for arbitrary choices of $\sigma_{1,i}$ we get the same quotients.

On the other hand, if $H_i$ is a quotient, it replaces the edges, loops and half-edges of $H_{i+1}$
by some quotients, so we can generate $H_i$ in this way. The reason is that according to
Corollary~\ref{cor:all_group_expansions}, we can generate all $\Gamma_i$ extending $\Gamma_{i+1}$ by
some choices $\sigma_{1,i}$ and $\tau$.

For the last statement, according to Lemma~\ref{lem:unique_quotients}, the edge and
loop-quotients are uniquely determined, so we are only free in choosing different half-quotients.
For different choices of half-quotients, we get different graphs $H_i$.\qed
\end{proof}

For instance suppose that $H_{i+1}$ contains a half-edge corresponding to the dipole from
Figure~\ref{fig:exponentially_many_quotients}. Then in $H_i$ we can replace this half-edge by one of
the four possible half-quotients of this dipole.

\begin{corollary} \label{cor:unique_expansion}
If $H_{i+1}$ contains no half-edge, then $H_i$ is uniquely determined. So for odd order of
$\Gamma_r$, the quotient $H_r$ uniquely determines $H_0$.
\end{corollary}

\begin{proof}
Implied by Proposition~\ref{prop:quotient_expansion} and Lemma~\ref{lem:unique_quotients} which
states that edge- and loop-quotients are uniquely determined. If the order of $\Gamma_r$ is odd, no
half-edges are constructed.\qed
\end{proof}

\heading{The Block Structure of Quotients.}
The following properties are key for identifying quotients of atoms in the input graph $H$. The
approach used in the meta-algorithm is to find a way how to expand $H_r$ by repeated
application of Proposition~\ref{prop:quotient_expansion} to $H_0$ which is isomorphic to the input
$H$.

A block atom $A$ of $G_i$ is always projected by (C1), and so it corresponds to a block atom of
$H_i$. Suppose that $A$ is a proper atom or a dipole, and let $\bo A = \{u,v\}$. For (C1) we get
$p(u) \ne p(v)$, and for (C2) and (C3) we get $p(u) = p(v)$. For (C1), $p(A)$ is isomorphic to an
atom in $H_i$. For (C2) and (C3) is $p(u)$ an articulation of $H_i$, and $p(A)$ corresponds to a
pendant star, or a pendant block with possible attached single pendant edges.

\begin{lemma} \label{lem:block_structure}
The block structure of $H_{i+1}$ is preserved in $H_i$ with possible some new pendant blocks
attached.
\end{lemma}

\begin{proof}
Edges inside blocks are replaced using (C1) by edge-quotients of block atoms, proper atoms and
dipoles which preserves 2-connectivity. The new pendant blocks in $H_i$ are created by replacing
pendant edges with the block atoms, loops by loop-quotients, and half-edges by half-quotients.\qed
\end{proof}
 
\section{Meta-algorithm} \label{sec:algorithm}

In this section, we establish the fixed parameter tractable algorithm of
Theorem~\ref{thm:metaalgorithm}. We show that for a class $\calC$ satisfying (P0) to (P3) we
can solve $\rcover(G,H)$ in time $\O^*(2^{e(H)/2})$. We use the property (P3) for essentially
3-connected graphs with colored pendant edges which code colors and lists of colors.

Let $k = |G|/|H|$, and we assume that $k \ge 2$. (If $k$ is not an integer, then clearly $G$ does
not cover $H$. If $k=1$, then we can test it using the algorithm for graph isomorphism given by (P1)
whether $G \cong H$.) The algorithm proceeds in the following major steps:
\begin{packed_enum}
\item We construct the reduction series for $G = G_0,\dots,G_r$ terminating with the unique
primitive graph $G_r$.  Throughout the reduction the central block is preserved, otherwise according
to Lemma~\ref{lem:preserved_center} there exists no semiregular automorphism of $G$ and we output
``no''. According (P0), the reduction preserves the class $\calC$, and also every atom belongs to
$\calC$.
\item Using (P2), we compute $\Aut(G_r)$ and construct a list of all subgroups $\Gamma_r$ of the
order $k$ acting semiregularly on $G_r$. The number of subgroups in the list is polynomial by (P2).
\item For each $\Gamma_r$ in the list, we compute $H_r = G_r / \Gamma_r$. We say that a graph $H_r$
is \emph{expandable} if there exists a sequence of extensions repeatedly applying
Proposition~\ref{prop:quotient_expansion} which constructs $H_0$ isomorphic to $H$. We test
the expandability of $H_r$ using dynamic programming while using (P1) and (P3).
\end{packed_enum}
It remains to explain details of the third step, and prove the correctness of the algorithm.

\subsection{Testing Expandability Using Dynamic Programming}

\heading{Catalog of Atoms.}
During the reduction phase of the algorithm, we construct the following \emph{catalog of atoms}
forming a database of all discovered atoms and their quotients. We are not very concerned with
a specific implementation of the algorithm, so this catalog is mainly used to simplify description.
For each isomorphism class of atoms represented by an atom $A$, we store the following information
in the catalog:
\begin{packed_itemize}
\item the atom $A$,
\item the corresponding colored edge of a given type representing the atom in the reduction,
\item the unique edge- and loop-quotients of $A$.
\end{packed_itemize}
For dipoles, according to Lemma~\ref{lem:dipole_quotients} we can have exponentially many
non-isomorphic half-quotients, and so we work with their half-quotients implicitly in the dynamic
programming.

If $A$ is not a dipole, we compute a list of all its pairwise non-isomorphic half-quotients, and
store them in the catalog in the following way. A half-quotient $Q$ might not be 3-connected, and so
we apply a reduction series on $Q$, and add all atoms discovered by the reduction to the catalog.
(We do not compute their half-quotients. They are never realized unless these atoms are directly
found in $G$ as well.) When the reduction series finishes, this half-quotient is reduced to a
primitive graph. We note that $\bo Q$, being a single vertex of the half-quotient, behaves like the
central block in the definition of atoms, i.e., it is never reduced.

Further, if a halvable dipole consists of exactly two edges of the same color, we compute its
half-quotient consisting of just the single loop attached, and we add this quotient to the catalog.
The reason is that that this quotient behaves exactly as a loop-quotient of some proper atom.

\begin{lemma} \label{lem:catalog_size}
The catalog contains polynomially many quotients and atoms.
\end{lemma}

\begin{proof}
First we deal with the number of atoms in $G_0,\dots,G_r$. Notice that by replacing an interior of
an atom, the total number of vertices and edges is decreased; the interiors of atoms in each $G_i$
contain at least two vertices and edges in total and are pairwise disjoint (implied by
Lemma~\ref{lem:nonintersecting_interiors}). Thus we have $\O(n+m)$ atoms in $G_0,\dots,G_r$.

For the number of quotients, let $A$ be a block or a proper atom. Each half-quotient of $A$ is created
by some semiregular action of an involution on $A$. According to (P2), there are polynomially many
half-quotients. For each half-quotient, we can have at most linearly many atoms in its reduction
series.  And by Lemma~\ref{lem:unique_quotients} we have the unique edge- and loop-quotient. So the
total number of atoms and their quotients is polynomial.\qed
\end{proof}

Throughout the algorithm, we repeatedly ask whether some atom or some of its quotients is contained
in the catalog. Each such query can be answered in polynomial time.

\begin{figure}[b!]
\centering
\includegraphics{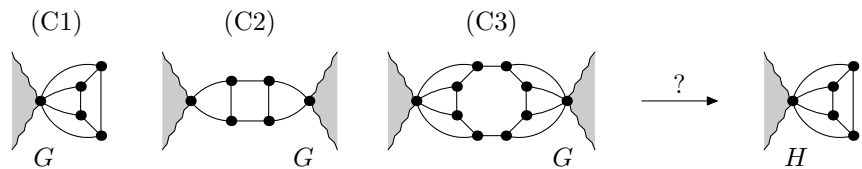}
\caption{For a pendant block of $H$, there are three possible preimages in $G$. It could be a block
atom mapped by (C1), or a proper atom mapped by (C2), or another proper atom mapped by (C3) (where
the half-quotient is created by $180^\circ$ rotation $\tau$).}
\label{fig:preimages_of_block}
\end{figure}

\heading{Preimages of a Pendant Block.} We now illustrate the fundamental difficulty in testing
whether $H_r$ is expandable to $H$, for simplicity we do it on pendant blocks. Suppose that $H$ has a pendant
block as in Figure~\ref{fig:preimages_of_block}.  Then there is no way to know whether this block
corresponds in $G$ to an edge-quotient of a block atom, or to a loop-quotient of a proper atom, or
to a half-quotient of another proper atom. It can easily happen that the catalog offers all three
options.  So without exploiting some additional information from $H$, there is no way to know what
is the preimage of this pendant block.

In our approach, we do not decide everything in one stage, instead we just remember a list of
possibilities. The dynamic programming deals with these lists and computes further lists for larger
parts of $H$.

\heading{Atoms in Quotients.}
We define atoms in the quotient graphs similarly as in Section~\ref{sec:atoms} with only one
difference. We choose one arbitrary block/articulation called the \emph{core} in $H_r$; for instance, we
can choose the central block/articulation. The core plays the role of the central block in the
definition of parts and atoms. Also, in the definition we consider half-edges and
loops as pendant edges, so they do not form block-atoms.

\begin{figure}[b!]
\centering
\includegraphics{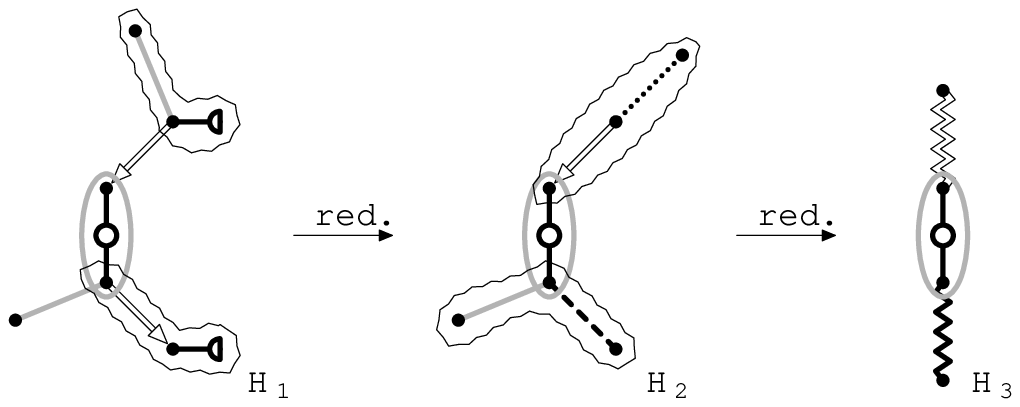}
\caption{The graph $H_1$ is one quotient of $G_1$ from Figure~\ref{fig:example_of_reduction}. We
further reduce it to $H_3$ with respect to the core block depicted in gray. Notice that $H_1$ and
$H_2$ only contain block atoms.}
\label{fig:quotient_reduction}
\end{figure}

We proceed with the reductions in $H_r$ further till we obtain a primitive quotient graph $H_s$, for
some $s \ge r$; see Figure~\ref{fig:quotient_reduction}. Notice that all atoms in
$H_r,\dots,H_{s-1}$ are necessarily block atoms since otherwise $G_r$ would contain some proper
atoms or dipoles and it would not be primitive. We add the newly discovered atoms to the catalog.

Now, the graph $H_s$ consists of the core together with some pendant edges, loops and half-edges.
Let $H_0,\dots,H_{s-1}$ be the graphs obtained by an expansion series of $H_s$ using
Proposition~\ref{prop:quotient_expansion}. Notice that the core is preserved as an
articulation/block in all these graphs. The core can be only changed by replacing of its colored
edges by edge-quotients.  Then the core in $H_0$ has to correspond to some block or articulation of
$H$. We test all possible positions of the core in $H$. (We have $\O(n)$ possibilities, so we
run the dynamic programming algorithm multiple times.) In what follows, we have the core fixed in
$H$ as well.

\heading{Overview of Dynamic Programming.}
Our goal is to apply a reduction series on $H$ defining $\calH_0,\dots,\calH_t$. As already
discussed above, we do not know which parts of $G$ project to different parts of $H$. Therefore each
$\calH_i$ is a set of graphs, and $\calH_t$ is a set of primitive graphs. We then
determine expandability of $H_r$ by testing whether $H_s \in \calH_t$.

Since each set $\calH_i$ can contain a huge number of graphs, we represent it implicitly in the
following manner. Each $\calH_i$ is represented by one graph $\calR_i$ with some colored edges and
with so-called \emph{pendant elements} attached to some vertices. Here, each pendant element can
represent a pendant edge, loop or half-edge at the same time. Further for each pendant element $x$,
we have a \emph{list} $\frakL(x)$ of possible realizations of the corresponding subgraph of $H$ by
the quotients from the catalog. Each graph of $\calH_i$ is created for $\calR_i$ by replacing the
pendant elements by some edges, loops and half-edges from the respective lists.

\begin{figure}[t!]
\centering
\includegraphics{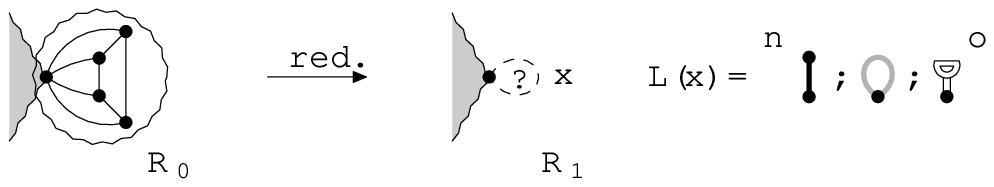}
\caption{Let $x$ be the pendant element corresponding to the pendant block
of $H$ depicted in Figure~\ref{fig:preimages_of_block}. Then $\frakL(x)$ contains three different
elements if all three atoms depicted in Figure~\ref{fig:preimages_of_block} are contained in the
catalog.}
\label{fig:list_example}
\end{figure}

Pendant elements of $\calR_i$ correspond to block parts of $H$ with pairwise disjoint interiors. (Recall
that block parts are defined in Section~\ref{sec:atoms}.) A pendant element $x$ contains an
edge/loop/half-edge in $\frakL(x)$ if and only if it is possible to expand this edge/loop/half-edge
to a graph isomorphic to the given block part. Further for each element of the list, we remember how
to do this expansion. For an example, see Figure~\ref{fig:list_example}.

Testing whether $H_s \in \calH_t$ is equivalent to testing whether there exists an \emph{embedding}
$H_s \hookto \calR_t$. Here, the embedding is an isomorphism $\pi : H_s \to \calR_t$ which maps pendant
edges, loops and half-edges of $H_s$ to the pendant elements of $\calR_t$. Further, we require that
the list of the pendant element $\pi(e)$ contains the mapped edge $e$.

These lists are used in the dynamic programming to compute $\calR_{i+1}$ from $\calR_i$.
According to Lemma~\ref{lem:catalog_size}, we have polynomially many atoms, and so the size of
each list is polynomial. We note that one list may contain many half-edges. 

\begin{lemma} \label{lem:list_properties}
Each list contains at most one edge and at most one loop. Further, if two lists contain the same edge
or loop, then they have to be equal.
\end{lemma}

\begin{proof}
If the pendant element of $\calR_i$ is fully expanded, it corresponds to one block part of $H$.
Suppose that an edge- or a loop-quotient is in the list. If it is fully expanded, then it has to be
isomorphic to this block part. But according to Lemma~\ref{lem:unique_quotients}, there is only one
way how to expand an edge- or a loop-quotient, because it can never contain half-edges.\qed
\end{proof}

\heading{Reductions with Lists.}
Suppose that we know $\calR_i$, and we want to apply one step of the reduction and compute
$\calR_{i+1}$. First, we find all atoms in $\calR_i$. (We define atoms with respect to the core
block, and we consider pendant elements as pendant edges.) To construct $\calR_{i+1}$ from
$\calR_i$:
\begin{packed_itemize}
\item We replace dipoles and proper atoms by edges of the corresponding colors from the catalog.
If the corresponding dipole or proper atom is not contained in the catalog, we halt the reduction
procedure.
\item We replace block atoms by pendant elements with constructed lists. If some list is empty, we
again halt the reduction.
\end{packed_itemize}
It remains to describe the construction of the lists for the created pendant elements.

Let $A$ be an atom in $\calR_i$, replaced by an edge/pendant element $e$, and we want to compute the
list for $e$. We call an atom $A$ as a \emph{star atom} if it consists of an articulation with
attached pendant edges, loops and pendant elements.

\begin{lemma} \label{lem:lists_for_nonstar_atoms}
Let $A$ be a non-star block atom in $\calR_i$. Then we can compute its list $\frakL(A)$ in polynomial time.
\end{lemma}

\begin{figure}[t!]
\centering
\includegraphics{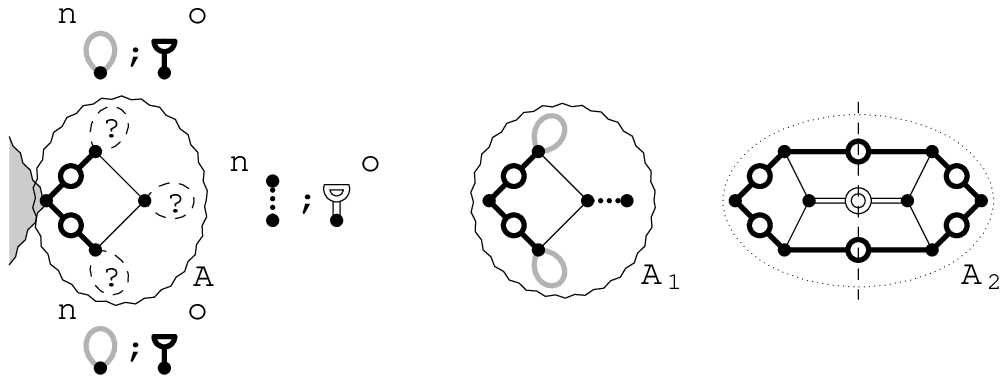}
\caption{On the left, a non-star block atom $A$ in $\calR_i$ with lists of its pendant elements
depicted. On the right, two possible quotient which can be embedded into $A$. So the list of the
pendant element replacing $A$ in $\calR_{i+1}$ contains a pendant edge corresponding to the block
atom $A_1$ and a half-edge corresponding to a half-quotient of the proper atom $A_2$. }
\label{fig:lists_for_nonstar_atoms}
\end{figure}

\begin{proof}
We iterate over all quotients in the catalog. For one such quotient $Q$, we determine whether $Q
\hookto A$ where $\bo Q$ is mapped to $\bo A$ as follows. Notice that $A$ is essentially
3-connected, and thus $Q$ has to be 3-connected as well. (Otherwise an embedding does not exist.)
So by (P3), we can test in polynomial time whether $Q \hookto A$ by coding colors of the pendant edges
of $Q$ by the colors of the vertices, and the lists of the pendant elements by lists of colors for
vertices of $A$. If $Q \hookto A$, we add the edge/loop/half-edge representing this quotient to the
list, and we remember the constructed mapping $Q \hookto A$. See
Figure~\ref{fig:lists_for_nonstar_atoms} for an example.\qed
\end{proof}

\begin{lemma} \label{lem:lists_for_star_atoms}
Let $A$ be a star atom in $\calR_i$. Then we can compute its list $\frakL(A)$ in time
$\O^*(2^{e(H)/2})$ where $e(H)$ is the number of edges in $H$.
\end{lemma}

\begin{proof}
Each star atom of $\calR_i$ corresponds either to a block atom isomorphic to a star, or to the
loop-quotient or a half-quotient of a dipole. Star atoms involve half-quotients of dipoles, and
Lemma~\ref{lem:dipole_quotients} states that a dipole can have exponentially many pairwise
non-isomorphic half-quotients. For a dipole, we iterate over all possible half-quotients which gives
$2^{e(H)/2}$ part in the complexity bound.

\emph{Case 1: Dipoles.}
First, we show how to deal with dipoles. We iterate over all dipoles in the catalog and try to
add them to the list. For each dipole, we first test whether the lists of the pendant elements
attached to the star atom $A$ are compatible with the unique loop quotient. Then we iterate over all
half-quotients of $D$. Let $D$ be one
dipole in the catalog with $\bo D = \{u,v\}$ and let $Q$ be one of its at most $2^{e(H)/2}$ possible
quotients.  Recall from the proof of Lemma~\ref{lem:dipole_quotients} that an edge of a dipole either
projects to a half-edge, or two edges of the same color and type project to one loop. So each $Q$
consists of loops and half-edges attached to $u$, and they have to be matched to the pendant
elements of $A$ with the corresponding lists.

We can reduce this problem to finding a perfect matching in bipartite graphs: Here, one part is
formed by the loops and the half-edges of $Q$, and the other part is formed by the pendant elements
of $A$. A loop/half-edge is adjacent to a pendant element, if and only if the corresponding list
contains this loop/half-edge. Each perfect matching defines one embedding $Q \hookto A$. 
We add the half-edge of the dipole $D$ to the list if there exists a perfect matching.

\emph{Case 2: Star Atoms.}
We iterate over all star atoms of the catalog, let $S$ be one of them. The star atom $S$ consists
of pendant edges, loops and half-edges attached to one vertex. Some of these half-edges correspond
to dipoles, and some to proper atoms. Let $e_1,\dots,e_d$ be the half-edges corresponding to the
dipoles $D_1,\dots,D_d$. We construct all quotients $Q$ of $S$ by replacing $e_1,\dots,e_d$ by all
possible choices of the half-quotients $Q_1,\dots,Q_d$ of $D_1,\dots,D_d$.  In total, we have at
most $2^{e(H)/2}$ different quotients $Q$ of $S$. For each $Q$, we test by the matching procedure
described above whether the edge representing the star atom $S$ should be added to the constructed
list. If yes, we add $S$ to the list.

The procedure computes the list correctly since we test all possible quotients from catalog, and for
each quotient we test all possibilities how it could be matched to $A$.  For each quotient $Q$, the
running time is clearly polynomial, and we have at most $2^{e(H)/2}$ quotients.\qed
\end{proof}

Algorithm~\ref{alg:computing_list} gives the pseudocode for computation of the list of a pendant
element replacing an atom $A$. If the returned list is empty, we halt the reduction. Either $H_r$ is
not expandable to $H$, or we have chosen a wrong core in $H$. The following diagram shows the
overview of the meta-algorithm:

\begin{equation} \label{eq:algorithm_diagram}
\begin{gathered}
\xymatrix{
G_0 \ar[r]^{{\rm red.}} \ar[d]_{\Gamma_0}&G_1 \ar[r]^{{\rm
red.}}\ar[d]_{\Gamma_1}&\cdots\ar[r]^{{\rm red.}}&G_r\ar[d]_{\Gamma_r}&&&\\%
H_0&H_1 \ar[l]_{{\rm exp.}}&\cdots \ar[l]_{{\rm exp.}}&H_r \ar@/^/[r]^{{\rm red.}} \ar[l]_{{\rm
exp.}}&H_{r+1} \ar@/^/[r]^{{\rm red.}} \ar@/^/[l]^{{\rm exp.}}&\cdots \ar@/^/[r]^{{\rm red.}}
\ar@/^/[l]^{{\rm exp.}}&H_s\vphantom{A_{s_d}} \ar@/^/[l]^{{\rm exp.}} \ar@{^{(}->}[d]^{?} \\%
\calR_0 \ar[r]^{{\rm red.}}&\calR_1 \ar[r]^{{\rm red.}}&\calR_2 \ar[r]^{{\rm red.}}&
\calR_3 \ar[r]^{{\rm red.}}&\calR_4 \ar[r]^{{\rm red.}}&\cdots \ar[r]^{{\rm red.}}&\calR_t&\\}
\end{gathered}
\end{equation}

\begin{algorithm}[t!]
\caption{The subroutine for computing lists of realizations} \label{alg:computing_list}
\begin{algorithmic}[1]
\REQUIRE An atom $A$ of $\calR_i$.
\ENSURE The list $\frakL(x)$ of the pendant element $x$ replacing $A$ in $\calR_{i+1}$.
\medskip
\STATE Initiate the empty list $\frakL(x)$.
\IF {$A$ is a non-star atom}
	\STATE Iterate over all quotients from the catalog.
	\STATE For each quotient $Q$, apply (P3) to test whether $Q \hookto A$.
	\STATE If some embedding exists, we add the edge/loop/half-edge of $Q$ to $\frakL(x)$
	together with the embedding.
\ENDIF
\medskip
\IF {$A$ is a star atom}
	\STATE Iterate over all dipoles and star atoms in the catalog.
	\FOR{each dipole $D$}
		\STATE Test whether the loop-quotient of $D$ matches the lists, if yes add the loop
		representing $D$ to $\frakL(x)$.
		\STATE Iterate over all half-quotients $Q$ of $D$.
		\FOR{each half-quotient $Q$}
			\STATE Test existence of a perfect matching between loops and half-edges of $Q$ and the
			pendant elements of $A$.
			\STATE If a perfect matching exists, add the half-edge of $D$ to $\frakL(x)$ together
			with its embedding, and proceed with the next dipole.
		\ENDFOR
	\ENDFOR
	\FOR{each star atom $S$}
		\STATE Compute all quotients $Q$ be replacing the edges $e_1,\dots,e_d$ corresponding to the
		dipoles by all possible combinations of half-quotients $Q_1,\dots,Q_d$.
		\FOR{each quotient $Q$}
			\STATE Test existence of a perfect matching exactly as above.
			\STATE If it exists, add the edge of $S$ to $\frakL(x)$ with its embedding, and proceed
			with the next star atom.
		\ENDFOR
	\ENDFOR
\ENDIF
\medskip
\RETURN The constructed list $\frakL(x)$.
\end{algorithmic}
\end{algorithm}

\begin{lemma} \label{lem:expandability_testing}
We can test whether $H_s \hookto \calR_t$ in polynomial time.
\end{lemma}

\begin{proof}
The graph $H_s$ consists of the core together with edges, loops and half-edges attached to it.  The
graph $\calR_t$ consists of the core together with pendant elements with computed lists attached to
it.  Similarly as in the proof of Lemma~\ref{lem:lists_for_nonstar_atoms}, since the core of
$\calR_t$ is essentially 3-connected, we use (P3) and test whether $H_s \hookto \calR_t$.\qed
\end{proof}

The following lemma states that we can test the expandability of $H_r$ using $\calR_t$.

\begin{lemma} \label{lem:expandability_correctness}
We have $H_s \hookto \calR_t$ for some choice of the core in $H$ if and only if $H_r$ is expandable to
$H_0$ which is isomorphic to $H$.
\end{lemma}

\begin{proof}
First suppose that $H_s \hookto \calR_t$ for some choice of the core. Then the embedding of $H_s$
gives a sequence of replacements of edges, loops and half-edges by edge-quotients, loop-quotients
and half-quotients respectively such that the resulting graph is isomorphic to $H$.  We can apply
these replacements in any order, since they modify the graph independently.

Therefore, we first replace pendant edges by the unique block atoms in $H_s$ to get $H_r$, this has
to be compatible with the sequence of replacements defined by $H_s \hookto \calR_t$. Then we do
replacements in the manner of Proposition~\ref{prop:quotient_expansion}, and construct the
expansions $H_{r-1},\dots,H_0$. Since we expand according to the embedding $H_s \hookto \calR_t$,
the constructed $H_0$ is isomorphic to $H$. 

On the other hand, suppose that $H_r$ is expandable to $H_0$ which is isomorphic to $H$. Then
according to Lemma~\ref{lem:block_structure}, the core of $H_r$ is preserved in $H$, so it has to
correspond to some block or to some articulation of $H$. We claim that $H_s \hookto \calR_t$ for this choice of the
core. There exists a sequence of replacements from $H_s$ which constructs $H_0$, and this sequence
of replacements have to be possible in $\calR_t$. Thus $H_s \hookto \calR_t$.\qed
\end{proof}

\subsection{Proof of The Main Theorem}

Now, we are ready to establish the main algorithmic result of the paper; see
Algorithm~\ref{alg:metaalgorithm} for the pseudocode.

\begin{proof}[Theorem~\ref{thm:metaalgorithm}]
We recall the main steps of the algorithm and discuss their time complexity. The reduction series
$G_0,\dots,G_r$ can be computed in polynomial time. The property (P2) ensures that there are
polynomially many semiregular subgroups $\Gamma_r$ of $\Aut(G_r)$ which can be computed in
polynomial time. (If $G_r$ is an edge or a cycle, then it is clearly true as well.)

We compute the quotient $H_r = G_r / \Gamma_r$, and we fix the central block/articulation of $H_r$ as
the core. Then we compute the reduction series $H_r,\dots,H_s$ by replacing only block atoms.
We compute $\calR_t$ for all choices of the core in time $\O^*(2^{e(H)/2})$ and we can
test $H_s \hookto \calR_t$ in polynomial time according to Lemma~\ref{lem:expandability_testing}. If
we succeed for at least one choice of the core, then $G$ regularly covers $H$. To certify this, we
construct the regular covering projection as follows: We proceed as in the proof of
Lemma~\ref{lem:expandability_correctness} and construct $\Gamma_{r-1},\dots,\Gamma_0$. Then
$\Gamma := \Gamma_0$ gives a regular covering projection $p: G \to H$. If $H_s \not\hookto \calR_t$
for all choices of the core, we proceed with the next choice of $\Gamma_r$. If we fail for all
choices of $\Gamma_r$, the algorithm outputs ``no''.

It remains to argue correctness of the algorithm. First suppose that the algorithm succeeds.
According to Lemma~\ref{lem:expandability_correctness} and
Proposition~\ref{prop:quotient_expansion}, we construct a semiregular subgroup $\Gamma$ of $\Aut(G)$
such that $G / \Gamma \cong H$ which proves that $G$ regularly covers $H$. On the other hand,
suppose that there exists a semiregular $\Gamma$ such that $H \cong G / \Gamma$. According to
Lemma~\ref{lem:group_reduction}, it corresponds to the unique semiregular $\Gamma_r$ on $G_r$ which
is one of the semiregular subgroups tested by the algorithm. Therefore $H_r$ has to be expandable to
$H_0$ isomorphic to $H$, and we detect this correctly according to
Lemma~\ref{lem:expandability_correctness}.\qed
\end{proof}

\begin{algorithm}[t!]
\caption{The meta-algorithm for regular covers -- $\rcover$} \label{alg:metaalgorithm}
\begin{algorithmic}[1]
\REQUIRE A graph $G$ of $\calC$ satisfying (P0), (P1), (P2) and (P3), and a graph $H$.
\ENSURE A regular covering projection $p: G \to H$ if it exists.
\medskip
\STATE Compute the reduction series $G_0,\dots,G_r$ ending with the unique primitive graph $G_r$.
\STATE During the reductions, we construct the catalog by introducing all detected atoms and their
quotients.
\medskip
\STATE Using (P2), we compute all semiregular subgroups $\Gamma_r$ of $\Aut(G_r)$.
\FOR{each semiregular $\Gamma_r$}
	\STATE Compute the quotient $H_r = G_r / \Gamma_r$.
	\STATE Choose, say, the central block/articulation of $H_r$ as the core.
	\STATE Compute the reduction series $H_r,\dots,H_s$ with respect to the core.
	\STATE Introduce newly discovered atoms to the catalog.
	\medskip
	\FOR{each guessed position of the core in $H$}
		\STATE Compute the reduction series $\calR_0,\dots,\calR_t$ as follows.
		\FOR{each block atom $A$ in $\calR_i$}
			\STATE Compute its list using Algorithm~\ref{alg:computing_list}.
		\ENDFOR
		\STATE To construct $\calR_{i+1}$, replace the proper atoms and dipoles of $\calR_i$ by
		colored edges, and the block atoms by pendant elements with the computed lists.
	\ENDFOR
	\medskip
	\STATE Test whether $H_s \hookto \calR_t$ by trying all possible mappings of the core of $H_s$
	to the core of $\calR_t$.
	\STATE If yes, use the embedding to compute the expansions $H_{s-1},\dots,H_0$ such that $H_0
	\cong H$. And construct the groups $\Gamma_{r-1},\dots,\Gamma_0$.
	\STATE The group $\Gamma := \Gamma_0$ defines the regular covering projection $p : G \to H$.
\ENDFOR
\medskip
\RETURN The regular covering projection $p$ if it is constructed, ``no'' otherwise.
\end{algorithmic}
\end{algorithm}

\begin{proof}[Corollary~\ref{cor:simple_cases}]
If $G$ is 3-connected, then it is primitive and only block atoms can appear. So the expansion runs
in polynomial time. If $|\Gamma| = |\Gamma_r|$ is odd, no half-edges occur in $H_r$, and so
according Corollary~\ref{cor:unique_expansion} the expansion gives the unique graph $H_0$. We can
just test whether $H_0 \cong H$, or we can compute the reduction series $\calR_0,\dots,\calR_t$
while ignoring half-quotients.\qed
\end{proof}

\subsection{More Details about Star Atoms and Their Lists} \label{sec:comb_int_star_atoms}

In this section, we give details and insides on lists of star atoms in $\calR_i$. We show
that this problem can be reduced to finding a certain generalization of a perfect matching which we call
\ivmatch. Here we describe a complete derivation of this problem, and in Conclusions we just give
its combinatorial statement.

\begin{figure}[t!]
\centering
\includegraphics{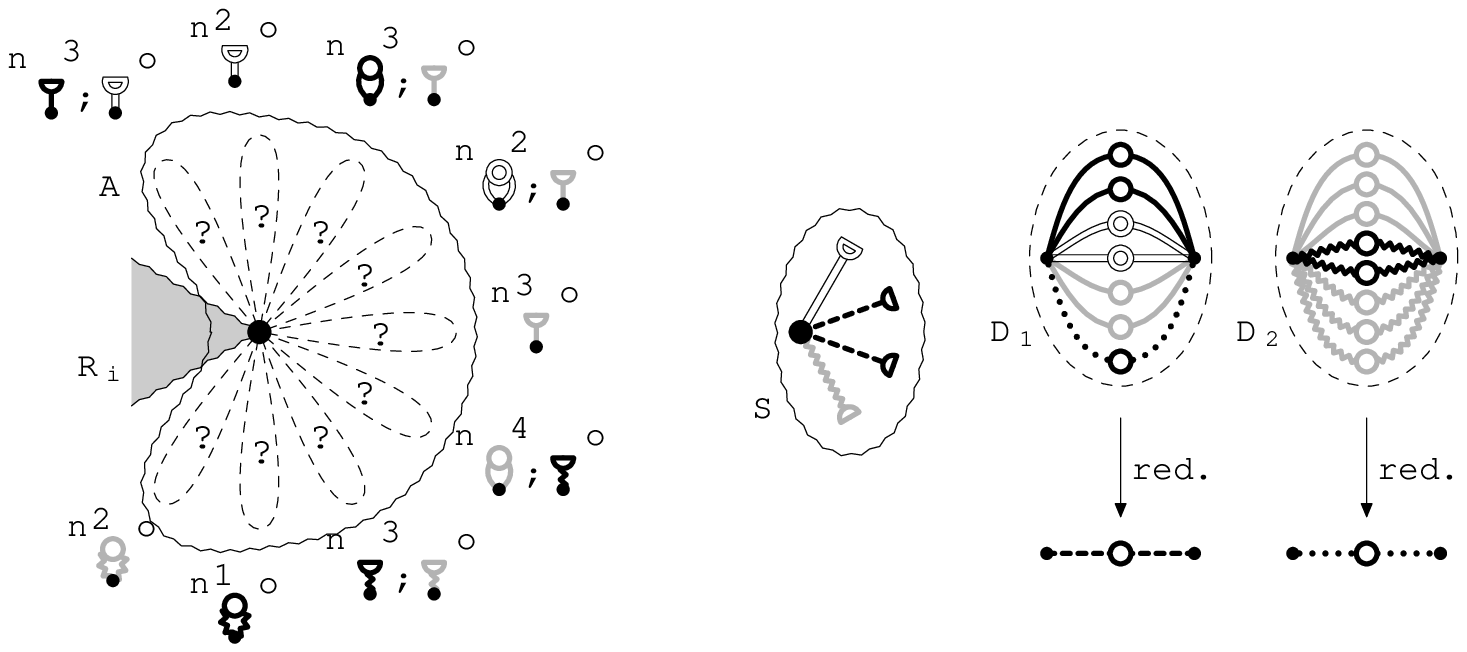}
\caption{On the left, a star atom $A$ in $\calR_i$ with 23 attached pendant elements, together with
their lists and multiplicities. On the right, a star atom $S$ from the catalog which belongs to
$\frakL(A)$. The bold dashed edge corresponds to the dipole $D_1$. The bold dotted edge in $D_1$
expands to another dipole $D_2$. The remaining colored edges correspond to proper atoms.}
\label{fig:testing_for_star_atoms}
\end{figure}

An instance of the problem is depicted in Figure~\ref{fig:testing_for_star_atoms}. Suppose that
$\calR_i$ contains a star atom $A$ with attached pendant elements, each with a previously computed
list. We want to determine the list $\frakL(A)$ which consists of some star atoms and half-quotients
of dipoles from the catalog. Let $S$ be a star atom from catalog, with attached half-edges, loops
and pendant edges. We want to test whether $S$ belongs to $\frakL(A)$. 

Before we do so, we establish basic properties of atoms concerning sizes. This properties allow us
to understand the structure.

\heading{Size Properties.} Let $A$ be an atom and let $Q$ be a quotient of this atom. We get the
following relations between the sizes of $Q$ and $A$, depending on the type of the quotient:
\begin{packed_itemize}
\item \emph{$Q$ is the edge-quotient:} Then $v(Q) = v(A)$ and $e(Q) = e(A)$.
\item \emph{$Q$ is the loop-quotient:} Then $v(Q) = v(A)-1$ and $e(Q) = e(A)$.
\item \emph{$Q$ is a half-quotient:} Then $v(Q) = v(A)/2$ and $e(Q) = e(A)/2$.
\end{packed_itemize}

Throughout each reduction, we calculate how many vertices and edges are in all the atoms replaced by
colored edges which we denote by $\hat v$ and $\hat e$. Initially, we put $\hat v(e) = 0$ and $\hat e(e) =
1$ for every edge $e \in E(G_0)$. For a subgraph $X$, we then define
$$
\hat v(X) := v(X) + \sum_{e \in E(X)} \hat v(e),
\qquad\text{and}\qquad
\hat e(X) := \sum_{e \in E(X)} \hat e(e).
$$
When an atom $A$ is replaced by an edge $e$ in the reduction, we put $\hat v(e) = \hat v(\int A)$
and $\hat e(e) = \hat e(\int A)$. So for a subgraph $X$ of $G_i$, the numbers $\hat v(X)$ and $\hat
e(X)$ are the numbers of vertices and edges when $X$ is fully expanded.

We similarly define $\hat v$ and $\hat e$ for quotients and their subgraphs; the difference here is
that the quotients might contain half-edges. Let $H(X)$ be the set of half-edges of a subgraph $X$.
Then for a half-edge $h \in H(X)$, created by halving an edge $e$, we put $\hat v(h) = \hat v(e)/2$
and $\hat e(h) = \hat e(e)/2$. For a subgraph $X$, we naturally define
$$
\hat v(X) := v(X) + \sum_{e \in E(X)} \hat v(e) + \sum_{h \in H(X)} \hat v(h),
\qquad\text{and}\qquad
\hat e(X) := \sum_{e \in E(X)} \hat e(e) + \sum_{h \in H(X)} \hat e(h).
$$
Finally, we also use $\hat v$ and $\hat e$ for pendant elements, created in the reduction series
$\calR_0,\dots,\calR_t$.

We get straightforwardly the following:

\begin{lemma} \label{lem:list_sizes}
For every pendant element $x$, the possible pendant edge, the possible loop and all half-edges of
the list $\frakL(x)$ have the same sizes $\hat v$ and $\hat e$ as $\hat v(x)$ and $\hat e(x)$
respectively.\qed
\end{lemma}

We apply this when $\frakL(x)$ is computed since we consider only quotients of the correct sizes.
Thus we can speedup the subroutine of Algorithm~\ref{alg:computing_list}. But for purpose of this
section, the following corollary is important:

\begin{corollary}
Let $x$ and $y$ be two pendant elements.
\begin{packed_enum}
\item[(i)] If $\frakL(x)$ and $\frakL(y)$ contain a half-edge in common, then $\hat v(x) = \hat
v(y)$ and $\hat e(x) = \hat e(y)$.
\item[(ii)] Let $\frakL(x)$ contains a loop of a color $c$ and a half-edge of a color $c'$. Then
$\frakL(y)$ cannot contain both the loop of the color $c'$ and the half-edge of the color $c$.
\end{packed_enum}
\end{corollary}

\begin{proof}
(i) Implied by Lemma~\ref{lem:list_sizes} directly.

(ii) Let $\frakL(x)$ contain a loop $e$ and $\frakL(y)$ contain a half-edge $h$ of the same color.
Then $\hat v(x) = \hat v(e)+1$ and $\hat v(y) = \hat v(e)/2+1$ for the vertices, and $\hat e(x) = \hat
e(e)$ and $\hat e(y) = \hat e(e)/2$ for the edges. Therefore $x$ is larger than $y$. For the same
reason, we deduce that $y$ is larger than $x$, which gives a contradiction.\qed
\end{proof}

The property (i) relates half-edges together. For pendant edges and loops, recall also
Lemma~\ref{lem:list_properties}.  The property (ii) states that there is a certain size hierarchy on
the pendant elements as we discuss below. We use this hierarchy to simplify the testing problem for
$A$ and $S$ as follows.

\begin{figure}[b!]
\centering
\includegraphics{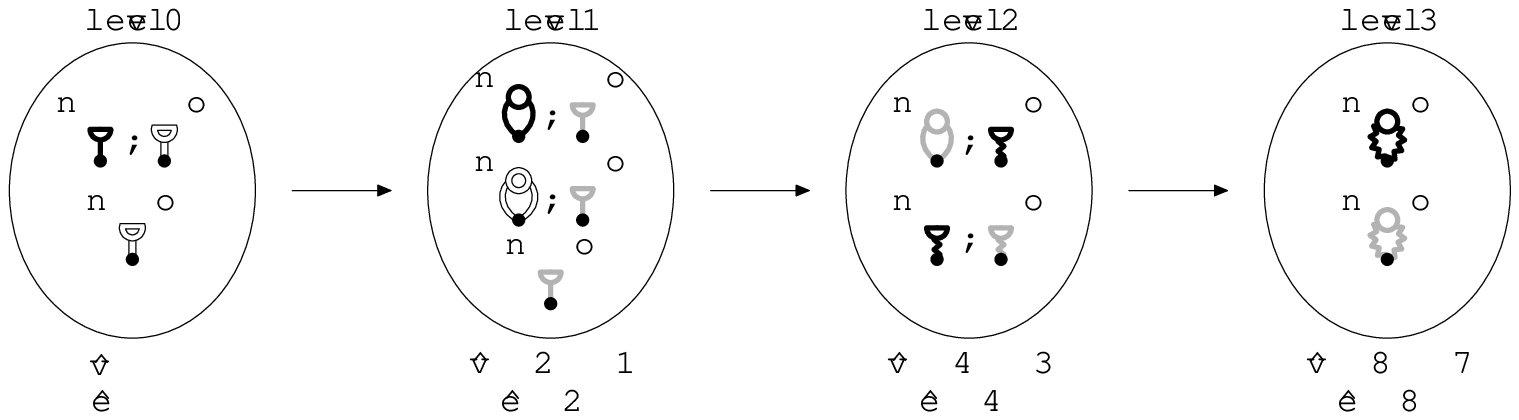}
\caption{A chain of pendant elements with four levels, obtained from the example in
Figure~\ref{fig:testing_for_star_atoms}, for some $\alpha$ and $\beta$. Notice that quotients
corresponding to one atom are placed in neighboring levels. All pendant elements of $A$ are placed
in this one chain, and we ignore their multiplicities.}
\label{fig:chains_of_pendant_elements}
\end{figure}

\heading{Chains of Pendant Elements.}
Suppose that we ignore pendant edges contained in the lists since they are easy to deal with.
Pendant elements are organized into independent chains, consisting of several levels. Each chain
starts with pendant elements $x$ of the level zero with $\hat v(x) = \alpha$ and $\hat e(x) =
\beta$. Further, it contains on the level $m>0$ the pendant elements $x$ with $\hat v(x) = 2^m
\alpha - (2^m-1)$ and $\hat e(x) = 2^m \beta$. See Figure~\ref{fig:chains_of_pendant_elements} for
an example of one chain from Figure~\ref{fig:testing_for_star_atoms}.

A star atom $A$ can contain multiple chains. It is important that different chains contain
completely different colors in their lists, so they behave completely independently. If $\frakL(x)$
contains a half-edge of a color $c$ and $\frakL(y)$ contains the loop of the same color $c$, then
$x$ belongs to a level $m$ and $y$ belongs to the level $m+1$ of the same chain.

\heading{Preprocessing Dipoles.}
We have a star atom $S$ with several pendant edges, loops and half-edges attached and we want to
test whether it embeds into a star atom $A$ in $\calR_i$. Let us denote the single vertex of both
$S$ and $A$ by $v$. As we discuss below, it is easy to deal with pendant edges and loops since they
correspond to unique edge- and loop-quotients. On the other hand, the half-edges are more complex
since we can have many different half-quotients of a half-edge. A half-edge can be of two types:
Either it corresponds to a half-quotient of a proper atom, or of a dipole; for example in
Figure~\ref{fig:testing_for_star_atoms} we have two half-edges corresponding to proper atoms, and
two half-edges corresponding to dipoles. In the case of a proper atom, it corresponds to exactly one
pendant element attached in $A$ to $v$. (Alternatively, one subtree of blocks attached to $v$ in
$H$.) In the case of a dipole, a half-quotient of this dipole may correspond to several different
pendants elements of $A$.

Recall Lemma~\ref{lem:dipole_quotients} describing the structure of every half-quotient of a dipole. The
resulting quotient has half-edges and loops attached to $v$. Again, edges of the dipole can be of
two types:
\begin{packed_itemize}
\item An edge can correspond to a proper atom, or alternatively it can be an original edge of $G$.
Then in a half-quotient we obtain a half-edge/loop from this edge which corresponds after the full
expansion to a subtree of blocks attached to $v$, and thus it corresponds to exactly one pendant
element of $A$.
\item Further, each dipole can also contain an edge corresponding to a dipole; for example in
Figure~\ref{fig:atoms_examples}, after one reduction step, we obtain on top a dipole consisting of
two parallel edges, one corresponding to a dipole and the other to a proper atom. A half-quotient of
this dipole would correspond to several pendant elements of $A$ which we want to avoid. According to
the definition of a dipole, there can be at most one edge corresponding to a dipole since every dipole
contains all parallel edges between the given two vertices. Therefore, we can just expand this edge
by replacing it with the edges of the dipole.  We obtain exactly the same half-quotients as before.
\end{packed_itemize}

\begin{figure}[b!]
\centering
\includegraphics{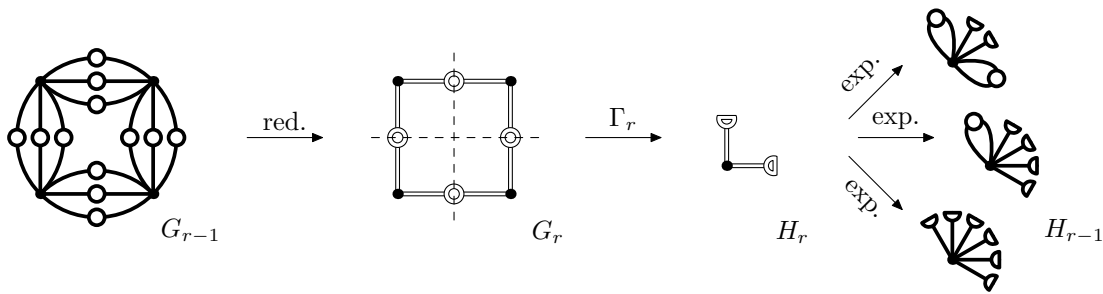}
\caption{By two reflections, the quotient $H_r$ consists of a star atom with two half-edges
corresponding to dipoles. All three expansions $H_{r-1}$ up to isomorphism are depicted. It is not
possible to expand a quotient $H_{r-1}$ with three attached loops.}
\label{fig:star_atom_with_two_dipoles}
\end{figure}

\begin{figure}[t!]
\centering
\includegraphics{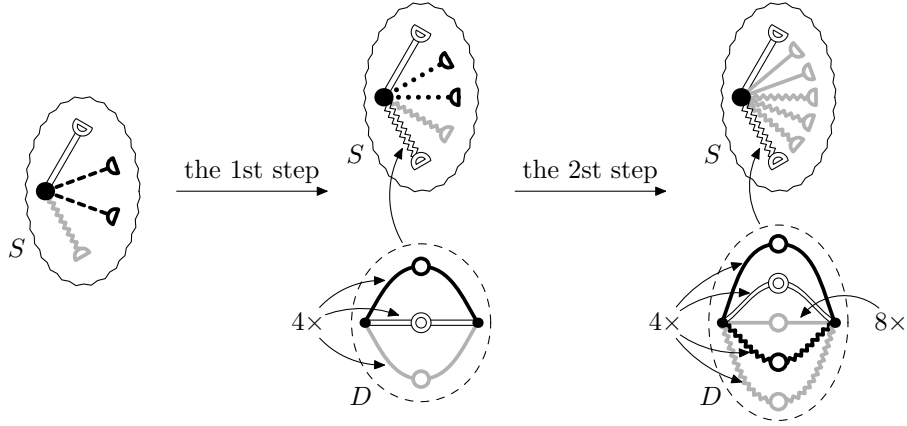}
\caption{In the first step, we expand both
dipoles $D_1$, each having two edges of each of the three color classes corresponding to proper
atoms. So the unified dipole $D$ now contains four colored edges for each of the three color classes.
Further $D_1$ has one halvable edge corresponding to a dipole $D_2$, so one half-edge of this
color is attached to $S$ directly for each $D_1$. In the second step, we expand the two half-edges corresponding to
$D_2$. Here we have two black curly edges, which are directly placed to $D$. But the remaining two color
classes have both odd sizes, so one edge from each is directly attached as a half-edge to $S$. The remaining
edges are placed to $D$. The resulting star atom $S$ together with the resulting unified dipole $D$
is depicted on the right.}
\label{fig:result_of_preprocessing}
\end{figure}

Further, $S$ may contain multiple half-edges corresponding to half-quotients of dipoles; this can be
obtained by factorization as depicted in Figure~\ref{fig:star_atom_with_two_dipoles}.
Suppose that each dipole is expanded as described, so it consists of some edges, each corresponding to a proper
atom or an original edge of $G$. We would like to unify these dipoles into one dipole $D$, containing
all the edges. But this might introduce additional quotients as in
Figure~\ref{fig:star_atom_with_two_dipoles}. If two dipoles both contain an odd number of edges of
one color, in the unified dipole we have a half-quotient consisting of only loops of this
color which is not possible in the case of two separated dipoles. There is an easy fix of the
problem, we pre-process each dipole and we remove one edge from each color class of odd size (of
necessarily halvable edges) and place it as the half-edge directly in $S$ attached to $v$. We surely
know that at least one half-edge of this color appears in every half-quotient of this dipole.
The resulting star atom $S$ belongs to $\frakL(A)$ if and only if the original $S$ belongs there.
In Figure~\ref{fig:result_of_preprocessing}, we illustrate this preprocessing for the example in
Figure~\ref{fig:testing_for_star_atoms}.

\heading{Dealing with Attached Edges and Loops.}
After the preprocessing step, the star atom $S$ contains several pendant edges, loops, half-edge
attached to $v$, with at most one half-edge corresponding to a dipole. If this dipole contains some
non-halvable edges, then they are paired in every half-quotient and form loops. So we can remove
them from the dipole and attach the corresponding number of loops directly to $v$ in $S$. After this
step, the dipole contains only even number of halvable edges in each color class.

Each pendant edge corresponds to a block atom attached to $v$. Each loop either corresponds to a
proper atom, or a dipole. Similarly as before, we expand these dipoles, replacing them by several
loops corresponding to their edges. After this expansion, each pendant edge/loop corresponds to
exactly one pendant element of $A$. But since their expansion is unique, they can be contained in
lists of only one type of pendant elements. Therefore we can arbitrarily assign pendant elements and
remove these pendant elements and loops from $S$ and the assigned pendant elements from $A$.

The resulting star atom $S$ contains only half-edges and at most one half-edge corresponds to a
dipole having only halvable edges corresponding to proper atoms. Hence the considered star atom
decomposes into the half-edges $H$ corresponding to proper atoms and at most one half-edge
corresponding to a halvable dipole $D$.

\heading{Reduction to the V-Matching Problem.} Suppose that $S$ contains a half-edge corresponding
to a dipole $D$, otherwise the half-quotient of $S$ is unique and we can easily match it $A$ using
perfect matching in a bipartite graph, as described in the proof of
Lemma~\ref{lem:lists_for_star_atoms}. Let $H$ be the set of the remaining half-edges attached in $S$
corresponding to proper atoms. In each half-quotient of the dipole $D$, we have color classes of
even sizes. For each color class of an even size $m$, we can choose an arbitrary number $\ell$ of
loops and the corresponding number of half-edges $m-2\ell$. If we know these values $m$ and $\ell$,
we can just test existence of a perfect matching as described in
Lemma~\ref{lem:lists_for_star_atoms}. Since we do not know the values $m$ and $\ell$, we need to
solve a generalization of perfect matching called \ivmatch\ in which we are free to choose these
values.

The input of \ivmatch\ gives a bipartite graph $B$ defined similarly as before. One part has a vertex
per a pendant element and the other part has a vertex for each edge of $D$ and half-edge of $H$. We
put an edge between $e \in D$ and $x$, if a half-quotient of $e$ or the loop-quotient created by
rotating two edges $e$ is in the list $\calL(x)$. We call the first case a \emph{half-incidence} and
the latter case a \emph{loop-incidence}. Further, we add a \emph{half-incidence} between $h \in H$
and $x$, if the half-edge $h$ is contained in the list $\calL(x)$. 

\begin{figure}[b!]
\centering
\includegraphics{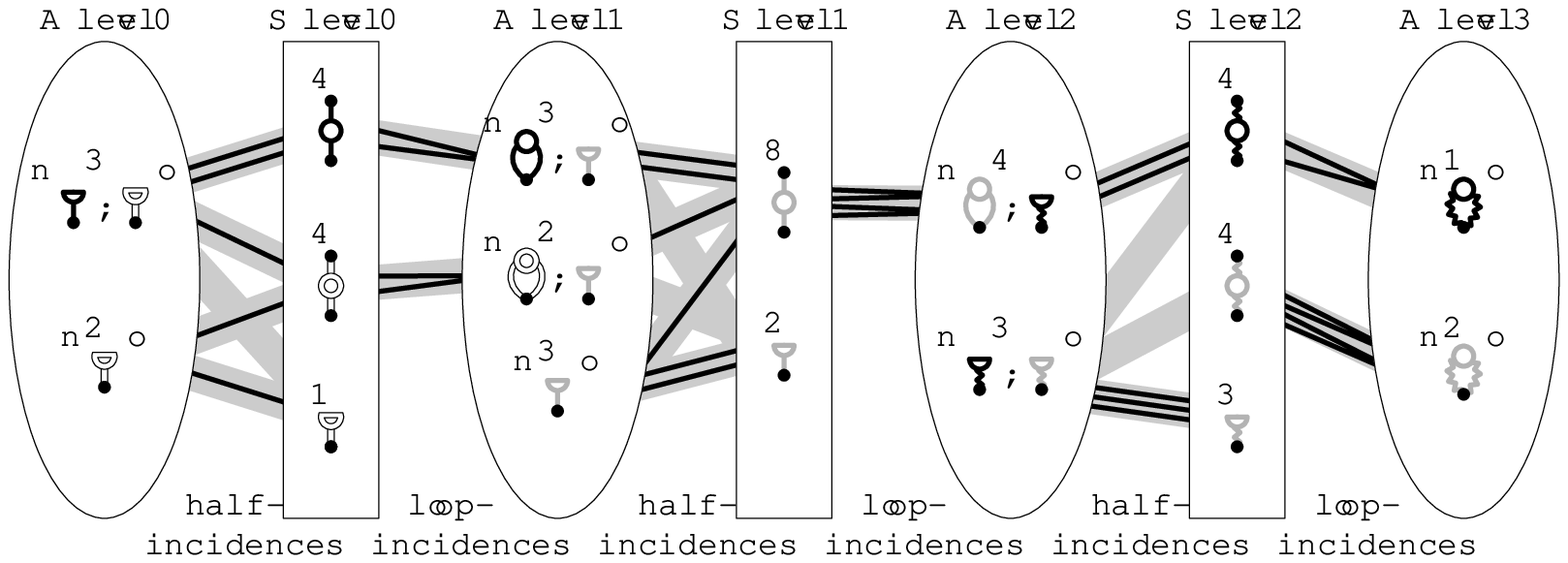}
\caption{An instance of the \ivmatch\ problem corresponding to the input $A$ in
Figure~\ref{fig:testing_for_star_atoms} and the preprocessed star atom $S$ in
Figure~\ref{fig:result_of_preprocessing}. The bipartite graph $B$ with edges depicted in gray and a
spanning subgraph $B'$ highlighted in bold. The part of pendant elements is in circles, the other
part is in boxes. The edges of $D$ are depicted by edges and the remaining half-edges corresponding
to proper atoms and depicted by half-edges.}
\label{fig:v-matching_problem_example}
\end{figure}

We ask whether there exists a \emph{spanning subgraphs} $B'$ of $B$, with each component of
connectivity a path of length one or two, as follows. Each pendant element $x$ is in $B'$ either
half-incident to exactly one vertex in the other part, or it is loop-incident to exactly two edges
$e,e' \in D$ of the same color class.  Further, each edge and half-edge is incident in $B$ to
exactly one pendant element. In what follows, we call $B'$ an \ivsubgraph\ of $B$.
See Figure~\ref{fig:v-matching_problem_example} for an example, with
several additional properties which we discuss now.

\heading{Special Properties of Inputs.}
The \ivmatch\ would likely be \cNP-complete in general, but there are additional properties already
depicted in Figure~\ref{fig:v-matching_problem_example}. We describe them in details, since they
might help in constructing a polynomial-time algorithm for this problem.

First, $B$ consists of separate connected components and for each we can solve the problem
separately. A connected component is induced by a chain of pendant elements and its incident
vertices representing $S$. In what follows, we assume that we just have a single chain.
This chain has its level structure which translates into a level structure for part with half-edges
of $H$ and edges of $D$ as follows.
\begin{packed_itemize}
\item Let $h$ be a half-edge of $H$ with $\hat v(h) = \alpha$ and $\hat e(h) = \beta$. Then this
half-edge can be in $B$ half-incident to pendant elements $x$ only of the level with $\hat v(x)=
\alpha$ and $\hat e(x) = \beta$.
\item Let $e$ be an edge of the dipole $D$ with $\hat v(e) = 2\alpha$ and $\hat e(e) = 2\beta$. It
can be half-incident to pendant elements $x$ only of the level with $\hat v(x) = \alpha$ and $\hat
e(x) = \beta$ and loop-incident to pendant elements $x$ only of the level with $\hat v(x) =
2\alpha-1$ and $\hat e(x) = \beta$.
\end{packed_itemize}
So we also have levels for the other part of $B$. Let us call these levels \emph{$A$ levels} and
\emph{$S$ levels}, respectively.  If we depict all levels from left to right according to their
order, alternating $A$ levels of pendant elements levels and $S$ levels of half-edges of $H$ and
edges of $D$. The graph $B$ has only edges going between consecutive levels as depicted in
Figure~\ref{fig:v-matching_problem_example}.

The second property says that $B$ can be viewed as a cluster graph. In each $S$ level, the edges of
$D$ and the half-edges of $H$ form clusters according to their color classes. Similarly in each $A$
level, the pendant elements form clusters according to equivalence classes of their lists. (We note
that two pendant elements with equal lists can correspond to non-isomorphic subgraphs in $H$. Then
their lists contain only half-edges.) Two clusters from different levels are either completely
adjacent, or not adjacent at all. In other words, the subgraph induced by all edges between two
clusters is either complete bipartite, or contains no edges.

Also each edge cluster is loop-adjacent to at most one pendant element cluster, since a loop-quotient
contained in a list uniquely determines the pendant element, as discussed above. We do not need the
condition that each list is loop-incident in the \ivsubgraph $B'$ to two edges of the same color class
since this is automatically achieved. On the other hand, there are no constraints for
half-incidences between clusters.

The third property reads as follows. If the highest level is an $A$ level, we can deal with this level
greedily. All pendant elements of this level has to be loop-incident in $B'$, so we can ignore
half-edges in their lists and assign them directly to the previous $S$ level.  For instance in
Figure~\ref{fig:v-matching_problem_example}, we know that we have two black shaken loops and one
gray shaken loops. Therefore we can remove this level together with the corresponding number of
edges of $D$. We can assume that the highest level is an $S$ level.  Similarly if the smallest level
is an $S$ level, we can match it greedily. If the smallest level is an $A$ level, this greedy
approach does not work. We know how many lists have to be realized by their half-quotients. But the
lists may contain multiple half-edge colors and therefore may be realized by different edges of $D$.
Therefore in the next $S$ level we may end up with different numbers of edges of $D$ which may
change solvability on higher levels.

Another property is that we know for each level of incidences, between two consecutive levels, how
many edges have to be placed in $B'$. We can just process the graph $B$ from left to right. On the
smallest $A$ level, we know how many half-incidences are in $B'$. Therefore we can compute the
remaining number of edges of $D$ which are loop-incident in $B'$. By processing in this way from
left to right, we get all the numbers.  Notice that we do not know how many half- and
loop-incidences are in $B'$ for each cluster, otherwise we could solve the problem directly by a
perfect matching in a modified graph.

Since the sizes of the levels are growing exponentially from left to right, we have at most
logarithmically many levels with respect to the size of the graph $H$. One can also show that the
number of vertices in each level is somewhat limited by the size of the level, but it does not seem
to give any useful bound.

\section{Applying the Meta-algorithm to Planar Graphs} \label{sec:planar_graphs}

In this section, we show that the meta-algorithm described in Section~\ref{sec:algorithm} applies to the class
of planar graphs. We describe the properties of planar graphs and their quotients in more details
than necessary. We do it to give a deeper insight into regular covers of planar graphs.

\subsection{Automorphism Groups of 3-connected Planar Graphs}

We review the well-known properties of planar graphs and their automorphism groups. These strong
properties are based on the Whitney's Theorem~\cite{whitney} stating that a 3-connected graph has a
unique embedding into the sphere.

\heading{Spherical Groups.}
A group is \emph{spherical} if it is a group of symmetries of a tiling of the sphere. The first class
of spherical groups are the subgroups of the automorphism groups of the platonic solids, i.e., $\gS_4$ for the
tetrahedron, $\gC_2 \times \gS_4$ for the cube and the octahedron, and $\gC_2 \times \gA_5$ for the
dodecahedron and the icosahedron. See Table~\ref{tab:number_of_subgroups} for the number of
conjugacy classes of the subgroups of these three groups. Note that conjugate subgroups $\Gamma$
determine isomorphic quotients $G / \Gamma$. The second class of spherical groups is formed by
infinite families $\gC_n$ and $\gD_n$.

\begin{table}[b!]
\centering
\begin{tabular}{|c|c|c|c|}
\hline
\multicolumn{4}{|c|}{$\gS_4$ of the order 24}\\
\hline
Order & Number & Order & Number\\
\hline
1 & 1 & 6 & 1\\
2 & 2 & 8 & 1\\
3 & 1 & 12 & 1\\
4 & 3 & & \\
\hline
\end{tabular}
\hskip 0.5cm
\begin{tabular}{|c|c|c|c|}
\hline
\multicolumn{4}{|c|}{$\gC_2 \times \gS_4$ of the order 48}\\
\hline
Order & Number & Order & Number\\
\hline
1 & 1 & 8 & 7\\
2 & 5 & 12 & 2\\
3 & 1 & 16 & 1\\
4 & 9 & 24 & 3\\
6 & 3 & & \\
\hline
\end{tabular}
\hskip 0.5cm
\begin{tabular}{|c|c|c|c|}
\hline
\multicolumn{4}{|c|}{$\gC_2 \times \gA_5$ of the order 120}\\
\hline
Order & Number & Order & Number\\
\hline
1 & 1 & 8 & 1\\
2 & 3 & 10 & 3\\
3 & 1 & 12 & 2\\
4 & 3 & 20 & 1\\
5 & 1 & 24 & 1\\
6 & 3 & 60 & 1\\
\hline
\end{tabular}
\caption{The number of conjugacy classes of the subgroups of the groups of platonic solids.}
\label{tab:number_of_subgroups}
\end{table}

\heading{Automorphisms of a Map.}
A map $\calM$ is a 2-cell embedding of a graph $G$ onto a surface $S$. For purpose of this paper,
$S$ is either the sphere or the projective plane.  A \emph{rotational scheme} at a vertex is a
cyclic ordering of the edges incident with the vertex.  When working with abstract maps, they are
graphs with a rotational scheme given for every vertex. An \emph{angle} is a triple $(v,e,e')$ where
$v$ is a vertex, and $e$ and $e'$ are two incident edges which are consecutive in the rotational
scheme of $v$.

An automorphism of a map is an automorphism of the graph which preserves the angles; in other words
the rotational scheme is preserved up to reflections. With the exception of paths and cycles,
$\Aut(\calM)$ is a subgroup of $\Aut(G)$. In general these two groups might be very different.
For instance, the star $S_n$ has $\Aut(S_n) = \gS_n$, but for any map $\calM$ of $S_n$ we just have
$\Aut(\calM) = \gD_n$.

\begin{lemma} \label{lem:map_groups}
For any map $\calM$, a permutation representation of the group $\Aut(\calM)$ can be computed in time
$\O(m^2)$ where $m$ is the number of edges of $\calM$.
\end{lemma}

\begin{proof}[Sketch]
There are $\O(m)$ angles in $\calM$. We fix one angle $(v,e,e')$, and test for each other angle
whether there is an automorphism mapping $(v,e,e')$ to it. The key observation is that if such an
automorphism exists, it is uniquely determined. We can just test in $\O(n+m)$ whether the map is
compatible with this prescribed mapping. The total running time is $\O(m^2)$.\qed
\end{proof}

If $\calM$ is drawn on the sphere, then $\Aut(\calM)$ is isomorphic to one of the spherical
groups~\cite{gross_tucker,coxeter}. In case $G$ is a $3$-connected planar graph, there exists the
unique embedding of $G$ onto the sphere. Then for any map $\calM$ of $G$, we have $\Aut(G) \cong
\Aut(\calM)$~\cite{whitney}.  Thus we get the following corollary of Lemma~\ref{lem:map_groups}.

\begin{corollary} \label{cor:aut_3_conn}
If $G$ is a $3$-connected planar graph, then $\Aut(G)$ is isomorphic to one of the spherical groups.
We can determine this group in polynomial time and find permutations which generate it in time
$\O(n^2)$.
\end{corollary}

We note that a linear-time algorithm for computing automorphism groups of planar graphs is
known~\cite{hopcroft_tarjan_planar_iso}.

\subsection{Automorphisms and Quotients of Primitive Graphs and Atoms}

To be able to apply the meta-theorem on planar graphs, we first investigate the automorphism groups
of atoms and primitive graphs, and their quotients.  Recall that a graph is essentially 3-connected
if it is a 3-connected graph with attached single pendant edges. 

\begin{lemma} \label{lem:planar_primitive_graph}
For a planar primitive graph $G$, the group $\Aut(G)$ is a spherical group and can be computed in
time $\O(n^2)$.
\end{lemma}

\begin{proof}
If $G$ is essentially 3-connected, this is directly implied by Lemma~\ref{lem:aut_ess_3-conn} and
Corrolary~\ref{cor:aut_3_conn}. If $G$ is $K_2$ or $C_n$ with attached single pendant edges, then
$\Aut(G)$ can be computed trivially as well.\qed
\end{proof}

For an atom $A$, Recall that $\Fix(A)$ is the point-wise stabilizer of $\bo A$ in $\Aut(A)$.
Further, by $\Aut_{\bo A}(A)$ we denote the set-wise stabilizer of $\bo A$ in $\Aut(A)$.

\begin{lemma} \label{lem:planar_atom_aut_groups}
We get the following automorphism groups for atoms:
\begin{packed_itemize}
\item For a star block atom, we get in general $\Aut_{\bo A}(A) = \Fix(A)$ which is a direct
product of symmetric groups. 
\item For a dipole, we get in general $\Fix(A)$ as a direct product of symmetric groups. For a
symmetric dipole, we have $\Aut_{\bo A}(A) = \Fix(A) \rtimes \gC_2$, for non-symmetric $\Aut_{\bo
A}(A) = \Fix(A)$.
\item For a proper atom $A$ of a planar graph, the group $\Aut_{\bo A}(A)$ is a subgroup of
$\gC_2^2$ and $\Fix(A)$ is a subgroup of $\gC_2$.
\item For a non-star block atom $A$ of a planar graph, $\Aut_{\bo A}(A) = \Fix(A)$ and it is a
subgroup of $\gD_n$.
\end{packed_itemize}
\end{lemma}

\begin{proof}
The edges of the same color class of a star block atom $A$ can be arbitrary permuted, so we get for
$\Aut_{\bo A}$ and $\Fix(A)$ the same group which is a direct product of symmetric groups.
For a dipole the situation is similar, just for a symmetric dipole we can permute the vertices in
$\bo A$, so we get a semidirect product with $\gC_2$.

Let $A$ be a proper atom with $\bo A = \{u,v\}$, and let $A^+$ be the essentially 3-connected graph
created by adding the edge $uv$. Since $\Aut_{\bo A}(A)$ preserves $\bo A$, we have $\Aut_{\bo A}(A)
= \Aut_{\bo A}(A^+)$, and $\Aut_{\bo A}(A^+)$ fixes in addition the edge $uv$. Because $A^+$ is
essentially 3-connected, then $\Aut_{\bo A}(A^+)$ corresponds to the stabilizer of $uv$ in
$\Aut(\calM)$ for the unique map $\calM$ of $A^+$. But such a stabilizer has to be a subgroup of
$\gC_2^2$.  For $\Fix(A)$ we further stabilize the vertices of $\bo A$, so it is a subgroup of
$\gC_2$.

Let $A$ be a non-star block atom. Then $\bo A = \{u\}$ is preserved, so we have one fixed vertex in
both $\Aut_{\bo A}(A)$ and $\Fix(A)$, so they are the same. Since $A$ is essentially 3-connected,
then $\Aut_{\bo A}(A)$ is a subgroup of $\gD_n$ where $n$ is the degree of $u$.\qed
\end{proof}

As a corollary of Proposition~\ref{prop:semidirect_product}, we can characterize automorphism groups
of planar graphs, similarly to Babai~\cite{babai1975automorphism,babai1996automorphism}.

\begin{corollary}
Each automorphism group of a planar graph is obtained from a spherical groups by repeated semidirect
products of direct products of groups of atoms from Lemma~\ref{lem:planar_atom_aut_groups}.
\end{corollary}

\begin{proof}
The primitive graph $G_r$ has a spherical automorphism group by
Lemma~\ref{lem:planar_primitive_graph}.  Recall $\Aut(G_i) \cong \Aut(G_{i+1}) \ltimes
\Ker(\Phi_i)$. The kernel $\Ker(\Phi_i)$ is a direct product of the groups $\Fix(A)$ for all atoms
$A$ in $G_i$, which are by Lemma~\ref{lem:planar_atom_aut_groups} subgroups of $\gC_2$ and $\gD_n$,
and direct products of symmetric groups.
\end{proof}

\heading{Geometry of Quotients.} We review in more details the geometry of quotients of planar graphs.
These are classical results from geometry of the sphere, and the reader can find the missing proofs and
details in~\cite{stillwell}.  We note that this precise understanding is not necessary for the
correctness of the meta-algorithm, but we believe it gives a deeper insight into the quotients of
planar graphs, in the direction of the Negami's Theorem~\cite{negami}.

We start by basic definitions from geometry. An automorphism of a 3-connected planar graph is called
\emph{orientation preserving}, if the respective map automorphism preserves the global orientation
of the surface. It is called \emph{orientation reversing} if it changes the global orientation of
the surface. A group of automorphisms of a 3-connected planar graph is either orientation preserving,
or it contains a subgroup of index two of orientation preserving automorphisms. (The reason is that
composition of two orientation reversing automorphisms is an orientation preserving automorphism.)

Let $\tau$ be an orientation reversing involution. The involution $\tau$ is called \emph{antipodal}
if it is a semiregular automorphism of a closed orientable surface $S$ such that
$S/\langle\tau\rangle$ is a non-orientable surface. Otherwise $\tau$ is called a
\emph{reflection}. A reflection of the sphere fixes a circle.

In particular, the half-quotient of the sphere by an antipodal involution is the projective plane
and the half-quotient by a reflection is a disk.  An orientation reversing involution of a
3-connected planar graph is called \emph{antipodal} if the respective map automorphism is antipodal
and it is called a \emph{reflection} if the respective map automorphism is a \emph{reflection}. A
reflection of a map on the sphere fixes always either an edge, or a vertex.

\begin{lemma}[\cite{stillwell}] \label{lem:quotients_of_planar_graphs}
Let $G$ be a 3-connected planar graph and $\Gamma$ be a semiregular subgroup of $\Aut(G)$.
Then the following can happen:
\begin{packed_head_enum}{(a)}
\item[(a)] The action of $\Gamma$ is orientation preserving and the quotient $G / \Gamma$ is planar,
\item[(b)] The action of $\Gamma$ is orientation reversing but does not contain an antipodal
involution. Then the quotient $G / \Gamma$ is planar and necessarily contains half-edges,
\item[(c)] The action of $\Gamma$ is orientation reversing and contains an antipodal involution.
Then $G / \Gamma$ is projective planar.\qed
\end{packed_head_enum}
\end{lemma}

As we show, if $G / \Gamma$ is projective non-planar, then necessarily also the quotient of
associated primitive graph $G_r / \Gamma_r$ is projective non-planar. The reason is that a
projection of an atom is always planar, and planarity is preserved by expansions.

\begin{lemma} \label{lem:planar_proper_half_quotients}
Let $A$ be a proper atom in planar graph and let $\bo A = \{u,v\}$.  Then there are at most two
involutory semiregular automorphisms $\tau \in \Aut_{\bo A}(A)$ transposing $u$ and $v$. Moreover, if there
are exactly two of them, then one is orientation preserving and the other one is a reflection. In
particular, the two respective quotient graphs may be non-isomorphic planar graphs.
\end{lemma}

\begin{proof}
The graph $A^+$ is an essentially 3-connected planar graph with a unique embedding into the plane.
Then any graph automorphism $\tau$ transposing $u$ and $v$ is a map automorphism fixing $e$. It is
easy to see that that either $\tau$ is a 180 degree rotation around the centre of $e$, or it is a
reflection~\cite{babai1975automorphism,babai1996automorphism}. According to
Lemma~\ref{lem:quotients_of_planar_graphs}, both possible quotients are planar.\qed
\end{proof}

The following lemma is straightforward, and completes the description of possible quotients of
primitive graphs.

\begin{lemma}
Let $H$ be a quotient of a cycle $C$ by a semiregular group of automorphisms. Then either $H$ is a
cycle, or a path, or a path with one pendant half edge, or $H$ is a path with two pendant half
edges. Depending on parity, just three of the above for cases happen.\qed 
\end{lemma}

\subsection{The Class of Planar Graphs Satisfies (P0) to (P3)}

In this section, we prove that planar graphs satisfy the properties (P0) to (P3), and thus the
meta-algorithm of Section~\ref{sec:algorithm} applies to them. The class of planar graphs clearly
satisfies (P0). The graph isomorphism can be tested in polynomial time for planar
graphs~\cite{hopcroft_tarjan_planar_iso} and for projective planar graphs~\cite{pp_iso}, which
implies (P1). For (P2), Corollary~\ref{cor:aut_3_conn} states that the automorphism group of a
3-connected planar graph is a spherical group, and so we can generate all at most linearly many
subgroups of given order and check which ones act semiregularly. It remains to prove (P3):

\begin{proposition}
The class of planar graphs satisfies the property (P3).
\end{proposition}

\begin{proof}
Let $H$ be a 3-connected projective planar graph. Lichtenstein~\cite{pp_iso} proved that there are
two possible cases: Either the number of possible projective maps of $H$ is
bounded by $c \cdot v(H)$, or $H$ contains an edge $uv$ such that $H \setminus \{u,v\}$
is a planar graph. The property (P3) gives two 3-connected graphs $G$ and $H$, and if $G \hookto H$,
then necessarily $G \cong H$ (while ignoring the colors), and thus the same case happens for both
$G$ and $H$. We deal with these two cases separately.

\emph{Case 1: Linear number of projective maps.} We want to show that the size of $\Aut(H)$ is
polynomial in $v(H)$. To see this, every $\pi \in \Aut(H)$ is either a map automorphism, or an
isomorphism between two maps of $G$. The number of automorphisms of a map is bounded by $4e(H)$.
Since we have a linear number of non-isomorphic maps, we get the bound $4c\cdot e(H)\cdot v(H)$,
which is $\O(v^2(H))$. We can construct all these automorphisms in polynomial time, for instance by
attaching specific gadgets to the vertices of $G$ and $H$, and testing whether the modified graphs
are still isomorphic. Therefore we can compute $\Aut(H)$ in polynomial time, and by composition with
some isomorphism from $G$ to $H$, we can test whether there exists an embedding $G \hookto H$.

\emph{Case 2: Removal of two vertices makes each of the graphs $G$ and $H$ planar.} We know that
there exists $uv \in E(G)$ and $u'v' \in E(H)$ such that $G \setminus \{u,v\}$ and $H \setminus
\{u',v'\}$ are both planar graphs. We test $4e^2(G)$ possible choices of the pairs of vertices, and
we want to test whether there exists an embedding $\pi$ such that $\pi(u) = u'$ and $\pi(v) = v'$.
Clearly, $c(u) \in \frakL(u')$ and $c(v) \in \frakL(v')$. It remains to deal with the planar
remainders $G' = G \setminus \{u,v\}$ and $H' = H \setminus \{u,v\}$. Notice that $G'$ and $H'$ are
not necessarily 3-connected.

Now, we can approach the problem exactly as with testing expandability of $H_r$ in the
meta-algorithm. We code the colors of the vertices of $G'$ by the colors of the single pendant edges
attached to the corresponding vertices of $G'$, and we code the lists of the colors by the list of the
pendant elements attached to the vertices in $H'$. We choose the central block/articulation of
both $G'$ and $H'$ as the core.

First, we proceed with a series of reductions $G'_0,\dots,G'_r$ of $G'$, replacing the atoms by
edges and placing them in the separate catalog. Here, the automorphisms are indeed color
preserving. Further, we have four partition classes of $V(G')$ according to the adjacencies to $u$
and $v$, and automorphisms preserves them. We end up with a primitive 3-connected graph $G'_r$.

Using the constructed catalog, we proceed similarly with the reductions $\calR'_0,\dots,\calR'_r$ of
$H'$. But there are two important differences. First, the starting graph $\calR'_0$ already contains
pendant elements with initial lists, given by $H'$. The second difference is that we have to
remember lists also for interior edges of $\calR'_i$, not only for pendant elements. The reduction
is done similarly, by replacing the atoms of $\calR'_i$ by edges and pendant elements in
$\calR'_{i+1}$. Let $A$ be an atom in $\calR'_i$, and to compute the list we iterate over all atoms
of $G'$ of the given type in the catalog. Let $B$ be one such atom from the catalog. Depending on
the type of $A$, we proceed as follows:

\begin{packed_itemize}
\item If $A$ is a block atom of $\calR'_i$, then according to Lemma~\ref{lem:planar_atom_aut_groups}
we have $\Aut_{\bo A}(A)$ is a subgroup of $\gD_n$, and so we can test all possible isomorphisms from $B$ to
$A$. If one such isomorphism defines an embedding of $B$ to $A$, we add the edge representing $B$ to
the list.
\item If $A$ is a proper atom of $\calR'_i$, then we similarly test according to
Lemma~\ref{lem:planar_atom_aut_groups} at most four possible isomorphisms from $B$ to $A$.
\item If $A$ is a dipole, we construct a bipartite graph. One part is formed by the colored edges of
$B$, the other part by the edges of $A$ with the associated lists of possible colors. The
adjacencies are given by containments of the colors in the lists. We test existence of a perfect
matching in this bipartite graph.
\end{packed_itemize}

Finally, when we reach the primitive graph $\calR'_r$, then according to
Lemma~\ref{lem:planar_primitive_graph} we have that $\Aut(\calR'_r)$ is a spherical group. Thus we
can test all possible isomorphisms from $G'_r$ to $\calR'_r$ whether they are compatible with the
lists of $\calR'_r$. The entire subroutine is clearly correct and runs in polynomial time.\qed
\end{proof}

\begin{proof}[Corollary~\ref{cor:planar_rcover}]
We just apply Theorem~\ref{thm:metaalgorithm} since planar graphs satisfy (P0) to (P3).\qed
\end{proof}

\section{Concluding Remarks} \label{sec:conclusions}

We conclude this paper by remarks to the meta-algorithm and open problems.

\heading{Possible Extensions of The Meta-algorithm.}
There are several possible natural extensions of the meta-algorithm. First, we can easily generalize
it for input graph $G$ and $H$ with half-edges, directed edges and halvable edges, and also for
colored graphs. Further, one can prescribe lists of possible images for the vertices $V(G)$ and of
possible pre-images of the vertices $V(H)$, the expandability testing can compute also with these
lists.

The most interesting extension would be to use the meta-algorithm to generate all quotients $H$ of
$G$.  Indeed, this cannot be achieved in polynomial time since there might be exponentially many
such quotients. But we can enumerate all labeled quotients with a polynomial delay, by using different
expansions determined in Proposition~\ref{prop:quotient_expansion}.

\heading{Open problems.} We already discussed that \rcover\ generalizes the graph isomorphism
problem. We ask the following:

\begin{problem}
Is the problem \rcover\ \cGI-complete?
\end{problem}

\noindent Similarly to the graph isomorphism problem, we are given two graphs $G$ and $H$ which
restrict each other. To solve \rcover, one needs to understand automorphism groups of graphs and
their semiregular subgroups. 

As possible next direction of research, we suggest to attack classes of graphs close to planar
graphs, for instance projective planar graphs or toroidal graphs. To do so, it seems that new
techniques need to be built. Also the automorphism groups of projective planar graphs and toroidal
graphs are not well understood.

Lastly, we indeed ask whether the slow subroutine for dipole expansion of the meta-algorithm can be
solved in polynomial time:

\begin{problem}
Can the complexity of the algorithm of Theorem~\ref{thm:metaalgorithm} be improved to be polynomial?
\end{problem}

In Section~\ref{sec:comb_int_star_atoms}, we have shown that the slow subroutine reduces to a generalized
matching problem called \ivmatch. Here, we describe a purely combinatorial formulation of the
\ivmatch\ problem. This reformulation can be useful to people studying combinatorial optimization
since they can attack the problem without understanding the regular coverings and the structural
results obtained in the entire paper.

The input of \ivmatch\ consists of a bipartite graph $B=(V,E)$ with a partitioning
$V_1,\dots,V_\ell$, with $\ell$ even, of its vertices which we call \emph{levels}, with all edges
between of consecutive levels $V_i$ and $V_{i+1}$, for $i=1,\dots,\ell-1$. The levels
$V_1,V_3,\dots$ are called \emph{odd} and the levels $V_2,V_4,\dots$ \emph{even}. Further each level
$V_i$ is partitioned into several \emph{clusters}, each consisting of a few vertices with identical
neighborhoods. There are three key properties:
\begin{packed_itemize}
\item The incidences in $B$ respect the clusters; between any two clusters the graph $B$ induces
either a complete bipartite graph, or an edge-less graph.
\item Each cluster of an even level $V_{2t}$ is incident with at most one cluster at the odd level
$V_{2t+1}$, and each cluster of $V_{2t+1}$ is incident with at most one cluster at $V_{2t}$; so we
have a matching between clusters at levels $V_{2t}$ and $V_{2t+1}$.
\item The incidences between the clusters of $V_{2t-1}$ and $V_{2t}$ can be arbitrary.
\end{packed_itemize}

\begin{figure}[t!]
\centering
\includegraphics{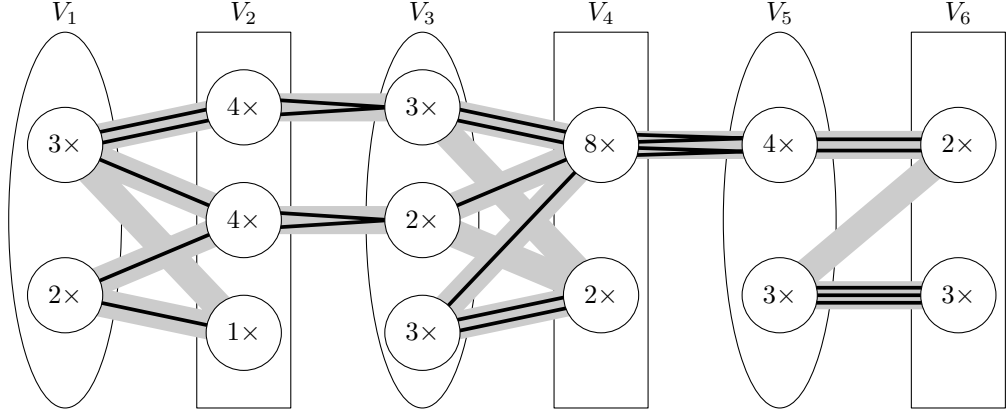}
\caption{An example input $B$, the clusters are depicted by circles together with their sizes.
The odd levels are drawn in circles and the even ones in rectangles. The edges of $B$ are depicted
by gray lines between clusters representing complete bipartite graphs. One spanning subgraph $B'$
solving the \ivmatch\ problem is depicted in bold.}
\label{fig:example_comb_v-matching}
\end{figure}

The problem \ivmatch\ asks whether there is a \emph{spanning subgraph} $B' = (V,E')$ with each
component of connectivity equal to a path of length one or two. We call this subgraph $B'$ an
\ivsubgraph\ of $B$. Each vertex of an odd level $V_{2t+1}$ is in $B'$ adjacent either to exactly
one vertex of $V_{2t+2}$, or to exactly two vertices of $V_{2t}$. Each vertex of an even level
$V_{2t}$ is adjacent to exactly one vertex of the levels $V_{2t-1} \cup V_{2t+1}$. In other words,
from $V_{2t-1}$ to $V_{2t}$ the edges of $E'$ form a matching, not necessarily perfect.  From
$V_{2t}$ to $V_{2t+1}$, the edges of $E'$ form independent $\vee$-shapes, with their centers in the
level $V_{2t+1}$. Figure~\ref{fig:example_comb_v-matching} shows an example.

\begin{problem}
What is the complexity of the \ivmatch\ problem?
\end{problem}

It would be interesting to know solutions even to more restricted versions of this problem. For
instance, what if the number of levels $\ell$ is small? We note that for instances having odd
numbers of levels, as the one in Figure~\ref{fig:example_comb_v-matching}, we can match the highest level
greedily; thus we can assume that the number of levels is always even. For two levels, the problem
is just the standard perfect matching problem for bipartite graphs. On the other hand, already the
case of four levels is very interesting and open. The instances of our algorithm needs to solve for the
problem \rcover\ have at most logarithmic number of levels in size of the smaller input graph $H$.

\bibliographystyle{plain}
\bibliography{algorithmic_regular_covers}

\end{document}